\documentclass[12 pt]{article}
\usepackage{listings}
\usepackage{amsmath}
\usepackage{amssymb}
\usepackage{graphicx}
\usepackage{longtable}
\usepackage{chapterbib}
\usepackage[timestamp,first]{draftcopy}
\usepackage[first]{draftcopy}
\draftcopySetScale{100}
\draftcopyPageX{200}
\draftcopyPageY{100}
\draftcopyBottomX{200}
\draftcopyBottomY{100}
\draftcopySetScaleFactor{6.0}
\usepackage{changebar}
\begin{document}
\title{\bf Developments from Programming the Partition Method for a Power Series Expansion}
\author{\bf Victor Kowalenko\\ ARC Centre of Excellence for Mathematics \\
and Statistics of Complex Systems\\ Department of Mathematics and Statistics\\
The University of Melbourne\\ Victoria 3010, Australia.}
\maketitle

\begin{abstract}
In a recent series of papers \cite{kow10}-\cite{kow11a} a novel method based on the coding of integer partitions 
has been used to derive power series expansions to previously intractable problems, where the standard Taylor/Maclaurin 
method fails. In this method the coefficient at each order $k$ of the resulting power series expansion is determined
by summing all the specific contributions made by each partition whose parts/elements sum to $k$. The specific
contributions are evaluated by assigning values to each element in a partition and then multiplying by a multinomial
factor, which depends on the frequencies of the elements in the partitions. This work aims to present for the first time
the theoretical framework behind the method, which is now known as the partition method for a power series expansion. 
To overcome the complexity in evaluating all the contributions from the partitions as the order $k$ increases, a 
programming methodology is created, thereby allowing for far more general problems to be considered than originally 
envisaged. This methodology is based on an algorithm called the bi-variate recursive central partition (BRCP) algorithm, 
which is, in turn, developed from a novel non-binary tree-diagram approach to scanning the integer partitions summing to a 
specific value. The main advantage of the BRCP algorithm over other means of generating partitions lies in the fact 
that the partitions are generated naturally in the multiplicity representation. By developing the theoretical 
framework for the partition method for a power series expansion it becomes apparent that scanning over all partitions 
summing to a particular value can be regarded as a discrete operation denoted here by the discrete operator $L_{P,k}[ \cdot]$. 
The summand inside this operator depends on the coefficients of an inner and outer series resulting from expressing the original 
function as a pseudo-composite function of two power series expansions. As a consequence, simple modifications to the program 
for the partition operator result in programs for other operators involving specific types of partitions such as
those with: (1) only odd or even elements, (2) a fixed number of elements, (3) discrete elements, (4) specific elements 
and (5) those with restrictions on the size of their elements. Another interesting modification results in the generation of
all the conjugate partitions for their original partitions by transposing the rows and columns of their Ferrers diagrams.
The operator approach is then applied to the theory of integer partitions, in particular to generalisations of
the generating functions for both discrete and standard partitions. The main generalisation involves the introduction of 
the parameter $\omega$, whose powers indicate the total number of elements in the partitions, while the coefficients
of the power series expansions become polynomials in $\omega$. Finally, power series expansions for more advanced 
infinite products involving quotients and products of the discrete and standard partition generating functions are derived,
culminating in the multi-parameter infinite product first studied by Heine.    

\vspace{1.0 cm}
{\em Keywords: Absolute convergence; Algorithm; Bivariate recursive central partition algorithm (BRCP); Coefficient; Conjugate 
partition; Conditional convergence; Discrete partition; Discrete partition number; Discrete partition polynomial; 
Divergent series; Divisor polynomials; Doub\-ly-restri\-cted partition; Equivalence; Even partition; Generating function; 
Infinite product; Multinomial factor; Multiplicity representation; Odd partition; Partition; Partition function; 
Partition function polynomial; Partition method for a power series expansion; Partition operator; Partition polynomial; 
Power series expansion; Programming methodology; Pseudo-composite function; Recurrence relation; Regularised value; 
Taylor/Maclaurin series; Transpose; Tree diagram}

\vspace{1.0cm}
\end{abstract}
\vspace{5cm}
\lstset{language=C, showstringspaces=false, showspaces=false}

\newpage
\section{Introduction}
This work grew out of a desire to develop a programming methodology on the partition method for a power series
expansion, which has featured prominently in a recent series of papers \cite{kow10}-\cite{kow11a} aimed at deriving 
power series expansions for previously intractable mathematical functions. Although still a relatively 
novel method, the method for a power series expansion was first introduced in the derivation of an asymptotic 
expansion for the particular Kummer or confluent hypergeometric function that emerges in the response theory of 
magnetised quantum plasmas such as the degenerate electron and charged Bose gases \cite{kow94}. Specifically, a large 
$|\alpha|$-expansion was obtained for $_1F_1(\alpha,\alpha+1;z)$, which is itself a variant of the incomplete gamma 
function $\gamma(\alpha,z)$. As a consequence, the physical properties of the magnetised charged Bose gas were later 
determined in the weak field limit in Ref.\ \cite{kow99}, the first time ever that the limit had been investigated for 
a plasma. Since then, the method for a power series expansion has been applied successfully to various mathematical 
functions such as $1/\ln^s(1+z)$, $\sec^s z$, $z^2 \csc^s z$, and the three Legendre-Jacobi elliptic integrals, 
$F(\psi,k)$, $E(\psi,k)$ and $\Pi(\psi,n,k)$, while more recently, it has been applied to a finite sum of inverse 
powers of cosines \cite{kow11a}, resulting in the development of a spectacular empirical method to solve the problem.

Initially, it was observed that in order to apply the partition method for a power series expansion the original function 
needed to be expressed as a composite function. In the present work, however, this condition will be relaxed to quotients 
of ``pseudo-composite" functions. The two functions involved in the construction of the quotients of the pseudo-composite functions 
must themselves be expressible in terms of power series expansions, which are referred to as the inner and outer series. Neither 
of these series, however, is required to be absolutely convergent. When the conditions for the quotients are met, a power series 
expansion can be obtained in which the coefficients at each order $k$ are calculated by summing the specific contributions 
due to each partition in which the parts or elements sum to $k$. The contribution made by each partition is in turn determined 
by: (1) assigning a value to each element in the partition, (2) multiplying these values by a multinomial factor composed of 
the factorial of the total number of elements, $N_k!$ divided by the factorial of the number of occurrences or frequencies of 
each element $i$ in the partition, $n_i!$, and (3) if necessary, carrying out a further multiplication with the coefficients 
of the inner series set at the number of elements in the partition.
 
As explained in the introduction to Ref.\ \cite{kow11}, the partition method for a power series expansion is able to produce 
power series where the standard technique or Taylor/Maclaurin series approach breaks down, but for those cases where
a power series expansion can be obtained by the standard technique, the power series expansions are identical, although from
a totally different perspective. As a consequence, the cross-fertilisation of both approaches is frequently responsible for 
new mathematical results and properties. A particularly fascinating property of the partition method for a power series 
expansion is that the discrete mathematics of partitions is being employed to derive power series expansions for continuous 
functions. 

Whilst Ref.\ \cite{kow11} was primarily concerned with the application of the partition method for a power series 
expansion to basic trigonometric and related functions, it was stated that a theoretical framework was being developed 
for situations involving more complicated functions where the assignment of values to the elements of the partitions is 
no longer specific, but quite general. This theoretical framework was necessary for the the development of a programming 
methodology whose purpose was to facilitate the summation process over the partitions. Therefore, both the theoretical framework 
and development of a programming methodology represent important topics in the present work. However, in the course of developing 
the theoretical framework it became apparent that the theory of partitions was also affected. Moreover, with the development 
of the new algorithm to facilitate the calculation of the coefficients in the partition method for a power series expansion, 
one could tackle various problems pertaining to integer partitions such as the evaluation of: (1) doubly-restricted partitions, 
(2) partitions with a fixed number of elements, (3) conjugate partitions via Ferrers diagrams, and (4) discrete/distinct 
partitions. As we shall see, such problems, which often require implementing different algorithms if they can be solved at 
all, can be addressed by making minor adjustments or modifications to the algorithm presented in Sec.\ 2 of this work. 
Therefore, by developing the theoretical framework and tackling the mathematical programming issues associated with the 
method, we have not only been to improve the method by making it more general than first thought, but we have also been 
able to solve other outstanding problems in the theory of partitions. In addition, the new operator approach that evolved 
in the process casts partitions in a different light and it is hoped that this approach will result in new advances in the 
future. Consequently, the work has progressed well and truly beyond its original conception resulting in its present 
size. 

Throughout this work the concept of regularisation of a divergent series will be employed. This concept is defined 
here as the removal of the infinity in the remainder so as to make the series summable or yield a finite limit. The 
finite limits obtained in this process are referred to as regularised values, whilst the statements in which they 
appear together with the series can no longer be regarded as equations. Instead, they are referred to as equivalence 
statements or equivalences, for short. A necessary property of the regularised value is that it must be unique, 
particularly if it is identical to the value one obtains when the series is absolutely convergent within a finite radius 
of the complex plane such as the geometric series. Much of this will become clearer as we progress further in the 
present work, but for those readers wishing to seek a greater understanding of the concept and its ramifications to 
asymptotics, they should consult Refs.\ \cite{kow11} and \cite{kow09a}-\cite{kow11b}.

Since all the partitions summing to a particular value are required to determine the coefficients in the partition
method for a power series expansion, Sec.\ 2 examines the current state of the art for generating partitions. Here,
various algorithms are compared with the bi-variate recursive central partition algorithm or BRCP for short, which 
has only been sketched out in Refs.\ \cite{kow10} and \cite{kow11}. The BRCP algorithm is based on a non-binary 
tree diagram representation for all the partitions summing to a specific value. Whilst it is acknowledged that some of 
the alternative algorithms can be faster than the BRCP algorithm in scanning the partitions and occasionally, when 
printing them out on a screen, the BRCP algorithm is the most efficient algorithm for implementation in the 
partition method for a power series expansion. This is because it generates the partitions in the multiplicity 
representation, which, as we shall see, is both the ideal and minimum amount of information required by the partition 
method for a power series expansion.

Sec.\ 3 presents the general mathematical theory underpinning the partition method for a power series expansion, which 
is discussed in terms of the quotient of two pseudo-composite functions in Theorem\ 1. Under certain conditions it
is also shown that the method can be used to derive a power series expansion for the inverse or reciprocal of the 
quotient. Consequently, a hybrid recurrence relation involving the coefficients from both power series expansions is 
obtained, while an example involving the reciprocal of a Bessel function to arbitrary order, viz.\ $1/J_{\nu}(z)$,
is presented to make the preceding material clearer to the reader. Later in the section, Theorem\ 1 is extended in a 
corollary, where the quotient of the pseudo-composite functions is taken to an arbitrary power, $\rho$. This extension 
only affects the multinomial factor in the method for a power series expansion by transforming it into the Pochhammer 
symbol involving $\rho$ and the number of elements in each partition.

Because the partition method for a power series expansion is based on summing over all the partitions summing to 
a specific value $k$, a discrete operator called the partition operator or $L_{P,k}[\cdot]$ is introduced. When this 
operator acts on unity, it gives the number of partitions or partition function, $p(k)$. On the other hand, if the 
operator acts upon only the multinomial factor mentioned above in Step\ 2, then it is found to yield a value of 
$2^{k-1}$, while if the phase factor of $(-1)^{N_k}$, where $N_k$ represents the number of elements, is included 
with the multinomial factor, then it vanishes. This means that the expressions for the coefficients given in Theorem\ 
1 can be expressed in terms of the new operator. One interesting property of the partition operator is that frequently 
the arguments inside the operator can be interchanged with the result the operator yields. As a consequence, not only 
are many of the results derived previously via the partition method for a power series expansion in Refs.\ 
\cite{kow10}-\cite{kow11} expressed in terms of the operator, but also the inverse relations are presented. The section 
concludes by showing that if the quotient of the pseudo-composite functions yields an infinitely differentiable function, 
then its derivatives can also be expressed in terms of the partition operator.  

Sec.\ 4 is devoted towards the creation of the programming methodology that enables the coefficients obtained from the 
partition method for a power series expansion to be calculated from the general theory in Sec.\ 3. This is necessary 
because beyond the first few orders, it becomes increasingly onerous to determine the coefficients by hand. The 
computational task is divided into two steps. In the first step a code utilising the BRCP algorithm of Sec.\ 2 is 
written in C/C++ so that its symbolic output can be processed by Mathematica \cite{wol92} in the second step. The 
reason for the second step is that often the coefficients in the resulting power series expansions are either rational 
or algebraic in nature and that neither of these forms can be handled effectively in C/C++ with its floating point 
arithmetic. In fact, the coefficients often become so small that they would practically vanish in C/C++. By importing the 
output into Mathematica, we can exploit its integer arithmetic routines to express the coefficients in integer form or
we can use the symbolic routines to express them, for example, as polynomials when required. 

The material in Secs.\ 2 to 4 serves as a platform for studying various problems in theory of partitions. First, we consider 
the issue of generating specific types/classes of partitions or different operators such as: (1) those with a fixed number 
of elements in them, (2) doubly restricted partitions where all the elements are greater one value and less than another, 
(3) discrete or distinct partitions where an element occurs only once in a partition and (4) partitions with specific elements 
in them. Solving these problems invariably means developing new algorithms or codes, but in the case of the BRCP algorithm we 
shall see that they can be solved with relatively minor modifications. This is due to the power and versatility of the 
of the tree diagram approach upon which the BRCP algorithm is based. Furthermore, the BRCP code can be adapted to exploit 
two-dimensional dynamic memory allocation when one wishes to determine conjugate partitions by means of Ferrers diagrams. 
For the benefit of the reader many of the codes discussed in Secs.\ 4 and 5 are presented in the appendix, where it can 
be seen that they are surprisingly compact.

In Secs.\ 6-8 the partition method for a power series expansion and the operator approach are employed in the derivation of 
generating functions from increasingly sophisticated extensions of the infinite product defined by $P(z)= \prod_{k=1}^{\infty} 
(1-z^k)$ and its inverse or $1/P(z)$. This famous product was found by Euler to yield the generating function whose 
coefficients are equal to the partition function $p(k)$. Before the study commences, however, Sec.\ 6 begins by 
showing in Theorem\ 2 that the generating function of $P(z)$ is absolutely convergent within the unit disk centred at the 
origin in the complex plane, but is divergent for all other values of $z$. That is, $P(z)$ represents the regularised value 
for divergent values of $z$ as given by Equivalence\ (\ref{fortyfoura}).  This is often postulated in the literature, but no 
formal proof of this important result has ever appeared previously. Next Theorem\ 1 is used to derive the generating function 
for $1/P(z)$ whose coefficients are given in terms of the partition operator acting with each element $i$ assigned a value 
of $p(i)$. In this case the coefficients $q(k)$, which are referred to as the discrete partition numbers, are only non-zero 
when $k$ is a pentagonal number, again a result that was first obtained by Euler. By inverting this method we obtain the 
partition function $p(k)$ in terms of the partition operator acting with each element $i$ assigned to $q(i)$. 

Because most of the discrete partition numbers vanish, the new result for the partition function simplifies dramatically 
when a program based on Secs.\ 4 and 5 is created. In fact, scanning over those partitions in which the elements 
are pentagonal numbers represents a totally different application from the examples studied in Sec.\ 5. Therefore, a code is 
developed which only determines the partitions whose elements are pentagonal numbers. However, if $P(z)$ and its inverse 
are expressed as exponentiated double sums, then Theorem\ 1 can be used to derive different expressions for the 
coefficients of the generating functions. In this instance the elements in the partitions are assigned values $\gamma_i$, 
which are obtained from summing the divisors or factors of $i$ divided by $i$. As a result, it is found that the only 
difference between the discrete partition numbers and the partition function via this alternative approach is that the former 
set of numbers possess an extra phase factor of $(-1)^{N_k}$ inside the partition operator, where again $N_k$ is the total 
number of elements in a partition.   
   
Sec.\ 7 begins with an extension of the inverse of $P(z)$, where the coefficients of $z^k$ in the product are now equal to 
$C_k$ rather than -1. Theorem\ 3 shows that a generating function can be obtained from this product where the coefficients are 
determined by applying the discrete partition operator $L_{DP,k}[\cdot]$ acting with each element $i$ assigned a value of 
$C_i$. Conversely, this theorem implies that any power series can be expressed as an infinite product. If the $C_k$ equal the
parameter $\omega$, then the coefficients of the generating function for the product $Q(z,\omega)$ not only become polynomials 
of degree $k$ in $\omega$ giving rise to the discrete partition polynomials $q(k,\omega)$, but also the powers of $\omega$ 
yield the number of elements involved in the discrete partitions. Then more identities involving the discrete partition polynomials 
are derived before a corollary to Theorem\ 3 appears, the latter dealing with the case where the factors in the product are taken 
to an arbitrary power $\rho_k$. Hence, one can derive generating functions for very complicated products with varying powers 
and/or with different factors accompanying the powers of $z$. 

Sec.\ 7 concludes by looking at the special case in the corollary to Theorem\ 3 where the arbitrary powers $\rho_k$ and coefficents 
$C_k$ are set equal to the constant value $\rho$ and the parameter $\omega$, respectively. In this case the coefficients of the 
generating function become $q(k,\omega,\rho)$ and are now polynomials of degree $k$ in both $\rho$ and $\omega$. Since the special 
cases of $\rho \!=\!2$ and $3$ feature in well-known products studied by Euler and Gauss, the properties and values of 
$q(k,\omega,2)$ and $q(k,\omega,3)$ are also examined in detail. 

Sec.\ 8 deals with the derivation of the generating functions for even more complicated products than those appearing in Sec.\ 
7. As a result of the success of introducing the parameter $\omega$ in the inverse of $P(z)$, the first example deals with
the introduction of $\omega$ into $P(z)$. This results in a generating function for the new product $P(z,\omega)$ whose coefficients 
$p(k,\omega)$ are polynomials of degree $k$ in $\omega$ that reduce to the partition function $p(k)$ when $\omega$ is set equal to unity. 
These partition function polynomials are expressed in terms of the partition operator acting with each element $i$ assigned to the 
discrete partition polynomials or rather $q(i,-\omega)$. They are found to possess many interesting properties, while their 
coefficients represent the number of partitions in which the number of elements is given by the power of $\omega$. Out of this 
analysis interesting recurrence relations are obtained for the number of discrete partitions or $q(i,1)$.

As was the case for $P(k)$, $P(k,\omega)$ and its inverse can also be expressed as an exponentiated double sum, both of which can 
be handled by Theorem\ 1. Thus, it is found that the partition function and discrete partition polynomials can be be expressed in 
terms of the partition operator with the main difference being that in the case of the former polynomials the elements 
$i$ are assigned to the polynomials $\gamma_i(\omega)$, while for the latter they are assigned to $-\gamma_i(\omega)$. These new 
polynomials represent the extension of the $\gamma_i$ mentioned above. Their coefficients are equal to divisors $d$ of $i$ divided 
by $i$, while each power of $\omega$ is equal to the reciprocal of the coefficient. Moreover, they reduce to the $\gamma_i$ when 
$\omega \!=\! 1$. Because this is a somewhat unusual situation involving divisors, a program is presented that evaluates the 
partition function and discrete partition polynomials in symbolic form so that they can be processed by Mathematica. This means 
that the final forms for the polynomials can be obtained by evaluating the divisor polynomials via the Divisors[k] routine in
Mathematica rather than having to create a separate program to solve this problem.  

Sec.\ 8 continues with the introduction of an arbitrary power $\rho$ into the product $P(z,\omega)$ and determining the coefficients 
$p(k,\omega,\rho)$ of the generating function for this extended product. These coefficients are given in terms of the partition
operator acting with the elements $i$ assigned to minus the partition function polynomials, viz.\ $-p(i,\omega)$, multiplied
by the Pochhammer symbol of $(-\rho)_{N_k}$. Next the generating function of the product of $Q(z,-\beta \omega)$ with $P(z,\alpha \omega)$
is studied. This infinite product denoted by $P(z,\beta \omega,\alpha \omega)$ combines the properties of discrete partitions with 
standard partitions. The coefficients of the resulting generating function, which are denoted by $QP_k(\omega,\alpha,\beta)$,
are polynomials of degree $k$ in $\omega$, while the powers of $\alpha$ and $\beta$ indicate the number of elements in the standard 
and discrete partitions respectively. With this result Sec.\ 8 concludes with the derivation of the generating function for 
Heine's product, which can be represented as the product of $P(z,\omega x,\omega)$ and $P(z,\omega y,\omega x y)$. The 
coefficients of the generating function for this infinite product, which arises in q-hypergeometric function theory, are denoted by 
$HP_k(\omega,x,y)$ and are obtained by summing the product of $QP_j(\omega,x,1)$ with $QP_{k-j}(\omega,y,xy)$ for $j$ ranging 
from 0 to $k$. Several coefficients in the last two examples are tabulated in order to display their complicated nature.

\section{Generating Partitions}
As indicated in the introduction the partition method for a power series expansion is composed of two major steps.
The first step involves determining all the partitions summing to an integer $k$, which represents the order 
of the variable in the resulting power series expansion. The second and more complicated step is to calculate 
the contribution that each partition makes to the coefficient of the $k$-th order term in the series expansion. This step
will be described extensively in the next two sections when the theoretical framework and programming methodology for the 
partition method for a power series expansion are presented. For now, however, this section is devoted to the problem of 
generating partitions in a suitable format to enable the second step of the partition method for a power series expansion to 
proceed. This means that we shall not only be interested in determining the parts or elements in a partition, but also with 
evaluating their number of occurrences or frequencies. So, whilst the generation of partitions is an interesting problem in 
its own right and continues to be the source for new algorithms as evidenced by Refs.\ \cite{zog98}-\cite{yam07},
it is required here in order to develop a programming methodology for the partition method for a power series expansion. 
Hence, we need to examine the existing algorithms for generating partitions to find which, if any, is the most suitable 
for implementation in the partition method for a power series expansion. Ultimately, we shall find that the novel algorithm 
sketched out in Ref.\ \cite{kow10} will prove to be the most suitable. Moreover, we shall see in Sec.\ 5 that partitions 
with specific properties can be determined by making modifications to the algorithm, which is often difficult to 
achieve, if not impossible, with the other partition-generating algorithms. That is, a completely different algorithm 
is usually required to solve for each specific property of partitions, whereas only minor changes to the bi-variate 
recursive central partition algorithm presented in this section are needed to generate partitions with these specific 
properties. In addition, as a result of the material in appearing Secs.\ 3 and 4, we shall be able to re-formulate the 
partition method for a power series expansion in terms of a partition operator denoted by $L_{P,k}[ \cdot]$. The 
modifications to the BRCP algorithm presented in Sec.\ 5 will mean that we are effectively programming different operators, 
which will, in turn, lead to the presentation of new and fascinating results when we study the various generating functions 
belonging to the theory of partitions in Sec.\ 6.    

As described in Ref.\ \cite{kow10}, when applying the partition method for a power series expansion there 
is actually no need to generate all the partitions at each order. For example, one can write down all the partitions 
necessary for evaluating the first few orders on a sheet of paper. Once they have been determined, one can then 
proceed to the second step of determining the specific contribution made by each partition to the coefficient of 
the series expansion. The problem occurs when we wish to evaluate the higher order terms, especially if our
ultimate aim is to derive an extremely accurate approximation to the original function. When we need
to go to higher orders, the complexity increases dramatically due to the exponential increase in the number 
of partitions. Then it is no longer feasible to write down all the partitions and carry out the calculations
to determine their contributions to the coefficient at a particular order. Consequently, a programming methodology 
is required for all orders despite the fact that this will ultimately become very slow for very high orders of the 
series expansion due to combinatorial explosion. Nevertheless, we shall see that in developing this programming methodology 
we shall uncover very interesting results for the first time in the theory of partitions. For example, it has already 
been stated that series expansion obtained via the partition method produces a power series that is identical to a 
Taylor/Maclaurin series when the latter can be evaluated. In these situations the development of a 
programming methodology for the partition method means that the higher order derivatives in such a series can 
be expressed in terms of a sum of the contributions from all the partitions summing to a particular order. 
This has profound implications in mathematics in that we now have a means of linking the continuous/differentiable 
property of a function with the discrete mathematics of partitions or number theory. 
 
When it was mentioned above that there was no need to employ an algorithm to generate all the partitions in the 
partition method for a power series expansion, it was meant that there was no need to write them down in 
for each value of $k$, since only the elements in each partition and their frequencies are required as the input for 
the second step of the partition method for a power series expansion. Representing a partition in this manner is known 
as the multiplicity representation, whereas we shall refer to the representation where each element in a partition is 
written down as the standard representation. As $k$ increases, the number of partitions summing to $k$ or $p(k)$  
increases exponentially. E.g., the number of partitions summing to 100, viz.\ $p(100)$, is $190\,596\,292$. As a result, 
it is no longer practical to write the partitions in the standard representation. Whilst the multiplicity representation 
is sufficient for the application of the power series for a power series expansion, if one wishes to determine high orders of 
the resulting power series expansion via the method, then one still needs to consider the various algorithms 
for generating partitions because it could turn out that an algorithm generating partitions in the standard 
representation may require only minor modifications to provide them in the multiplicity representation. Furthermore, 
the generation of integer partitions continues to attract interest up to the present time. Therefore, we shall
review the existing algorithms, but ultimately our aim will be to determine the most appropriate for the partition method 
for a power series expansion. On the other hand, those with only an interest in generating partitions may find the other 
algorithms more suitable in which case they are urged to obtain more information by consulting the list of references.   

Having justified the need for generating all the partitions summing to an arbitrary integer, we now
turn to the issue of finding an appropriate algorithm that will expedite the process, but will do so
in appropriate form for the second step of the partition method for a power series expansion. For a time there 
seemed to be only one useful algorithm for generating partitions. This was McKay's algorithm \cite{mck70}, which 
was basically a succession rule whereby partitions were generated in linear time. It was developed further 
by Knuth \cite{knu005}, who used the fact that if the last element greater than unity is a two, then the next 
partition can be determined very quickly. This modification means that each partition takes almost a constant 
amount of time to be generated. The Knuth/McKay algorithm, which is implemented in C/C++ below, generates the 
partitions summed to a global integer $n$ in a particular form of the standard representation known as reverse 
lexicographic order. Consequently, the elements in a partition, say $a[1]$ up to $a[k]$, are printed out according to $a[1] 
\geq a[2] \geq \cdots \geq a[k]$, while the first element of each new partition is less than or equal to the 
first element of the preceding partition. For example, the partitions summing to 5 in reverse lexicographic 
order are: \newline 
$ 5| \\
4 |1 | \\
3 | 2 | \\
3 | 1 | 1 | \\
2 | 2 | 1 | \\
2 | 1 | 1 | 1 | \\
1 | 1 | 1 | 1 | 1 |\\ $
The rules for generating partitions in reverse lexicographic order can be obtained from Refs.\ \cite{knu005}
and \cite{ski90}. Briefly, if the partition is not composed entirely of ones, then it ends with a value of 
$x \!+\!1$ followed by zero or more ones. The next smallest partition in lexicographic order is obtained
by replacing the segment of the partition $\{ \dots,x+1,1,...,1\}$ by $\{\cdots,x,...,x,r\}$, where the 
remainder $r$ is less than or equal to $x$.

\begin{lstlisting}{}

/* This program determines partitions in reverse
   lexicographic order following McKay/Knuth 
   algorithm as discussed on p. 38 of Fascicle 3 
   of Vol. 4 of D.E. Knuth's The Art of Computer 
   Programming. */

#include <stdio.h>
#include <memory.h>
#include <stdlib.h>

int main(int argc,char *argv[])
{
int *a, i, m, n, q, x;
if(argc != 2) printf("execution: ./knuth <partition#>\n");
else
                n=atoi(argv[1]);
                a=(int *) malloc((n+1)*sizeof(int));
P1:     a[0]=0;
                m=1;
P2:     a[m]=n;
                q=m-(n==1);
P3:     for(i=1; i<=m; i++) printf("%i|", a[i]);
                printf( " \n" );
                if(a[q] != 2) goto P5;
P4:     a[q--]=1;
                a[++m]=1;
                goto P3;
P5:     if(q==0) goto end;
                x= a[q]-1;
                a[q]=x;
                n= m-q+1;
                m= q+1;
P6:     if (n<= x) goto P2;
                a[m++]=x;
                n -=x;
                goto P6;
end: ;
    }
printf("\n");
free(a);
return(0);
}

\end{lstlisting}

By current standards the above code is considered to be slow for generating the partitions at each order due 
to the excessive unconditional branching. From a computational point of view it is also very non-structured 
and hence, does not accord with modern programming practice. A significantly faster algorithm for 
generating partitions in reverse lexicographic order has been developed by Zoghbi and Stojmenovic in Ref.\ 
\cite{zog98}. Actually, these authors present two algorithms in their paper, but the second, which generates 
the partitions in lexicographic order, is slower than the first. Nevertheless, if one runs the above code 
against a C/C++-coded version of the first algorithm, then one finds that it takes 1362 CPU seconds to print 
out the partitions summing to 80 on the screen of a Sony VAIO laptop with 2 GB RAM compared with 1399 CPU seconds 
using the Knuth/McKay code given above. On the other hand, if the partitions are directed to an output file, 
then it takes only 30 CPU seconds with the Knuth/McKay code compared with 28 CPU seconds with the Zoghbi/Stojmenovic 
code.

\begin{lstlisting}{}

/* This program determines partitions in ascending order
following the algorithm given at J. Kelleher's web-page.*/

#include <stdio.h>
#include <memory.h>
#include <stdlib.h>

int main(int argc,char *argv[])
{
int *a, i, m, n, ydummy, xdummy, count;
if(argc!= 2) printf("execution: ./kelleher <partition#>\n");
else{
    n=atoi(argv[1]);
    a=(int *) malloc((n+2)*sizeof(int));

    for(i=0; i<=(n+1);i++) a[i]=0;
    a[1]=n;
    count=1;
    while (count != 0){
          xdummy=a[count-1]+1;
          ydummy=a[count]-1;
          count=count-1;
          while (xdummy <= ydummy){
                a[count]=xdummy;
                ydummy=ydummy-xdummy;
                count=count+1;
                                  }
          a[count]=xdummy+ydummy;
          for (i=0;i<= count;i++){
                printf("%i",a[i]);
                if (i<count) printf("|");
                                 }
                printf("\n");
                     }
    }
printf("\n");
free(a);
return(0);
}
\end{lstlisting}

Lexicographic ordering of the partitions ordering of the partitions is not the only method of generating partitions in 
ascending order. Kelleher and O'Sullivan in Refs.\ \cite{kel05} and \cite{kel09} have developed algorithms for generating 
partitions in ascending order, but not in lexicographic order, one of which has been created as the C/C++-code
appearing above. When this code is run for the partitions summing to 5, the following output is obtained:\newline
$1|1|1|1|1 \\
1|1|1|2 \\
1|1|3 \\
1|2|2 \\
1|4 \\
2|3 \\
5\\ $
Now, by inverting the reverse lexicographic order of the partitions summing to 5 given earlier, we see that the partition
\{1,1,3\} in Kelleher/O'Sullivan code appears before \{1,2,2\}, while in McKay/Knuth code it appears after the latter
partition. A similar situation occurs with the \{1,4\} and \{2, 3\} partitions. Yet, the partitions for both codes are 
arranged in ascending order. Moreover, Kelleher and O'Sullivan are able to take advantage of the different ordering between 
to develop an even more efficient version of the above code. As a result, they state that their optimised version is much 
superior to the Zoghbi/Stojmenovic code. Interestingly, if the above code is run to yield the $15\,796\,476$ partitions 
summing to 80 on the same Sony laptop as the previous code, then it takes 1342 CPU seconds to output the partitions. Yet, if
one does the same with the more efficient version of their code, then it takes the same amount of time to generate the 
partitions. On the other hand, if the partitions are directed to an output file, then it is found that the above code takes 
29 CPU seconds, while the more efficient version takes 28 CPU seconds, the same time taken as the Zoghbi/Stojmenovic code.

According to Ref.\ \cite{kel11}, Kelleher has found on his computer system that partitions summing to 80 are generated 
at a rate of $1.30 \times 10^{8}$ per second using the first algorithm, while with the second algorithm the rate is 
$2.87 \times 10^{8}$ per second. For the Zoghbi/Stojmenovic code the rate is $1.26 \times 10^{8}$ per second, while 
with the Knuth/McKay code the partitions are generated at a rate of $1.73 \times 10^{8}$ per second. Hence, the reason why 
there was no marked difference in performance in the C versions of the Kelleher codes is attributed to the manner in which the 
partitions were printed out.  

Unfortunately, the above codes do not utilise the recursive nature of partitions, which can be observed by realising 
that the partitions summing to $k \!+\! 1$ include all the partitions summing to $k$ with an extra element of unity in 
them in addition to other partitions possessing elements greater than unity. In fact, according to p.\ 45 of Ref.\ 
\cite{knu005}, the number of partitions summing to $k$ with $m$ elements, which is denoted by 
$\begin{vmatrix} k \\ m\end{vmatrix}$, obeys the following recurrence relation:
\begin{eqnarray}
\begin{vmatrix} k \\ m\end{vmatrix} =  \begin{vmatrix} k-1 \\ m-1 \end{vmatrix} 
+ \begin{vmatrix} k-m \\ m \end{vmatrix} \;\; . 
\label{one}\end{eqnarray}
In addition, $\begin{vmatrix} k+m \\ m\end{vmatrix}$ represents the number of partitions summing to $k$ with at most $m$
elements. If $P(k,m)$ represents the partitions summing to $k$ with at most $m$ parts, then with the aid of the above 
recurrence relation we obtain
\begin{eqnarray}
P(k,m)= P(k,m-1) + P(k-m,m) \;\;.
\label{onea}\end{eqnarray}
This result is derived on p.\ 96 of Ref.\ \cite{com74}.

To incorporate recursion into the process of generating the partitions, we need to construct an algorithm that utilises the 
special tree diagrams, which first appeared in Ref.\ \cite{kow94}, but have since been applied to other problems or situations in 
Refs.\ \cite{kow10}-\cite{kow11}. These trees are different from the binary tree approach in the recent work of Yamanaka et al
\cite{yam07}, which seeks to generate each partition in the standard representation in a constant time rather than an average
time. Their work represents a further development on Fenner and Loizou \cite{fen80}, who seemed to have been the first to develop 
a binary tree representation for partitions based on a lexicographic ordering. To construct the special tree diagrams, one begins 
by drawing branch lines to all pairs of numbers that can be summed to the seed number $k$, where the first number in the 
tuple is an integer less than or equal to $[k/2]$. Here [$x$] denotes the greatest integer less than or equal to $x$.  
For example, in Fig.\ 1, which displays the tree diagram for the seed number equal to 6, we draw branch lines to (0,6), (1,5), 
(2,4) and (3,3). Whenever a zero appears in the first element of a tuple, the path stops, as evidenced by (0,6). For the other 
pairs, one draws branch lines to all pairs with integers that sum to the second number under the prescription that the first 
member of each new tuple is now less than or equal to half its second member. Hence, for (1,5) we get paths branching out to (0,5), 
(1,4) and (2,3), but not (3,2) or (4,1). This recursive approach continues until all paths are terminated with a tuple containing 
a zero as indicated in the figure. 

\begin{figure}
\vspace*{7.5 cm}
\begin{picture}(440,300)
\thicklines
\put(0,130){6}
\put(0,115){idx(6,1)}
\put(15,130){\line(1,-3){40}}
\put(58,250){(0,6)}

\put(15,130){\line(1,-1){40}}
\put(58,170){(1,5)}
\put(45,150){idx(5,1)}
\put(90,170){\line(1,1){40}}
\put(133,210){(0,5)}
\put(90,170){\line(1,0){40}}
\put(133,170){(1,4)}
\put(120,150){idx(4,1)}
\put(160,170){\line(1,1){40}}
\put(203,210){(0,4)}
\put(160,170){\line(1,0){40}}
\put(203,170){(1,3)}
\put(190,150){idx(3,1)}
\put(230,170){\line(1,1){40}}
\put(273,210){(0,3)}
\put(230,170){\line(1,0){40}}
\put(273,170){(1,2)}
\put(260,150){idx(2,1)}
\put(300,170){\line(1,1){40}}
\put(343,210){(0,2)}
\put(300,170){\line(1,0){40}}
\put(343,170){(1,1)}
\put(330,150){idx(1,1)}
\put(370,170){\line(1,0){40}}
\put(413,170){(0,1)}

\put(160,170){\line(1,-1){40}}
\put(203,130){(2,2)}
\put(190,110){idx(2,2)}
\put(230,130){\line(1,0){40}}
\put(273,130){(0,2)}
\put(90,170){\line(1,-1){40}}
\put(133,130){(2,3)}
\put(120,110){idx(3,2)}
\put(160,130){\line(1,-1){40}}
\put(203,90){(0,3)}

\put(15,130){\line(1,1){40}}
\put(58,90){(2,4)}
\put(45,70){idx(4,2)}
\put(85,90){\line(1,0){40}}
\put(128,90){(0,4)}
\put(85,90){\line(1,-1){40}}
\put(128,50){(2,2)}
\put(115,30){idx(2,2)}
\put(155,50){\line(1,0){40}}
\put(198,50){(0,2)}

\put(15,130){\line(1,3){40}}
\put(58,10){(3,3)}
\put(45,0){idx(3,3)}
\put(85,10){\line(1,0){40}}
\put(128,10){(0,3)}
\end{picture}
\caption{Tree diagram of the partitions summing to 6 for the BRCP algorithm.}
\label{figone} \end{figure}
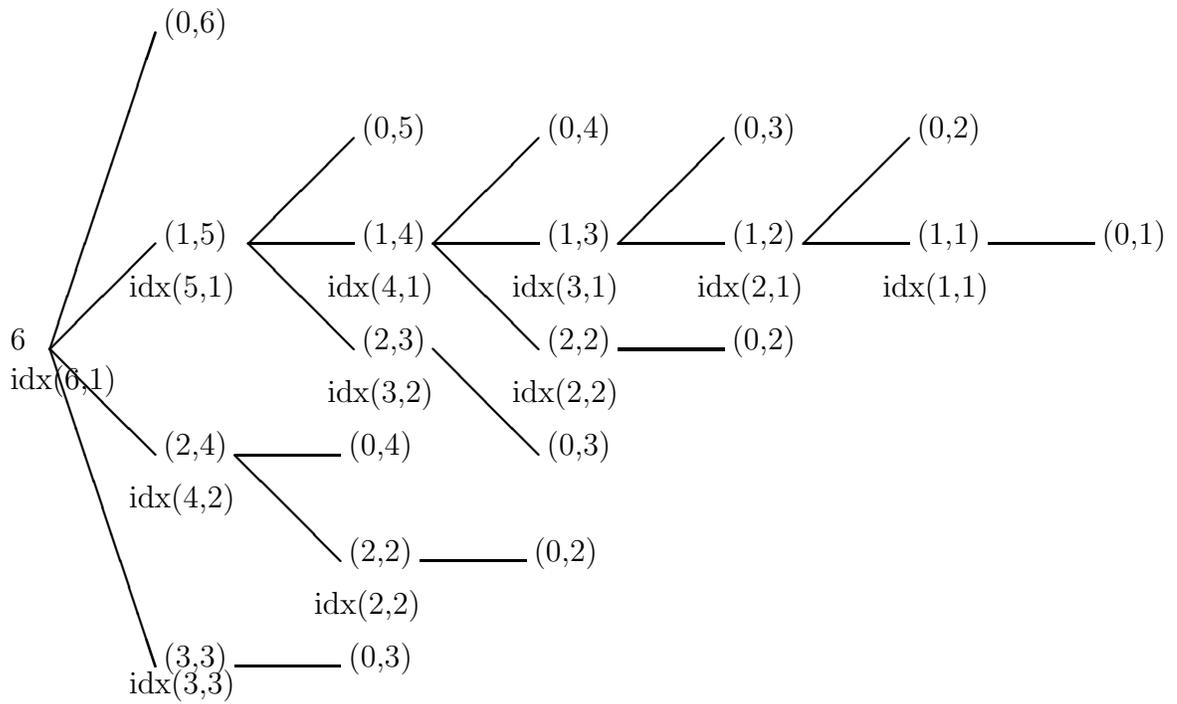

It is obvious that all the first members plus the second member of the final tuple in each 
path represents a partition for $k$. E.g., the path in the figure consisting of (1,5), 
(1,4), (2,2) and (0,2) represents the partition $\{1,1,2,2\}$, while that consisting of
(1,5), (1,4), (1,3) and (0,3) represents the partition $\{1,1,1,3\}$. In addition, the number
of branches along each path represents the number of elements in each partition, whilst those
tuples with zeros in vertical columns represent the partitions with the same number of 
elements in them. For $k \!>\!3$ the last path in the figure consists of ([$k$/2]+1,[$k$/2]) 
and (0,[$k$/2]) when $k$ is odd and ($k/$2,$k$/2) and (0,$k$/2) when $k$ is even. In both even 
and odd $k$ cases the tree diagram terminates at what we shall call the central partition, which 
represents the partition of $\{[k/2]+1,[k/2]\}$ for odd values of $k$ and $\{k/2,k/2\}$ for even 
values of $k$. Unfortunately, this is not all that is required to produce the final tree diagram 
in the figure. Duplicated paths involving permutations of the same partition must also be 
removed so that each partition appears only once in the final diagram. When this removal 
process is carried out, one will eventually end up with the tree diagram displayed in Fig.\ \ref{figone}.   
To determine $\begin{vmatrix} k\\m\end{vmatrix}$, we simply count all the tuples with zero in 
them $m$ branches from the seed number. Hence, $\begin{vmatrix} 6\\3\end{vmatrix}$ equals
the number of tuples with zero in them in the third column, which comes to 3. In addition,
the number of partitions with exactly $m$ parts is the same number as the number of partitions
whose largest element is $m$, which will be become apparent when we discuss conjugate partitions. 

In Ref.\ \cite{kow10} it was stated initially that duplicated partitions could be removed by
the introduction of a search algorithm. Such an algorithm would result in a major increase 
in the complexity of the above graphical approach for generating partitions, which would,
in turn, make the other algorithms presented earlier not only more appealing, but also
far more efficient. To avoid the introduction of a search algorithm, it was stated later
in the same reference that the generation of the partitions could be improved by viewing 
the tree diagrams from a different perspective. Thus, one notes first that the entire tree 
emanating from (1,5) in the diagram is the same tree that one would obtain if the graphical 
method was applied to 5 instead of 6, thereby exhibiting the recursive nature of the method. 
Similarly, the tree emanating from (1,4) is the same tree diagram one would obtain by applying 
the graphical method to 4 instead of 6. This continues all the way down to the last partition 
whose elements are only composed of ones. Moreover, only the partitions with unity in them 
emanate from (1,5) in the diagram, whereas the partitions emanating from (2,4) only 
possess elements greater than or equal to 2. Similarly, the partitions emanating from (3,3) 
only possess elements that are greater than or equal to 3. In fact, in the last instance since 3 
is half of 6, there will only be threes involved along the path from (3,3). Thus, we see that 
the last path represents the central partition $\{3,3\}$, which would have been $\{4,3\}$ had we 
constructed a tree diagram for partitions summing to 7.

From the tree diagram we see that the second number in a tuple decrements with each rightward horizontal
movement or right branch, while the first element of each tuple increments with downward vertical movement. 
In short, the trees are two-dimensional. That is, two variables are required to construct them, a property
first observed by D. Balaic. This is particularly interesting as it means that we are describing a rare 
instance of bi-variate recursion. Hence, there is no need for the introduction of a search algorithm to 
remove duplicated partitions. Since it has been noted that the tree diagrams terminate at the central 
partition, we shall refer to the algorithm that generates the partitions based on the tree diagrams  
as the bi-variate recursive central partition or BRCP algorithm. 

Before presenting an elementary code utilising the BRCP algorithm, let us investigate how the recursive
properties of partitions are included in the tree diagrams such as Fig.\ \ref{figone}. From the 
figure we see that the total number of partitions $p(k)$ can be obtained by summing all the partitions 
of $k$ that can be separated into $m$ elements, where $m$ ranges from 1 to $k$ and $k \!=\! 6$ in 
this instance. Since there is only one partition with one element and one with only $k$ elements, by 
scanning over the columns in the tree diagram we obtain the trivial equation of
\begin{eqnarray}
p(k) =2+ \sum_{m=2}^{k-1} \; 
\begin{vmatrix} k \\ m\end{vmatrix} \;\; , 
\label{oneb}\end{eqnarray}
which is valid for $k \!\geq\! 3$. On the other hand, if we scan the rows of the tree diagram, then we 
see that the total number of partitions can also be obtained by letting $p(k,m)$ represent the number 
of partitions whose elements are greater than or equal to $m$. As a result, we arrive at
\begin{eqnarray}
p(k) = 1+ \sum_{m=1}^{[k/2]} p(k-m,m)  \;\; , 
\label{onec}\end{eqnarray}
where, again, $[x]$ is the greatest integer less than or equal to $x$. This result is given in Ref.\ 
\cite{wik11} except that the variables in $p(k,m)$ have been interchanged and will be used later in Sec.\ 4 
when we introduce the partition operator. 

It has already been stated that $\begin{vmatrix} k+m \\ m\end{vmatrix}$ represents the number of partitions 
summing to $k$ with at most $m$ elements, i.e.\ $P(k,m)$. If we put $m \!=\! 2$ and $k \!=\! 4$, then from the 
diagram that there are three tuples with zeros in them in the column two branches away from the seed number, 
viz.\ \{1,5\}, \{2,4\} and \{3,3\}. Hence, $P(4,2) \!=\!3$. From Eq.\ (\ref{one}) we see that $P(4,2)$ is also 
equal to the sum of $\begin{vmatrix} 5 \\ 1\end{vmatrix}$ and $\begin{vmatrix} 4 \\ 2\end{vmatrix}$. If we treat 
the five in the tuple \{1,5\} in the tree diagram as a seed number, then $\begin{vmatrix} 5 \\ 1\end{vmatrix}$ 
is equal to one corresponding to the tuple \{0,5\}. Furthermore, if we now treat the four in the tuple \{1,4\} 
from \{1,5\} as a seed number, then we find that two branches further to the right $\begin{vmatrix} 4 \\ 2\end{vmatrix}$ 
equals 2 corresponding to the tuples \{0,3\} and \{0,2\}. Note that we could not have used the four in the tuple 
\{2,4\} in the figure as the seed number because this tree diagram gives all the partitions summing to 4, whose elements 
are greater than unity. Including the first branches emanating from the seed number of 4, we see that these 
correspond to the partitions of \{1,3\} and \{2,2\}. Hence, summing the $\begin{vmatrix} 5 \\ 1\end{vmatrix}$ and 
$\begin{vmatrix} 4 \\ 2\end{vmatrix}$, we find once more that $P(4,2)$ is equal to 3, confirming that the tree 
diagrams do indeed possess the recursive properties of partitions.

According to Knuth \cite{knu005}, $\begin{vmatrix} k \\ m \end{vmatrix}$ also represents the number of partitions
summing to $k$, whose largest element is $m$. This connection can be observed by using Ferrers diagrams, which
are studied later in this work. As an example, let us consider $\begin{vmatrix} 6 \\3 \end{vmatrix}$, which can be 
determined by summing all those tuples with a zero in the vertical column three branches from the seed number 
(the standard approach) and is, therefore, equal to 3. Meanwhile, the largest element of a partition always appears 
in the final tuple ending a path in the tree diagram. Hence, the number of partitions whose largest element is 3 can 
be determined by summing all the paths ending with the tuple of (0,3). There are three of these in the tree diagram 
with the first occurring at the top of the fourth column from the seed number, the second at the bottom of the third 
column and the third at the bottom of the second column.  
  
With regard to the BRCP algorithm, an elementary version in C/C++, which first appeared in Ref.\ \cite{kow10}, is  
\begin{lstlisting}{}
void idx(int k,int j) {
  printf("%d",k);
  k=k-j;
  while (k > = j){
    printf(",%d(",j);
    idx(k--, j++);
    printf(")");
                 }
                      }
\end{lstlisting}
Note that the order of the variables in the above code is interchanged compared with the tuples in the tree 
diagram. For $k\!=\!4$ the output from this code is $4,1(3,1(2,1(1))),2(2)$. By processing the commas and 
parentheses, we obtain the partitions in the order they appear in the tree diagram, viz.\ \{4\}, \{1,3\}, 
\{1,1,2\}, \{1,1,1,1\} and \{2,2\}. Although the output is very compact, the code in this form is not suitable 
for the implementation into the second step of the partition method for a power series expansion. In fact, 
the above code has to be adapted in order to solve various problems connected with the theory of partitions 
studied later in Sec.\ 5. Even if we want to list the partitions on separate lines in a similar manner to the
other codes, modifications are required. Nevertheless, the above code does represent the simplest implementation of the 
BRCP algorithm. It is not only more structured and hence, more elegant than the reverse lexicographic algorithm of McKay
and Knuth, but it is also more powerful or versatile. For example, one single call to ${\bf idx(6,1)}$ results 
in all the other calls to the routine as shown in the tree diagram. In fact, the total number of calls to ${\bf idx}$ 
yields the total number of partitions $p(k)$, which is an important quantity in its own right. By introducing 
a counter for the number of calls to ${\bf idx}$ in the above code, we obtain the total number of partitions 
required to construct a tree diagram without the need to create new routine. We shall observe later in this work,
especially in Sec.\ 5, how making only a few changes can result in a host of the special properties 
being determined from partitions.

Let us consider the generation of partitions summing to a particular value as we have done for the other codes 
codes presented in this section. To demonstrate the versatility of the BRCP code, we shall not print the 
partitions as in the same manner as these codes. Instead, we shall generate the partitions in the multiplicity
representation, which is required for the second step of the partition method for a power series expansion. 
Consequently, a function prototype called ${\bf termgen}$ needs to be introduced into the BRCP code called
${\bf partgen}$ as displayed below. When it is called, it will compute the frequencies of the elements by 
counting the same elements in each partition, the latter being represented by the array called ${\it part}$. 
In order to facilitate the call to ${\bf termgen}$, ${\bf idx}$ has undergone minor modification so that in 
processing the partitions the ones are counted first, the twos next and so on. As a result, the following 
output is obtained for $k \!=\!5$:\newline
$1: 1(5)\\ 
2: 1(1) 1(4) \\
3: 2(1) 1(3) \\
4: 3(1) 1(2) \\
5: 5(1) \\
6: 1(1) 2(2) \\
7: 1(2) 1(3) \\$ 
In the code given below the variable ${\it term}$ represents a rolling count of the partitions as they 
are being determined by ${\bf idx}$. In the above output the first value printed out on each line
is the value of ${\it term}$, which is followed by a colon. Hence, the final value of ${\it term}$ 
represents the total number of partitions for each value of $k$ or the variable ${\it tot}$ in the code. 
Each line of output only gives the nonzero values of the frequencies of the elements accompanied by the 
values of the elements presented in parentheses. For example, 1(5) denotes the partition \{5\} where 
$n_5 \!=\!1$ and all the other $n_i$ equal zero, while 3(1) 1(2) represents \{1,1,1,2\}, in which case 
$n_1 \!=\! 3$ and $n_2 \!=\! 1$. As expected, the final partition is the central partition for 
$k \!=\! 5$, viz.\ \{2,3\}. 

As a comparison, the code given below was run on the same Sony laptop as the other codes, where 
it was found that it took 1561 CPU seconds to compute all the partitions summing to 80 according to the 
above format. Hence, the execution time compared with the other codes has increased, primarily due to 
the extra processing of the partitions. However, if the output is directed to a file, which can, in turn,
be used as input to the partition method for a power series expansion, then it only takes 26 CPU seconds 
to execute, which makes it the best performing code in this mode. 

Whilst the other algorithms/codes discussed in this section may prove to be faster than even an optimised 
version of the BRCP algorithm in other situations, they do not possess the versatility or flexibility of 
the latter. We shall use this versatility when we embark on programming the partition method for a power 
series expansion in Sec.\ 4. It should also be noted that Refs.\ \cite{knu005} and \cite{yam07} present 
extra algorithms for generating partitions according to a specific number of parts, whilst the latter reference 
present another algorithm which outputs doubly-restricted partitions or where the elements lie in a specified 
range. In Sec.\ 5 we shall see that only minor modifications to the BRCP algorithm are required to solve these 
problems. That is, there is no need to create an entirely different algorithm to solve such problems.    

\begin{lstlisting}{}
#include <stdio.h>
#include <memory.h>
#include <stdlib.h>

int tot,*part;
long unsigned int term=1;

void termgen()
{
int freq,i;
printf("%ld: ",term++);
for(i=0;i<tot;i++){
                  freq=part[i];
                  if(freq) printf("%i(%i)",freq,i+1);
                  }
printf("\n");
}


void idx(int p,int q)
{
part[p-1]++;
termgen();
part[p-1]--;
p -= q;
while(p >= q){
        part[q-1]++;
        idx(p--, q);
        part[q++ -1]--;
             }
}

int main(int argc,char *argv[])
{
int i;
if(argc !=2) printf("partgen <sum of the partitions>\n");
else{
    tot=atoi(argv[1]);
    part=(int *) malloc(tot*sizeof(int));
    if(part==NULL) printf("unable to allocate array\n\n");
    else{
         for(i=0;i<tot;i++) part[i]=0;
         idx(tot,1);
         free(part);
         }
     }
printf("\n");
return(0);
}


\end{lstlisting}

Another interpretation of the tree diagram in Fig.\ \ref{figone} is to realise that the first branch from
the seed number is the only single element partition, viz.\ \{6\}. The partitions summing to $n \!-\!1$ 
or 5 in this case appear along the second branch to \{1,5\} in the tree diagram. If the partitions summing 
to $n-1$ have already been stored in an array, then all we need to do is increment the number of ones by one in 
all these partitions to get the partitions for $n$. The next branch emanating from the seed number represents 
the tree diagram for $n\!-\!2$, but now all the elements in the partitions are greater than unity. So, if the 
partitions summing to $n\!-\!2$ have been stored previously, then we disregard those partitions where the number 
of ones is non-zero and increment the number of twos in the remaining partitions by one to get the partitions 
summing to $n$. The next branch from the seed number represents the tree diagram for $n\!-\!3$ except that 
those partitions, where the number of ones and twos are non-zero, are now neglected. To get the partitions that 
sum to $n$, we increment the number of threes by one in the remaining partitions. This process continues until 
$[n/2]$ tree diagrams have been processed. Although the new interpretation may lead to less processing of 
the partitions, it comes at the expense of having to store all previous partitions in memory. Nevertheless, 
we shall return to this interpretation in the next section. 

\section{The Partition Method for a Power Series Expansion}
As mentioned earlier, the second step in the partition method for a power series expansion is the more
important step. It involves coding partitions so that each makes a distinct contribution to the coefficient
in the resulting power series expansion. As discussed in Refs.\ \cite{kow09} and \cite{kow11} these
contributions need not be numerical in nature as was the case in the computation of the reciprocal logarithm numbers 
in Ref.\ \cite{kow10}. They can also be functions or polynomials, which will become evident later in this work.
 
The contribution to a coefficient made by each partition in a tree diagram is not only dependent upon the total number of 
elements or parts in the partitions, $N$, but also on the elements in the partitions. Originally, when the method
was devised, the elements were set equal to $l_i$, while $n_i$ represented the number of occasions or frequency each 
$l_i$ appeared in a particular partition. Therefore, if there were $j$ elements in a partition, then 
$\sum_{i=1}^j n_i \!=\! N$, while $\sum_{i=1}^j n_i l_i \!=\! k$, which represents the order of the power series expansion. 
Later, it was decided to let $n_1$ represent the number of ones in a partition, $n_2$ represent the number of twos, $n_3$, 
the number of threes and so on. Then $\sum_{i=1}^k n_i \!=\!N$. The reason why the former approach was adopted 
initially was that it eliminated the redundancy caused by the fact that often, many of the $n_i$ were equal to
zero in the partitions. For example, $n_6$ equals unity in only one partition in the tree diagram, while for this 
partition $n_1$ to $n_5$ equal zero. On the other hand, the problem with the first approach is that $N$ varies 
and this means that writing down general formulae or expressions involving all the partitions is far more awkward. 
Furthermore, the redundancy due to the fact that many of the $n_i$ are equal to zero in the partitions has no 
effect on the tree diagram and consequently, on the the BRCP algorithm. Hence, we shall adopt the second approach when 
discussing partitions, which seems to be the generally accepted approach used by mathematicians \cite{knu005}.

In Refs.\ \cite{kow10}-\cite{kow11} the partition method was applied to specific instances where the standard method
of deriving Taylor/Maclaurin series expansions breaks down. However, Taylor/Maclaurin series expansions are simpler
to derive when the original function is differentiable. In cases where both the partition method can be applied and a 
Taylor/Maclaurin series can be derived, they yield identical results even if the resulting expansion is divergent. 
In addition, it was described extensively in both Refs.\ \cite{kow09} and \cite{kow11} how the method can be 
extended to situations where the coefficients may be dependent upon a variable rather than provide a pure number. 
Therefore, it should be possible to develop a general theorem describing the partition method for a
power series expansion as is the case with Taylor/Maclaurin series. Of course, such a theorem will need to indicate 
under what conditions the method is valid. Moreover, in Refs.\ \cite{kow09} and \cite{kow11} it was described 
how the method could be inverted to yield the power series expansion of the reciprocal function. That is, if the 
method can be applied to a function $f(x)$, then it could often be applied to $1/f(x)$. A theorem on the method
would also need to indicate under what conditions it can be inverted.    
 
Before we present the theorem describing the partition method for a power series expansion, we need to introduce 
some preliminaries. First, in the literature a composite function $g \circ f(x)$ is defined as being equal to 
$g(f(x))$. Here we shall define a pseudo-composite function $g_a \circ f$ as being equal to $g(a f(x))$, where $a$ 
need not necessarily be a number. Next we need to explain the concept of regularisation. For more detailed descriptions
of this process the reader is referred to Refs.\ \cite{kow09}-\cite{kow11} and \cite{kow09a}-\cite{kow95}. According 
to this concept, when a power series representation for a function is divergent, particularly as happens with
an asymptotic expansion, it needs to be regularised in order to yield meaningful values that are representative
of the original function. As a consequence, when it is uncertain that a power series representation is convergent 
or when it is known to be divergent, we cannot use the equals sign in a mathematical statement. Instead, we introduce 
the less stringent equivalence symbol and refer to the resulting expression as an equivalence statement. For example, in 
Refs.\ \cite{kow09} and \cite{kow09a}-\cite{kow001} it is shown that the geometric series, i.e. $\sum_{k=0}^{\infty} z^k$ 
is absolutely convergent for $|z| \!<\! 1$, conditionally convergent for $\Re \, z \!<\!1$ and $|z| \!\geq\!1$, 
undefined for $\Re\, z \!=\! 1$ and divergent for $\Re\, z \!>\! 1$. In the last two instances it is simply invalid 
to say that the series is equal to anything. However, through the process of regularisation it is found that the 
regularised value of the geometric series is the same value of $1/(1-z)$ that one obtains as the limit value when 
the series is either conditionally or absolutely convergent. Furthermore, this value is bijective and hence, 
unique for $z$ lying in the principal branch of the complex plane. Therefore, when the series is divergent, we 
can only say that it is equivalent to its regularised value. Because the regularised value is equal to the limit 
value of the series when it is convergent, which is not always the case as discussed in Ch.\ 4 of Ref.\ \cite{kow09a}, 
we replace the equals sign by the equivalence symbol for all values of $z$. That is, we can express the geometric 
series as an equivalence statement or equivalence for short by the following statement:
\begin{eqnarray}
\sum_{k=N}^{\infty} z ^k \equiv \frac{z^N}{1-z}\quad, \quad \forall\;  z\;\;.
\label{two}\end{eqnarray}     

Now we introduce the main theorem in this work.\newline
{\bfseries Theorem\ 1}. 
Given that the function $f(z)$ can be expressed in terms of a power series referred to here as the inner 
power series, in which $f(z) \equiv \sum_{k=0}^{\infty} p_k y^k$ and  $y \!=\! z^{\mu}$, and that the 
function $g(z)$ can be expressed in terms of another or outer power series, i.e.\ $g(z) \equiv h(z) 
\sum_{k=0}^{\infty} q_k z^k$, where $h(z)$ is an arbitrary function or number, then for non-zero values 
of $p_0$ there exists a power series expansion for the quotient of the pseudo-composite functions 
$g_a \circ f$ and $h_a \circ f$ given by
\begin{eqnarray}
\frac{g_a \circ f}{h_a \circ f} \equiv  \sum_{k=0}^{\infty} D_k \, y^k  \;\;.
\label{three}\end{eqnarray}
In the above equivalence the first few coefficients $D_k$ are given by
\begin{eqnarray}
D_0 = F(ap_0) \;\;, \;\; D_1=a F^{(1)}(ap_0) \, p_1\;\;, 
\label{three-a}\end{eqnarray} 
and 
\begin{eqnarray}
D_2= \frac{a^2}{2}\; F^{(2)}(ap_0) \, p_1^2+ aF^{(1)}(ap_0) \,  p_2 \;\;. 
\label{three-b}\end{eqnarray}
A general formula for the coefficients can be derived by analysing the partitions summing to $k$. This yields 
\begin{eqnarray}
D_k = \sum_{\scriptstyle n_1, n_2,n_3, \dots,n_k=0 \atop{\sum_{i=1}^k in_i =k}}^{k,[k/2],
[k/3],\dots,1} a^N F^{(N)}(a p_0) \prod_{i=1}^k \frac{p_i^{n_i} }{n_i !} \;\;.
\label{four}\end{eqnarray} 
In this equation $\sum_{i=1}^k n_i \!=\! N$, while $F^{(N)}(ap_0)$ represents the $N$-th derivative of the 
function $F(a p_0)$, which, in turn, represents the regularised value of the power series expansion  
$\sum_{j=0}^{\infty}q_j (ap_0)^j$. That is, $F(ap_0) \equiv \sum_{j=0}^{\infty} q_j (ap_0)^j$. For the 
important case where $p_0 \!=\! 0$, the coefficients are given by
\begin{eqnarray}
D_k = \sum_{\scriptstyle n_1, n_2,n_3, \dots,n_k=0 \atop{\sum_{i=1}^k in_i =k}}^{k,[k/2],
[k/3],\dots,1} q_N \, a^N N!  \prod_{i=1}^k \frac{p_i^{n_i} }{n_i !} \;\;,
\label{five}\end{eqnarray}
Moreover, for $D_0 \!\neq\!0$, the inverted quotient of the pseudo-composite functions can also be expressed 
in terms of a power series and is given by
\begin{eqnarray}
\frac{h_a \circ f}{g_a \circ f} \equiv \frac{1}{D_0}  \sum_{k=0}^{\infty} E_k y^k \;\;.
\label{six}\end{eqnarray} 
Here the coefficients $E_k$ are found to be 
\begin{eqnarray}
E_k = \sum_{\scriptstyle n_1, n_2,n_3, \dots,n_k=0 \atop{\sum_{i=1}^k=k}}^{k,[k/2],
[k/3],\dots,1} (-1)^N N! \; D_0^{-N} \prod_{i=1}^k \frac{D_i^{n_i} }{n_i !} \;\;.
\label{seven}\end{eqnarray}
Finally, the coefficients $D_k$ and $E_k$ satisfy the following recurrence relation: 
\begin{eqnarray}
\sum_{j=1}^k D_j\,E_{k-j} = 0 \;\;. 
\label{sevena}\end{eqnarray}
 
{\bfseries Remark 1}. As in the case of a Taylor/Maclaurin series the power series given by Equivalences\ (\ref{three}) and 
(\ref{six}) can be either (1) convergent for all values of the variable, (2) absolutely convergent within a finite radius of 
convergence or (3) asymptotic, which is defined here as a power series expansion with zero radius of absolute convergence.  

{\bfseries Remark 2}. The second result for the $D_k$, viz.\ Eq.\ (\ref{five}), is similar in form to the definition on p.\
134 of Ref.\ \cite{com74} for the partial or second type of Bell polynomial. In fact, the latter are a special case of the
above theorem, which can be obtained by setting $q_k \!=\! a \!=\! 1$ and $p_k \!=\! x_k/k!$.

{\bfseries Proof}. Since $g(z)$ can be expressed in terms of a power series expansion and the function $h(z)$, 
we have 
\begin{eqnarray}
\frac{g_a \circ f}{h_a \circ f}  \equiv \left( q_0 +  \sum_{k=1}^{\infty} q_k \, a^k f(z)^k \right) \;\;.
\label{eight}\end{eqnarray}
Introducing the power series expansion for $f(z)$ into the above result yields
\begin{align}
\frac{g_a \circ f}{h_a \circ f} & \equiv \left( q_0 +  q_1 \,a \sum_{k=0}^{\infty} p_k\, y^k
+ q_2 \,a^2 \Bigl( \sum_{k=0}^{\infty} p_k \,y^k \Bigr)^{\!2} \right.
\nonumber\\
& \left.  + \;\; q_3\,a^3 \Bigl( \sum_{k=0}^{\infty} p_k\, y^k \Bigr)^{\! 3}+
 q_4\,a^4 \Bigl( \sum_{k=0}^{\infty} p_k\, y^k \Bigr)^{\! 4} + \cdots \right) \;\;.
\label{nine}\end{align}
Isolating the zeroth order term of the power series expansion for $f(z)$ in the above result yields
\begin{align}
\frac{g_a \circ f}{h_a \circ f} & \equiv \left( q_0 +  q_1 \,a \Bigl( p_0 + \sum_{k=1}^{\infty} 
p_k\, y^k \Bigr)+ q_2 \,a^2 \Bigl( p_0+ \sum_{k=1}^{\infty} p_k \,y^k \Bigr)^{\!2} \right.
\nonumber\\
& \left.  + \;\; q_3\,a^3 \Bigl( p_0 + \sum_{k=1}^{\infty} p_k\, y^k \Bigr)^{\! 3}
+  q_4\,a^4 \Bigl( p_0+ \sum_{k=1}^{\infty} p_k\, y^k \Bigr)^{\! 4} + \cdots \right) \;\;.
\label{ten}\end{align}
Expanding in descending powers of $p_0$ yields
\begin{align}
\frac{g_a \circ f}{h_a \circ f} & \equiv \sum_{j=0}^{\infty} q_j \left(a\, p_0\right)^{\!j} 
+ \sum_{j=1}^{\infty} j\, q_j\, a^j p_0^{j-1} \sum_{k=1}^{\infty} p_k\, y^k + \sum_{j=2}^{\infty} \binom{j}{2} q_j \,a^j\, 
p_0^{j-2} 
\nonumber\\
& \times \;\;\left( \sum_{k=1}^{\infty} p_k \,y^k \right)^{\!2} 
+ \sum_{j=3}^{\infty} \binom{j}{3} q_j\,a^j\, p_0^{j-3} \left( \sum_{k=1}^{\infty} p_k\, y^k \right)^{\! 3}
+ \cdots  \;\;.
\label{eleven}\end{align}

Let us now represent the regularised value of the series $\sum_{k=0}^{\infty} q_k (a p_0)^k$ by $F(ap_0)$. That is,
\begin{eqnarray}
\sum_{k=0}^{\infty} q_k (a p_0)^k \equiv F(a p_0) \;\;.
\label{twelve}\end{eqnarray}
We can use this definition to simplify the sums on the rhs of the Equivalence\ (\ref{eleven}), but first we note
that
\begin{eqnarray}
\sum_{j=i}^{\infty} \binom{j}{i} q_j a^j p_0^{j-i} = \frac{a^i}{i!} \sum_{j=1}^{\infty} q_j(a p_0)^{j-1} 
\prod_{l=0}^{i-1}(j-l)= \frac{a^i}{i!} \frac{d^i}{dz^i} \sum_{j=0}^{\infty} q_j z^j \Bigg\vert_{z=ap_0}.
\label{thirteen}\end{eqnarray}
Introducing Equivalence\ (\ref{twelve}) into the above result yields
\begin{eqnarray}
\sum_{j=i}^{\infty} \binom{j}{i} q_j a^j p_0^{j-i} \equiv \frac{a^i}{i!} \frac{d^i}{dz^i} F(z) \Bigg\vert_{z=ap_0} \!\!\!\!
= \frac{a^i}{i!} \, F^{(i)}(ap_0) \;\; .
\label{fourteen}\end{eqnarray}
Consequently, Equivalence\ (\ref{eleven}) can be expressed as
\begin{eqnarray}
\frac{g_a \circ f}{h_a \circ f} \equiv \sum_{k=0}^{\infty} \frac{a^k}{k!} \, F^{(k)}(a p_0) \left( \sum_{j=1}^{\infty} p_j y^j
\right)^{\!k} \;\;.
\label{fifteen}\end{eqnarray}

If Equivalence\ (\ref{fifteen}) is expanded in powers of $y$, then we obtain
\begin{align}
& \frac{g_a \circ f}{h_a \circ f}  \equiv F(ap_0)+ a\, F^{(1)}(ap_0) p_1 y + \Bigl( \frac{a^2}{2} \, F^{(2)}(ap_0) \, p_1^2
+ a\, F^{(1)}(ap_0) p_2 \Bigr) y^2
\nonumber\\
& + \;\; \Bigl( \frac{a^3}{3!}\, F^{(3)}(ap_0) + a^2 F^{(2)}(ap_0) \, p_1 p_2+ a F^{(1)}(a p_0)\, p_3 \Bigr) y^3 + O(y^4)
\;\;.
\label{sixteen}\end{align}
From this result we see that the coefficients of the zeroth, first and second order terms in $y$ correspond to the 
results for $D_0$, $D_1$ and $D_2$ given in the theorem. Moreover, we see that
\begin{eqnarray}
D_3 = \frac{a^3}{3!}\, F^{(3)}(ap_0)\,p_1^3 + a^2 F^{(2)}(ap_0) \, p_1 p_2+ a F^{(1)}(a p_0)\, p_3 \;\;.
\label{seventeen}\end{eqnarray}

The first few coefficients are relatively easy to write down, but beyond that it becomes progressively more difficult.
It is at this stage we need to introduce partitions into the analysis to facilitate the derivation of a general 
expression for the coefficients of the powers of $y$ in Equivalence\ (\ref{sixteen}). 
 
If we look closely at the result for $D_3$, then we see that it is the sum of three separate contributions, which is due to
the fact that there are only three partitions summing to 3, namely \{1,1,1\}, \{1,2\} and \{3\}. As stated in Ref.\ 
\cite{kow10}-\cite{kow11}, in order to evaluate the specific contribution due to each partition, a value must be 
assigned to the elements appearing in the partitions. In particular, Ref.\ \cite{kow09} states that these assigned values 
depend upon the coefficients of the inner series, i.e.\ the power series expansion for $f(x)$, which becomes the variable in 
the outer power series for $g(z)$. Hence, each element $i$ in a partition is assigned a value of $p_i$. Furthermore, 
there is also a multinomial factor associated with each partition. This factor is not only dependent 
upon the frequencies $n_i$ of the elements in a partition, but also on their sum, $N \!=\! \sum_{i=1}^k n_i$. It arises
from the fact that the value of $k$ on the rhs of Equivalence\ (\ref{fifteen}) is used to render each partition. That is, 
$k$ in Equivalence\ (\ref{fifteen}) corresponds to $N$ in the partition method for a power series expansion. Often the 
multinomial factor simply becomes $N!/n_1! \, n_2! \cdots n_k!$, but in Equivalence\ (\ref{fifteen}) there is an extra 
factor of $a^k\, F^{(k)}(a p_0)/k!$ in the terms. Again, as $k$ in Equivalence\ (\ref{fifteen}) plays the role of $N$
in the partition method for a power series expansion, this means that the standard multinomial factor must be multiplied by 
a factor of $a^N\, F^{(N)}(a p_0)/N!$ for each partition. Therefore, the contribution from a partition to the overall 
coefficient is given by
\begin{eqnarray}
C\Bigl(n_1,n_2, \dots ,n_k \Bigr)= a^N\,F^{(N)}(a p_0) \prod_{i=1}^k \frac{p_i^{n_i}}{n_i !} \;\;.
\label{eighteen}\end{eqnarray}
We see from this result that we do not actually require the partitions themselves, but their frequencies in order to
evaluate the coefficients in the power series for the quotient of the pseudo-composite functions. Nevertheless, it 
should be noted that each set of frequencies identifies a distinct partition.    

To make the preceding material more understandable, let us evaluate the fourth order term in $y$ in Equivalence\ 
(\ref{fifteen}) or $D_4$, which is determined by considering all the partitions summing to 4. There are 5 of these: 
\{1,1,1,1\}, \{1,1,2\}, \{2,2\},\{1,3\} and \{4\}. For the first partition, $n_1 \!=\! 4$, while the other $n_i$ 
are zero. Therefore, from Eq.\ (\ref{eighteen}) we have
\begin{eqnarray}
C(4,0,0,0)=a^4\,F^{(4)}(a p_0)\,p_1^4/4! \;\;.
\label{nineteen}\end{eqnarray}
In the case of the second partition $n_1 \!=\! 2$, $n_2 \!=\! 1$ and the other $n_i$ vanish, while for the third 
partition, only $n_2$(=2) does not vanish. According to Eq.\ (\ref{eighteen}), the contributions due to both  
partitions are
\begin{eqnarray}
C(2,1,0,0)=a^3\,F^{(3)}(a p_0)\,p_1^2\, p_2/3! \;\;,
\label{twenty}\end{eqnarray}
and 
\begin{eqnarray}
C(0,2,0,0)=a^2\,F^{(2)}(a p_0)\,p_2^2/2! \;\;.
\label{twentyone}\end{eqnarray}
For the fourth partition $n_1 \!=\! n_3 \!=\! 1$ with $n_2$ and $n_4$ equal to zero, while in the final partition 
only $n_4$ (=1) is non-zero. Hence, these partitions yield
\begin{eqnarray}
C(1,0,1,0)=a^2\,F^{(2)}(a p_0)\,p_1\, p_3/2! \;\;,
\label{twentytwo}\end{eqnarray}
and 
\begin{eqnarray}
C(0,0,0,1)=a\,F^{(1)}(a p_0)\,p_4 \;\;.
\label{twentythree}\end{eqnarray}
Hence, $D_4$ is given by the sum of five or $p(4)$ contributions.

To derive a general formula for the coefficients $D_k$, we need to sum over all the partitions summing to $k$. 
This entails summing over all values that the frequencies of the elements can take. Each frequency, $n_i$, 
can only range from 0 to $[k/i]$, since this is the maximum number that the element $i$ can appear in a 
partition summing to $k$. In addition, the partitions are constrained by the condition that 
\begin{eqnarray}
\sum_{i=1}^k n_i i=n_1+2n_2+3n_3+\cdots+ k n_k=k \;\;.
\label{twentyfour}\end{eqnarray}
From this result we see that $n_k$ equals zero for all partitions except \{k\}, in which case it will equal unity, 
while the other $n_i$ will equal zero for this partition. That is, the $n_i$ are more often zero than non-zero. 
Consequently, much redundancy occurs when summing over the allowed values of $n_i$.  More succinctly, the sum over 
all partitions summing to $k$ can be expressed as 
\begin{eqnarray}
D_{k} = \sum_{\scriptstyle n_1, n_2,n_3, \dots,n_{k}=0 \atop{\sum_{i=1}^{k}i n_i =k}}^{k,[k/2], [k/3],\dots,1}  
C(n_1,n_2, \dots,n_k) \;\;.
\label{twentyfive}\end{eqnarray}
where $N \!=\! \sum_{i=1}^k n_i$. If Eq.\ (\ref{eighteen}) is introduced into the above equation, then we obtain the
general formula given by Eq.\ (\ref{six}). 

Now we consider the case of $p_0 \!=\!0$. This means that Equivalence\ (\ref{ten}) reduces to 
\begin{eqnarray}
\frac{g_a \circ f}{h_a \circ f}  \equiv  q_0 + \sum_{k=1}^{\infty}  q_k \,a^k \Bigl( \sum_{j=1}^{\infty} 
p_j\, y^j \Bigr)^k \;\;.
\label{twentysix}\end{eqnarray}
Expanding the first few powers in $y$ yields
\begin{align}
\frac{g_a \circ f}{h_a \circ f} & \equiv  q_0 + q_1 \,a p_0\,y+ \Bigl( q_1 a p_2+ q_2 a^2 p_1^2 \Bigr) y^2 
+ \Bigl( q_1 a p_3 
\nonumber\\
& + 2 q_2 a^2 p_1\, p_2+ q_3 a^3 p_1^3 \Bigr) y^3 + O \!\left(y^4 \right) \;\;.  
\label{twentyseven}\end{align} 
Hence, we obtain a power series expansion in $y$ except now the coefficients $D_k$ take a  different form. In particular, 
we see that the first few coefficients are given by
\begin{eqnarray}
D_0 = q_0\;\;, \;\; D_1=q_1 a \, p_0\;\;, \;\; D_2= q_1 a p_2 +q_2 a^2 p_1^2 \;\;,
\label{twentyeight-a}\end{eqnarray} 
and 
\begin{eqnarray}
D_3= q_1\,a p_3 + 2 q_2 \; a^2 p_1 \,p_2 + q_3 a^3 p_1^3 \;\;. 
\label{twentyeight-b}\end{eqnarray}
Furthermore, we can introduce the partition method to derive a general formula for coefficients since Equivalences\ 
(\ref{fifteen}) and (\ref{twentysix}) are isomorphic. In fact, the only difference between the two equivalences is that 
the terms in the outer series in the latter equivalence are multiplied by $q_k a^k$ as opposed to  $a^k F^{(k)}(ap_0)/k!$
in Equivalence\ (\ref{fifteen}). This means that the only change to the partition method will occur in the multinomial factor,
which will now be $q_N a^N N!/n_1! n_2!\cdots n_k!$. Consequently, the contribution by each partition to the $D_k$ is given 
by  
\begin{eqnarray}
C\Bigl(n_1,n_2, \dots ,n_k \Bigr)= q^N\,a^N\, N! \prod_{i=1}^k \frac{p_i^{n_i}}{n_i !} \;\;.
\label{twentynine}\end{eqnarray}
Introducing this result into Eq.\ (\ref{twentyfive}) gives Eq.\ (\ref{four}).

Although it has been stated that the expansion given by Equivalence\ (\ref{three}) can be asymptotic, this does not
mean that it will be divergent for all values of $y$. For if it were, then inverting the equivalence implies that the 
inverted quotient of the pseudo-composite functions vanishes for all values of $y$. Therefore, there must be a 
region in the complex plane where the equivalence symbol can be replaced by an equals sign. In Ref.\ \cite{kow09} it
was found that an asymptotic series possesses zero radius of absolute convergence, but it is either conditionally convergent
or divergent depending upon which sector in the complex plane the variable is situated. For the values of $y$ where
Equivalence\ (\ref{three}) is convergent, we can invert the quotient of the pseudo-composite functions $g_a \circ f$ 
and $h_a \circ f$. Therefore, provided $D_0 \! \neq \!0$, we find that
\begin{eqnarray}
\frac{h_a \circ f}{g_a \circ f} = \frac{1}{D_0} \; \frac{1}{ \left(1+ \sum_{k=1}^{\infty} ( D_k /D_0 ) \, y^k \right)}  \;\;.
\label{thirty}\end{eqnarray}
The rhs of this equation can now be regarded as the regularised value of the geometric series where the variable equals
to $-\sum_{k=1}^{\infty} (D_k/D_0)y^k$. According to Refs.\ \cite{kow09} and \cite{kow09a}-\cite{kow002}, the geometric 
series or $1+z+z^2+ \cdots$ is either conditionally or absolutely convergent for $\Re\,z \!<\! 1$. For all other values
of $z$ in the principal branch of the complex plane, it is either divergent or undefined and must be regularised. Hence, 
treating the rhs of the above equation as the limit of the geometric series means that $\Re\, \sum_{k=1}^{\infty}
(D_k/D_0)y^k \!>\! -1$.  For all other values of $y$ the series will be either divergent or undefined, the latter
case occurring when $\Re\, \sum_{k=1}^{\infty} (D_k/D_0) y^k \!=\!-1$. Therefore, Eq.\ (\ref{thirty}) can be expressed
as
\begin{eqnarray}
\frac{h_a \circ f}{g_a \circ f} \equiv  \frac{1}{D_0} \; \sum_{k=0}^{\infty} \left(- \sum_{j=1}^{\infty} 
( D_j /D_0 ) \, y^j \right)^{\!k}  \;\;.
\label{thirtyone}\end{eqnarray}

Equivalence\ (\ref{thirtyone}) is isomorphic to Equivalence\ (\ref{fifteen}), which means in turn that we can apply
the partition method again. In this instance the coefficients of $y^j$ in the inner series or rather the values to
be assigned to the elements $j$ in the partitions are equal to $-D_j/D_0$ instead of $p_j$.  In addition, the multinomial
factor becomes the standard value of $N!/n_1! n_2! \cdots n_k!$ for the partitions summing to $k$. Hence, the contribution
from each partition is given by
\begin{eqnarray}  
C\Bigl(n_1,n_2, \dots ,n_k \Bigr)= (-1)^N\,N!\,D_0^N\, \prod_{i=1}^k \frac{D_i^{n_i}}{n_i !} \;\;.
\label{thirtytwo}\end{eqnarray}
Introducing the above into Eq.\ (\ref{twentyfive}) yields the result given by Eq.\ (\ref{seven}). 

To derive the final result, we use Equivalences (\ref{three}) and (\ref{six}) to obtain
\begin{eqnarray}
1= \frac{g_a \circ f}{h_a \circ f} \; \frac{h_a \circ f}{g_a \circ f} \equiv \sum_{k=0}^{\infty} D_k y^k
\Bigl( \frac{1}{D_0} \sum_{k=0}^{\infty} E_k \,y^k \Bigr) \;\;.
\label{thirtythree}\end{eqnarray}
Since $E_0 \!=\!1$ from Eq.\ (\ref{seven}), we can separate the zeroth order term in $y$, which cancels the term of 
unity on the lhs of the equivalence. After multiplying the series together, we are left with
\begin{eqnarray}
\sum_{k=1}^{\infty} y^k \sum_{j=0}^{k} D_j E_{k-j} \equiv 0 \;\;.
\label{thirtythreea}\end{eqnarray}
As indicated earlier in the proof, there will be some values of $y$, actually a region in the complex plane, where 
the equivalence symbol can be replaced by an equals sign. Otherwise, either Equivalence\ (\ref{three}) or (\ref{six})  
would be zero for all values of $y$. Given that there will be an infinite number of values of $y$ where an equals 
sign applies, the lhs of the above result will also vanish for these values of $y$. For this to occur, it means 
that the inner series, which is independent of $y$, must also vanish. Hence, we arrive at Eq.\ (\ref{sevena}), thereby 
completing the proof of the theorem.

It should be noted that Theorem\ 1 has avoided the issue of determining the radius of absolute convergence for the power
series on the rhs of Equivalence\ (\ref{three}). This is because although one can determine a value for the radius 
of absolute convergence when deriving the resulting power series, it is often only an estimate, not the supremum as
demonstrated by various examples in Ref.\ \cite{kow11}. Furthermore, one can introduce a divergent power series as the 
inner series and despite the fact that the outer series may also possess a finite radius of absolute convergence, the 
resulting power series appearing in Equivalence\ (\ref{three}) may be convergent over the entire complex plane. For 
example, the power series expansion for cosecant derived via the partition method for a power series expansion in 
Ref.\ \cite{kow11} has a radius of absolute convergence equal to $\pi$. This inner series became the variable for 
the outer series, which was given by the geometric power series, whose radius of absolute convergence is unity. Yet, 
the resulting power series in Equivalence\ (\ref{three}) was merely another representation of the standard Taylor/Maclaurin
series for sine, which is convergent for all values of the variable. So, this is an example where both the inner and outer
series were both absolutely convergent within different radii in the complex plane, but the resulting series expansion obtained 
via the partition method for a power series expansion was convergent for all values of the variable.    

To make the preceding material clearer, let us consider an example. Power series expansions for $\csc z$ and $\sec z$ have 
already been obtained in Ref.\ \cite{kow11} via the partition method for a power series expansion. Because $\cos z$ and 
$\sin z$ can be expressed as Bessel functions of half-integer order, namely $J_{-1/2}(z)$ and $J_{1/2}(z)$, we can also
say that power series expansions have been developed for the reciprocal or inversion of these special functions for 
$\nu \!=\! -1/2$ and $\nu \!=\! 1/2$. Whilst various extensions of these results are presented in Ref.\ \cite{kow11}, one 
which has been overlooked is the derivation of a power series expansion for the reciprocal of a Bessel function to arbitrary 
order $\nu$. Therefore, if we introduce the standard power series expansion for Bessel functions given by No.\ 8.440 in Ref.\ 
\cite{gra94} into the denominator, then we obtain
\begin{align}
&J_{\nu}(z)^{-1} = \frac{(2/z)^{\nu}\, \bigl(1/\Gamma(\nu+1) \bigr)} {(1-z^2/4(\nu+1)+ z^4/16 \cdot 2!(\nu+2)
(\nu+1)-\cdots)} \;\;.
\label{thirtythreea1}\end{align}
The above result is not valid for $\Re \, \nu \!=\! -1$ since in this case the leading term in the series expansion for
$J_{\nu}(z)$ vanishes. Then we would need to examine this case by itself, which is left for the reader to consider.
From Theorem\ 1 we see that the the inner power series is simply the Taylor/Maclaurin power series expansion for 
Bessel functions. Hence, $y \!=\! z^2$, and $p_0 \!=\!0$, while for $k \geq 1$,  $p_k = (-1)^k/2^{2k} (\nu+1)_k k!$, where 
$(\nu+1)_k$ is the Pochhammer notation for $\Gamma(k+\nu+1)/\Gamma(\nu+1)$. Moreover, we can regard the denominator 
on the rhs of Eq.\ (\ref{thirtythreea1}) as the regularised value of the geometric series, which means that coefficients of the 
outer series, viz.\ $q_k$, are equal to $(-1)^k$ for $k \geq 1$, while the function $h(z) = (2/z)^{\nu}\, \Gamma(\nu+1)$. 
Since $a \!=\! 1$, the pseudo-composite functions become composite functions with the quotient in Theorem\ 1 equal to 
$(z/2)^{\nu}/\Gamma(\nu+1)\,J_{\nu}(z)$. By using Equivalence (\ref{three}) we arrive at
\begin{eqnarray}
\frac{(z/2)^{\nu}}{\Gamma(\nu+1)\, J_{\nu}(z)} \equiv \sum_{k=0}^{\infty} h_k(\nu) \,\Bigl(\frac{z}{2} \Bigr)^{2k} \;\;,
\label{thirtythreea2}\end{eqnarray}
where the coefficients $h_k(\nu)$ are determined from Eq.\ (\ref{five}) and are given by
\begin{align}
&h_k(\nu)=(-1)^k \sum_{\scriptstyle n_1, n_2,n_3, \dots,n_k=0 \atop{\sum_{i=1}^{n_k}i n_i =k}}^{k,[k/2],
[k/3],\dots,1} \!\! \! \! \! (-1)^N N! \prod_{i=1}^k \Bigl( \frac{1}{(\nu+1)_i \, i!} \Bigr)^{n_i} \; \frac{1}{n_i!} \;\;.
\label{thirtythreea2a}\end{align}
From this result we see that $h_0(\nu) \!=\! 1$, $h_1(\nu) \!=\! 1/(\nu+1)$, and $h_2(\nu) \!= \!(\nu+3)/2(\nu+1)^2 (\nu+2)$,
while for $k \!=\!3$, Eq.\ (\ref{thirtythreea2}) yields
\begin{align}
h_3(\nu)&= \left( \Bigl( \frac{\Gamma(\nu+1)}{\Gamma(\nu+2)}\Bigr)^{\!3} \!\!- 2! \; \frac{\Gamma(\nu+1)}
{\Gamma(\nu+2)} \; \frac{\Gamma(\nu+1)} {2! \cdot\Gamma(\nu+3)} + \frac{\Gamma(\nu+1)}
{3! \cdot \Gamma(\nu+4)} \right)
\nonumber\\
&=  \;\; \left( \frac{\nu^2+8\nu+19}{3! \cdot(\nu+1)^3 (\nu+2) (\nu+3)} \right) \;\;.
\label{thirtythreea3}\end{align}

Since $J_{1/2}(z)$ is related to $\sin z$, the power series expansion for $1/J_{1/2}(z)$ is identical
to the one derived in Ref.\ \cite{kow11} for cosecant, whose coefficients were expressed in terms of the cosecant
numbers denoted by $c_k$.  These numbers were later found to be related to the Riemann zeta function. Hence,
we obtain 
\begin{eqnarray}
h_k(1/2) \!=\!2^{2k} c_k = 2 \left( 2^{2k} -2 \right) \, \frac{\zeta(2k)}{\pi^{2k}}\;\;.
\label{thirtythreea4}\end{eqnarray}   
Similarly, because $J_{-1/2}(z)$ is related to $\cos z$, the expansion for secant, which is derived in terms of special
numbers called the secant numbers or $d_k$ in Ref.\ \cite{kow11}, is related to $\nu \!=\! -1/2$ in Equivalence\ 
(\ref{thirtythreea1}). It was then found that the secant numbers could be expressed as the difference of specific 
values of Hurwitz zeta function. Thus, we arrive at
\begin{eqnarray} 
h_k(-1/2) \!=\! 2^{2k} d_k = \frac{1}{\pi^{2k+1}}\; \Bigl( \zeta\left( 2k+1,1/4 \right)-\zeta \left(2k+1,3/4
\right) \Bigr) \;\;. 
\label{thirtythreea5}\end{eqnarray}

It was also observed in Ref.\ \cite{kow11} that the absolute convergence of the power series series expansions for both 
cosecant and secant were determined by the distance from the origin to the first zeros of these functions. In the case 
of cosecant the expansion is absolutely convergent for $|z| \!<\! \pi$, while in the case of secant it is absolutely 
convergent for $|z| \!<\! \pi/2$. Therefore, the expansion given in Equivalence\ (\ref{thirtythreea2}) will be absolutely 
convergent for $|z|$ less than the magnitude of the first zero for $J_{\nu}(z)$. For these values of $z$ we can replace 
the equivalence symbol by an equals sign, thereby producing an equation. Consequently, we can use the ``equation form" 
of Equivalence\ (\ref{thirtythreea2}) to demonstrate that the coefficients $h_k(\nu)$ can also be evaluated by 
recurrence relations. Because there is an infinite number of values of $z$ where the ``equation form" is valid, a 
recurrence relation can be be obtained simply by multiplying Equivalence\ (\ref{thirtythreea2}) by the cnovergent 
power series expansion for $J_{\nu}(z)$ and setting the resulting product equivalent to unity. Since $z$ is still fairly 
arbitrary within the radius of absolute convergence, we can equate like powers of $z$ of the resulting equation. This yields
\begin{eqnarray}
h_k(\nu) =\sum_{j=0}^{k-1} \frac{(-1)^{k-j+1}\; \Gamma(\nu+1)\, h_j(\nu)}{(k-j)!\,\Gamma(k-j+\nu+1)}\;\;.
\label{thirtythreea6}\end{eqnarray}

Other recurrence relations can be developed by using integral results involving Bessel functions. For example, one 
can express No.\ 3.768(9) in Ref.\ \cite{gra94} as
\begin{eqnarray}
\int_0^1 dx\;x^{\mu-1} (1-x)^{\mu-1} \, \frac{\cos(ax)}{J_{\mu-1/2}(a/2)} =\sqrt{\pi}\,  a^{1/2-\mu} 
\,\Gamma(\mu) \cos \!\left(\frac{a}{2}\right) \;\;. 
\label{thirtythreea7}\end{eqnarray}
We have already seen that there exists a finite radius of absolute convergence for $z$, where Equivalence\ (\ref{thirtythreea}) 
becomes an equation. This means that there is a region in the complex plane where the expansion for the reciprocal of
the Bessel function can be introduced into the above result without the necessity to replace the equals sign by an
equivalence symbol. Furthermore, by introducing the power series expansion for cosine into Eq.\ (\ref{thirtythreea7}) 
we arrive at
\begin{align}
\int_0^1 dx \; x^{\mu-1} & (1-x)^{\mu-1} \sum_{k=0}^{\infty} a^{2k} \sum_{j=0}^k
\frac{(-1)^{k-j}}{2^{4j}} \frac{x^{2k-2j}}{(2k-2j)!}\,h_j(\mu-1/2)
\nonumber\\
& \equiv \frac{\sqrt{\pi}\,\Gamma(\mu)}{2^{2\mu-1} \Gamma(\mu+1/2)} \sum_{k=0}^{\infty} 
\frac{(-1)^{k}}{(2k)!} \left(\frac{a}{2}\right)^{\!2k} \;\;.
\label{thirtythreea8}\end{align}
The integral on the lhs of the above equivalence is merely the integral representation for the beta function. Hence, 
we can introduce the gamma function product for the beta function. Since $a$ is fairly arbitrary, we can equate 
like powers of $a$ on both sides of the equation, thereby obtaining 
\begin{align}
\sum_{j=0}^{k} \frac{(-1)^{k-j}}{2^{4j}\, (2k-2j)!}& \, \frac{\Gamma(2k-2j+\nu+1/2)}{\Gamma(2k-2j+2\nu+1)} 
\;  h_j(\nu) 
\nonumber\\
& =  \frac{\sqrt{\pi}}{2^{2k+2\nu}} \frac{(-1)^k}{(2k)!} \frac{1}{\Gamma(\nu+1)} \;\;,
\label{thirtythreea9}\end{align}
where $\mu-1/2$ has been replaced by $\nu$.

If we put $k\! = \!4$ into either recurrence relation and introduce the values of $h_k(\nu)$ for $k \!=\! 1$ to 3 
given earlier, then we find that
\begin{eqnarray}
h_4(\nu) = \frac{\nu^4+17 \nu^3+117 \nu^2 + 379 \nu +422 }{4! \,(\nu+1)^4 \, (\nu+2)^2 \,(\nu+3) 
(\nu+4)} \;\;.
\label{thirtythreea10}\end{eqnarray}
This agrees with the value given in Table\ \ref{table1}, which displays the values of $h_k(\nu)$ up to
$k \!=\! 7$.  These values have been evaluated by introducing the recurrence relation into the Sum routine 
in Mathematica \cite{wol92}. In the next section we shall develop a programming methodology where summing
the contributions from the partitions in Eq.\ (\ref{thirtythreea2a}) will be just as expedient as using 
the recurrence relations.

\begin{table}
\begin{tabular}{|c|c|} \hline
$k$ &  $h_k(\nu)$  \\ \hline
$0$ & $ 1$ \\  
$1$ &  $ \frac{1}{(\nu+1)}$  \\ 
$2$ & $ \frac{\nu+3}{2(\nu+1)^2 (\nu+2)}  $ \\ 
$3$ & $ \frac{\nu^2 +8 \nu+19}{3! \cdot (\nu+1)^3(\nu+2)(\nu+3)}$  \\ 
$4$ & $ \frac{\nu^4 +17 \nu^3+117 \nu^2 +379 \nu+422}
{4! \cdot (\nu+1)^4(\nu+2)^2(\nu+3)(\nu+4)}$ \\ 
$5$ & $ \frac{\nu^5 +26 \nu^4+294 \nu^3+1816 \nu^2+5969 \nu+7302}
{5! \cdot (\nu+1)^5(\nu+2)^2(\nu+3)(\nu+4)(\nu+5)}$ \\ 
$6$ & $ \frac{\nu^8 +42 \nu^7+811 \nu^6+9412 \nu^5+71155 \nu^4+349786 \nu^3+1043637 
\nu^2 +1674616 \nu+1091052}
{6! \cdot (\nu+1)^6 (\nu+2)^3 (\nu+3)^2 (\nu+4) (\nu+5)(\nu+6)}$ \\  
$7$ & $ \frac{\nu^9 +55 \nu^8 +1417 \nu^7 +22535 \nu^6 +243311 \nu^5 +
1827401 \nu^4 +9292435 \nu^3 +29539597 \nu^2 +51572980 \nu +36978156}
{7! \cdot (\nu+1)^7 (\nu+2)^3 (\nu+3)^2 (\nu+4) (\nu+5) (\nu+6) (\nu+7)}$
\\ \hline
\end{tabular}
\caption{Coefficients for the power series expansion of the reciprocal of the Bessel function of 
order $\nu$.}
\label{table1}
\end{table}

From the table it can be seen that the $h_k(\nu)$ possess common properties. For example, there is
always a factor of $k! (1+\nu)^k$ in the denominator. In fact, if the denominator of a function 
$f(x)$ is denoted by $DN(f(x))$, then we have 
is given by
\begin{eqnarray}
DN\bigl(h_k(\nu) \bigr) =k! \,(\nu+1)^k\, (\nu+2)^{[k/2]}\cdots (\nu+k-1)^{[k/(k-1)]}\,(\nu+k)\;\;,
\label{thirtythreea11}\end{eqnarray}
where $[x]$ denotes the greatest integer less than or equal to $x$. This result also follows from Eq.\
(\ref{thirtythreea2a}) when we examine the upper limits in the summations. On the other hand, the highest order 
term in $\nu$ in the numerator is always $k$ less than the highest order term in the denominator. Therefore, 
for large $|\nu|$, $h_k(\nu) \approx \nu^{-k}/k!$.

\begin{figure}
\begin{center}
\includegraphics{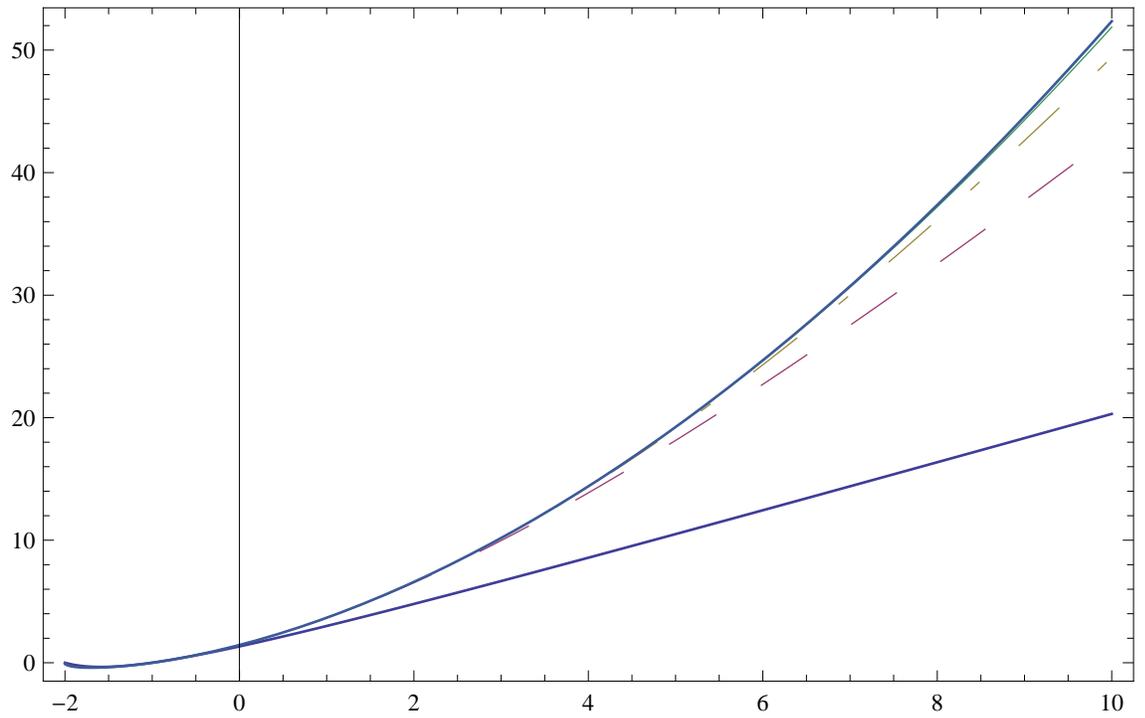}
\end{center}
\caption{$h_k(\nu)/h_{k+1} (\nu)$ versus $\nu$ for $k= 1,5,8,12$ and 18}
\label{figtwo}
\end{figure}

By applying d'Alembert's ratio test as described on p.\ 24 of Ref.\ \cite{cop76} to the series on 
the rhs of Equivalence\ (\ref{thirtythreea2}), we find that for all $k$ the series is only absolutely 
convergent when $|z| < 2 \sqrt{|h_k(\nu)/h_{k+1}(\nu)|}$. Therefore, if there is a supremum for this ratio 
that applies to all values of $k$, then it represents the radius of absolute convergence. Fig.\ \ref{figtwo} 
presents a graph of $h_{k}(\nu)/h_{k+1}(\nu)$ versus $\nu$ for $-2 \!<\! \nu \!<\! 10$ and various values 
of $k$. The figure shows that the larger $k$ is, the greater the value of the ratio $h_k(\nu)/h_{k+1}(\nu)$. 
Therefore, the limit as $k \to \infty$ represents represents the supremum for the ratio when $\nu \!>\! -2$ 
with the nearest singularities to the origin occurring at    
\begin{eqnarray}
z_{0,\nu} = \pm 2 \lim_{k \to \infty}\sqrt{\frac{h_{k}(\nu)}{h_{k+1}(\nu)}} \;\;.
\label{thirtythreea12}\end{eqnarray}
The figure also shows that as $k$ increases, $h_k(\nu)/h_{k+1}(\nu)$ becomes a better approximation to 
$\lim_{k \to \infty} h_k(\nu)/h_{k+1}(\nu)$ over an ever-increasing range of values for $\nu$. Eg., we see
that $h_5(\nu)/h_6(\nu))$ is an accurate approximation for  $-2 \!<\! \nu \!<\! 4.5$, while $h_8(\nu)/
h_9(\nu)$ is an even better approximation for $-2 \!<\! \nu \!<\! 7$. 

According to p.\ 372 of Ref.\ \cite{abr65}, for $\nu \!>\! -1$, the first zero for $J_{\nu}(z)$ is real, but 
for $\nu \!<\!-1$, excluding negative integers, or $\nu$ complex, the first zero is complex. Since the 
first zero represents the singularity of $1/J_{\nu}(z)$ nearest to the origin, we expect the rhs to be
real only for $\nu > -1$. Table\ \ref{table2} presents a sample of the values for the first zero of $J_{\nu}(x)$ 
for various values of $\nu$ obtained by using the special numerical routine known as BesselJZero in Mathematica
and by setting $k \!=\! 17$ in Eq.\ (\ref{thirtythreea12}). For $\nu \!>\! -1$  we see that the values obtained
via both approaches agree. In particular, for $-1 \!<\! \nu \!<\! 1$ the values obtained from Eq.\ 
(\ref{thirtythreea12}) were found to be accurate to at least 10 decimal places. As expected, as $\nu$ becomes greater 
than unity, the accuracy of Eq.\ (\ref{thirtythreea12}) begins to wane compared with the values from BesselJZero. 
The values presented in the second and third columns for these values of $\nu$ only agree to the first decimal
place. Mathematica is able to evaluate the zeros to arbitrary precision, as evidenced by the values for
$\nu$ equal to -1 and -3/2 in the second column of the table. On the other hand, if one wishes to
obtain more accurate values of the zeros via Eq.\ (\ref{thirtythreea12}), then one will have to determine the 
$h_{k}(\nu)$ for much greater values of $k$.

\begin{table}
\begin{tabular}{|c|c|c|} \hline
$\nu$ &  ${\rm BesselJZero}$ & $\pm 2\sqrt{h_{17}(\nu)/h_{18}(\nu)}$  \\ \hline
$0$ & $2.404\,825\,557\,695\,7$ & $\pm 2.404\,825\,557\,695\,4$ \\ 
$1$ & $3.831\,705\,97$ & $\pm 3.831\,705\,96$ \\ 
$2$ & $5.135\,622\,3$ & $ \pm 5.135\,622\,1$ \\ 
$5$ & $8.771\,4 $ & $8.771\,3$ \\ 
$10$ & $14.47$ & $ 14.46$ \\ 
$-1/3$ & $1.866\,350\,858\,873\,4$ & $1.866\,350\,858\,873\,8$ \\ 
$-4/5$ & $0.936\,806\,664\,511$ & $0.936\,806\,664\,510$ \\ 
$-3/2$ & $2.798\,386\,045\,783\,878$ & $\pm 1.199\,678\,640\,257\,655\,i$ \\ 
$-1/3+i$ & ${\rm Unable \; to\;  evaluate}$ & $\pm 2.076\,341\,434\,394\,476 $ \\
$$ & $$ & $\pm 1.556\,637\,759\,994\,043\,i$ \\ 
$3/2-i$ & ${\rm Unable \; to\;  evaluate}$ & $\pm 4.529\,756\,943\,967\,303 $ \\
$$ & $$ & $\mp 1.293\,935\,107\,111\,323\,i$ \\ \hline
\end{tabular}
\caption{Evaluation of the first zero of $J_{\nu}(x)$ using BesselJZero in Mathematica and Eq.\
(\ref{thirtythreea12}) with $k \!=\!17$.}
\label{table2}
\end{table}

For $\nu \!\leq\! -1$, BesselJZero[$\nu$,1] gives the first zero of $J_{\nu}(z)$, but it is the first zero situated on the 
positive real axis, not necessarily the closest to the origin. To see this more clearly, for $\nu \!=\! -3/2$ Mathematica 
gives as its first zero $z \!=\! 2.798\,386\cdots$, but this is not the first zero in the complex plane. According to Eq.\ 
(\ref{thirtythreea12}), the closest zero is given by $z \!=\! 1.199\,687 \cdots i$ as displayed in the adjacent column 
of the table. The latter value can be confirmed as a zero by introducing it into the BesselJ routine in Mathematica, whereupon 
one finds that it yields the tiny value of $(7.35 \cdots+7.35 \cdots i)\times 10^{-14}$. 

It also appears that for $\Re \, \nu \!>\! -2$, Eq.\ (\ref{thirtythreea12}) still gives the nearest zero 
to the origin of the complex plane. For example, the remaining values in the table represent the zeros
obtained for complex orders of Bessel functions. Zeros for such functions cannot be obtained via
BesselJZero since the routine can only handle real numbers. However, Eq.\ (\ref{thirtythreea12}) can
yield them. When the value obtained from Eq.\ (\ref{thirtythreea12}) for $J_{-1/3+i}(z)$ is introduced
into the BesselJ routine in Mathematica, a complex value with a magnitude of the order of $10^{-12}$ is obtained, 
while a complex value with a magnitude of the order of $10^{-8}$ is obtained when the final value in the third 
column of the table is introduced into the routine. This last example demonstrates again that the accuracy of Eq.\
(\ref{thirtythreea12}) wanes as $|\nu|$ increases when $k$ is fixed, which we have already observed in Fig.\
\ref{figtwo}.

The reader may well ask if it is possible to adapt the preceding method to evaluate the next Bessel zero 
or even higher order zeros. This does not seem to be possible at this stage unless the specific value for 
$z_{0,\nu}$ can be determined. For example, by taking the logarithm of No.\ 8.544 in  Ref.\ \cite{gra94}
and differentiating, we eventually arrive at
\begin{eqnarray}
\frac{J_{\nu}^{'}(z)}{J_{\nu}(z)} -\frac{\nu}{z} =-2z \sum_{k=0}^{\infty} \frac{1}{z_{k,\nu}^2-z^2} \;\;. 
\label{thirtythreea13}\end{eqnarray} 
In the above equation, $J_{\nu}^{'}(z)$ can be written as a series in powers of $z/2$ with coefficients $a_k
\!=\! 1/k! \, (\nu+1)_k$, while the rhs of Equivalence\ (\ref{thirtythreea12}) can be used to replace the reciprocal
of the Bessel function. Multiplying both series yields 
\begin{eqnarray}
\frac{J_{\nu}^{'}(z)}{J_{\nu}(z)} -\frac{\nu}{z} \equiv \frac{\nu}{z} \sum_{k=0}^{\infty} \Bigl( \frac{z}{2} 
\Bigr)^k \, e_{k}(\nu) \;\;, 
\label{thirtythreea14}\end{eqnarray}   
where $e_{k}(\nu) =\sum_{j=0}^k a_{k-j}\,h_k(\nu)$. To carry out the singularity analysis for the next zero, we 
need to remove the $k \!=\!0$ term from the rhs of Eq.\ (\ref{thirtythreea13}). Hence, the power series for 
this term would need to be introduced into the lhs with the above power series, yielding another power series with coefficients
expressed in terms of $e_{k}(\nu)$ and powers of $z_{0,\nu}$. If an approximation for the latter is introduced, 
then the coefficients would only be approximate. For large values of $k$ they would still be affected by $z_{0,\nu}$
and thus, it would not be possible to isolate the next zero. On the other hand, if the exact result for the $z_{0,\nu}$,
which is known for $\nu \!=\! \pm 1/2$, then the series would represent the power series for the rhs of 
Eq.\ (\ref{thirtythreea13}) with the sum beginning at $k \!=\! 1$ or with the singularity at the second zero. This 
situation is discussed immediately below Eq.\ (83) in Ref.\ \cite{kow11}. It should be noted that there are typographical
errors there since the power of $1/(2\pi)$ in the expression for $c_k^{*}$ should be $k$, not $k+1$ and the limit 
in the next line below should refer to $c_{k+1}^{*}/c_k^{*}$ , not $c_{k+1}/c_{k}$.  

Before investigating how the BRCP algorithm can be implemented to evaluate the coefficients $D_k$ and $E_k$ in 
Theorem\ 1, we now prove an interesting corollary to the theorem, but before we can do this, the following lemma is 
required:

{\bfseries  Lemma 1}. For $\alpha$ complex, regularisation of the binomial series yields
\begin{eqnarray}
_1 {\mathcal F}_0 (\alpha;z) =\sum_{k=0}^{\infty} \frac{\Gamma(k+ \alpha)}{\Gamma(\alpha) \, k!}\;  z^k 
\begin{cases} =(1-z)^{-\alpha} \;\;, & \quad \Re\, z <1\;\;, \\
\equiv (1-z)^{-\alpha} \;\;, & \quad \Re \, z \geq 1. \end{cases}
\label{thirtythreeb}\end{eqnarray}

{\bfseries Remark}. This lemma represents the generalisation of the regularisation of the geometric series, which
reduces to the latter for $\alpha \!=\! 1$.

{\bfseries Proof}. The proof of Lemma\ 1 appears immediately below Proposition\ 1 in Ch.\ 4 of Ref.\ \cite{kow09a}. 
Although the proof is concerned with real values of $\alpha$, there is no reason why $\alpha$ cannot be complex. This is
due to the fact that the proof involves differentiating and integrating the integral representation for the beta function,
viz.
\begin{eqnarray}
B(k+ \alpha, 1-\alpha) = \int_0^1 dt \; \frac{t^{k+\alpha-1}}{(1-t)^{\alpha}} \;\;.
\label{thirtythreec}\end{eqnarray}
This integral is not only defined for real values of $\alpha$, but also for complex values of $\alpha$. Its convergence is 
limited by the values for the real part of $\alpha$, not its imaginary part. Since differentiation and integration only affect 
the real part of $\alpha$, the proof of Proposition\ 1 in Ref.\ \cite{kow09a} can be extended to complex values of $\alpha$. 
Moreover, the result in Lemma\ 1 can be simplified by replacing the equals sign with an equivalence symbol. Then we 
have the one equivalence statement, which is valid for all values of $z$ and $\alpha$.

One final remark is in order. The reader should note that a similar notation to the standard notation for generalised
hypergeometric functions has been introduced by expressing the binomial series as $_1{\mathcal F}_0$. The reason for the
slight variation in notation is that the generalised hypergeometric function notation of $_pF_q( \alpha_1, \dots,\alpha_p; \beta_1, 
\dots, \beta_q;z)$ is only valid when the series is absolutely convergent, i.e. when $p \leq q $ or for $|z| \!<\!1$ when 
$p \!=\! q$ \cite{pru90}. The conditions in the lemma are more general. That is, the binomial series has been taken into the 
regions of the complex plane where it becomes either conditionally convergent or divergent. Hence, it has been necessary 
to alter the notation to avoid any confusion with the standard notation. This completes the proof of the lemma.

We are now in a position to generalise the derivation of the coefficients in Theorem\ 1 by  considering the following 
corollary.\newline
{\bfseries Corollary 1 to Theorem 1}. Given the same conditions on the pseudo-composite functions $g_a \circ f$ and
$h_a \circ f$ as in Theorem 1, there exists a power series expansion for the quotient of the pseudo-composite functions
raised to an arbitrary power $\rho$ which is given by
\begin{eqnarray} 
\Bigl( \frac{g_a \circ f}{h_a \circ f} \Bigr)^{\rho} \equiv  \sum_{k=0}^{\infty} D_k(\rho) \, y^k  \;\;,
\label{thirtyfour}\end{eqnarray}
For $k \! \geq \! 1$, the generalised coefficients for the $D_k$ of Theorem\ 1 or $D_k(\rho)$ in the above 
result are given by
\begin{eqnarray}
D_k (\rho)= \sum_{\scriptstyle n_1, n_2,n_3, \dots,n_k=0 \atop{\sum_{i=1}^k in_i =k}}^{k,[k/2],
[k/3],\dots,1}  (-\rho)_N \, D_0^{-N}  \prod_{i=1}^k \frac{\left(-D_i \right)^{n_i} }{n_i !} \;\;,
\label{thirtyfive}\end{eqnarray}
and
\begin{eqnarray}
D_k (-\rho)= \sum_{\scriptstyle n_1, n_2,n_3, \dots,n_k=0 \atop{\sum_{i=1}^k in_i =k}}^{k,[k/2],
[k/3],\dots,1} N! \, D_0(\rho)^{-N-1}  \prod_{i=1}^k \frac{\left(-D_i(\rho) \right)^{n_i} }{n_i !} \;\;,
\label{thirtysix}\end{eqnarray}
where $(\rho)_N$ denotes the Pochhammer notation for $\Gamma(N \!+\! \rho)/\Gamma(\rho)$. In addition, for $\rho \!=\! 
\mu \!+\! \nu$ the coefficients satisfy the following recurrence relation: 
\begin{eqnarray}
D_k(\rho) =\sum_{j=0}^k D_{j}(\mu) D_{k-j}(\nu) \;\;.
\label{thirtyseven}\end{eqnarray} 

{\bfseries Remark 1}. By putting $\rho \!=\! -1$ in Eq.\ (\ref{thirtyfive}) we see immediately that the $D_k(-1)$ become the $E_k$
in Theorem\ 1 given by Eq.\ (\ref{seven}).  Furthermore, the $D_k(1)$ represent another representation for the $D_k$ given by
either Eq.\ (\ref{four}) or (\ref{five}), the latter being valid when $p_0 \!=\!0$.

{\bfseries Remark 2}. Note that when $\rho$ equals a non-negative integer $j$, Eq. (\ref{thirtyfive}) simplifies dramatically
because $(-\rho)_N$ vanishes for $N \!>\! j$, reflecting the fact that we are dealing with a finite polynomial of order $j$. 
For all other values of $j$ the coefficients $D_k(\rho)$ represent polynomials in $\rho$ of order $k$.

{\bfseries Proof}. Since the pseudo-composite functions $g_a \circ f$ and $h_a \circ f$ are subject to the conditions
in Theorem\ 1, we know that there exists a power series expansion where
\begin{eqnarray} 
\Bigl( \frac{g_a \circ f}{h_a \circ f} \Bigr)^{\rho} \equiv  \Bigl( \sum_{k=0}^{\infty} D_k \, y^k\Bigr)^{\rho}  \;\;.
\label{thirtyeight}\end{eqnarray}
We also know that there will be a region in the complex plane where the equivalence symbol can be replaced by an equals
sign. Separating the zeroth order term in above result yields
\begin{eqnarray} 
\Bigl( \frac{g_a \circ f}{h_a \circ f} \Bigr)^{\rho} \equiv  D_0^{\rho} \Bigl( 1+ \sum_{k=1}^{\infty} (D_k/D_0 )\, 
y^k\Bigr)^{\rho}  \;\;.
\label{thirtynine}\end{eqnarray}
We can treat the series on the rhs as the variable in the regularised value of the binomial series. Then according to Lemma\ 1
we have
\begin{eqnarray} 
\Bigl( \frac{g_a \circ f}{h_a \circ f} \Bigr)^{\rho} \equiv  D_0^{\rho} \sum_{k=0}^{\infty} \frac{\Gamma(k- \rho)}{\Gamma(-\rho)\, k!}
\; \Bigl(-\sum_{j=1}^{\infty} (D_j/D_0 )\, y^j \Bigr)^{k}  \;\;.
\label{forty}\end{eqnarray}
Now the above equivalence is isomorphic to Equivalence\ (\ref{fifteen}), which means that it is in the form where the 
partition method can be applied. As stated in the proof to Theorem\ 1, the values to be assigned to elements $i$ in the 
partitions are given by coefficients of $y^i$ in the inner series. Therefore, in this case each element $i$ is assigned 
a value of $-D_i/D_0$.  To evaluate the contribution from each partition, we need to multiply the product of all the 
assigned values by a factor consisting of the multinomial factor $N!/n_1!\, n_2! \cdots n_k!$ and $\Gamma(N-\rho)/N!\, 
\Gamma(-\rho)$. The latter term arises from the fact that $k$ in the coefficient of the outer series plays the role
of $N$ in the partition method. Therefore, by introducing the Pochhammer notation of $(\rho)_k$ for 
$\Gamma(N \!+\! \rho)/\Gamma(\rho)$, we find that the contribution due to each partition is 
\begin{eqnarray}   
C\Bigl(n_1,n_2, \dots ,n_k \Bigr)= (-1)^N\, (-\rho)_N\, D_0^{-N} \prod_{i=1}^k \frac{D_i^{n_i}}{n_i!} \;\;,
\label{fortyone}\end{eqnarray}
As expected, for $\rho \!=\! -1$ the above result reduces to Eq.\ (\ref{thirtytwo}). Furthermore, by summing over all
partitions summing to $k$, we obtain the total coefficient in the resulting power series expansion, which is
given by Eq.\ (\ref{thirtyfive}).

As a result of establishing Equivalence\ (\ref{thirtyfour}), we have
\begin{eqnarray} 
\Bigl( \frac{g_a \circ f}{h_a \circ f} \Bigr)^{-\rho} \equiv  \sum_{k=0}^{\infty} D_k(-\rho) \, y^k  \;\;,
\label{fortyonea}\end{eqnarray}
We also know that there will be a region in the complex plane where the equivalence symbol can be replaced by an equals 
sign in Equivalence\ (\ref{thirtyfour}). For these values of $y$ we can invert the equivalence, thereby obtaining
\begin{eqnarray} 
\Bigl( \frac{h_a \circ f}{g_a \circ f} \Bigr)^{\rho} = \frac{1} {D_0(\rho)} \; \frac{1}
{\left(1+\sum_{k=1}^{\infty} D_k(\rho) \, y^k \right)}  \;\;.
\label{fortyoneb}\end{eqnarray}
Theorem\ 1 can again be applied to this result by treating the rhs as the regularised value of the geometric series with
the variable equal to the series in the denominator. In this case $h(z) \!=\! 1$, $q_k \!=\! (-1)^k$, $a \!=\! 1$ 
and $p_k \!=\! D_k(\rho)/D_0(\rho)$ for $k \!\geq \! 0$, while $p_0 \!=\!0$.  The coefficients of the resulting power series 
expansion can be determined by Eq.\ (\ref{five}). Moreover, the resulting power series expansion will be equal to the power 
series expansion on the rhs of Equivalence\ (\ref{fortyonea}) since both have the same regularised value. Hence, we can equate
like powers of both expansions, which yields
\begin{eqnarray}
D_k(-\rho)= D_0(\rho)^{-1} \sum_{\scriptstyle n_1, n_2,n_3, \dots,n_k=0 \atop{\sum_{i=1}^k in_i =k}}^{k,[k/2],
[k/3],\dots,1} N! \, (-1)^N  \prod_{i=1}^k \frac{\left(D_i(\rho)/D_0(\rho) \right)^{n_i} }{n_i !} \;\;.
\label{fortyonec}\end{eqnarray}
By taking the factor of $(-1)^N$ inside the product and $1/D_0(\rho)$ outside of it, we obtain Eq.\ 
(\ref{thirtysix}). 

Since $\rho \!=\! \mu \!+\! \nu$, we have
\begin{eqnarray} 
\Bigl( \frac{g_a \circ f}{h_a \circ f} \Bigr)^{\rho} = \Bigl( \frac{g_a \circ f}{h_a \circ f} \Bigr)^{\mu}  
\Bigl( \frac{g_a \circ f}{h_a \circ f} \Bigr)^{\nu}  \;\;.
\label{fortytwo}\end{eqnarray}
We also note that there will be a region of complex plane for $y$ in Equivalence\ (\ref{thirtyfour}) where we can replace
the equivalence symbol by an equals sign. Therefore, for these values of $y$, Eq.\ (\ref{fortytwo}) yields
\begin{eqnarray}
\sum_{k=0}^{\infty} D_k(\rho)\, y^k = \sum_{k=0}^{\infty} D_k(\nu)\, y^k \sum_{k=0}^{\infty} D_k(\mu) \, y^k 
=\sum_{k=0}^{\infty} y^k\, \sum_{j=0}^k D_j(\nu) D_{k-j}(\mu) \;\;.
\label{fortythree}\end{eqnarray}
Since there is an infinite number values of $y$ for which the above equation holds, we can equate like powers of $y$.
As a consequence, we obtain Eq.\ (\ref{thirtyseven}), which completes the proof of the corollary. 

To make the preceding material more concrete, let us consider a couple of examples. In order to simplify these examples, 
we now regard the sum over partitions as a discrete operator, which will be denoted by $L_P[ \cdot]$. We shall refer to 
this form for the sum over partitions summing to $k$ as the partition operator. That is, the partition operator is defined as
\begin{eqnarray}
L_{P,k}[\cdot] =  \sum_{\scriptstyle n_1, n_2,n_3, \dots,n_{k}=0 \atop{\sum_{i=1}^{k}i n_i =k}}^{k,[k/2], [k/3],\dots,1}  \;\;. 
\label{fortyfour}\end{eqnarray}
For the situation where the partition operator acts on unity, it yields $p(k)$ or the number of partitions summing to $k$.
Hence, $L_{P,k}[1]=p(k)$. The number of partitions summing to $k$ can also be obtained from the following generating 
function \cite{knu005}:
\begin{eqnarray}
P(z) = \prod_{m=1}^{\infty} \frac{1}{1-z^m} \equiv \sum_{k=0}^{\infty} p(k) z^k \;\;.
\label{fortyfoura}\end{eqnarray}  
Note the appearance of the equivalence symbol in the above result since the lhs can become divergent. This is because the  
lhs has been treated as an infinite product of regularised values of the geometric series in obtaining the power series on 
the rhs. This will become clearer in Sec.\ 6, where we shall also extend the above result. 

By applying Theorem\ 1 to the simple case where $p_0 \!=\!0$, $p_k \!=\! b^k$, $h(z) \!=\!1$, $a \!=\!1$, $q_k \!=\!1$ 
and $y \!=\! z$, we find that the quotient of the composite functions becomes
\begin{eqnarray}
g(f(z)) = 1+ \frac{bz}{1-2 bz} \equiv 1 + \frac{1}{2} \sum_{k=0}^{\infty} (2bz)^{k+1} \;\;.
\label{fortyfourb}\end{eqnarray}
Hence, via Eq.\ (\ref{five}) we arrive at
\begin{eqnarray}
L_{P,k}\Bigl[ N_k! \prod_{i=1}^k \frac{1}{n_i!} \Bigr]= 2^{k-1}  \;\;,
\label{fortyfourc}\end{eqnarray}
where, as before, $N_k \!=\! \sum_{i=1}^{k} n_i$. The $k$ subscript has been introduced here for the first time as it will 
become apparent that we shall need to sum the $n_i$ to different limits shortly. If we choose $p_k \!=\! (-b)^k$ and 
$q_k \!=\! (-1)^k$ instead, then following the same procedure we obtain 
\begin{eqnarray}
L_{P,k}\Bigl[ (-1)^{N_k} N_k! \prod_{i=1}^k \frac{1}{n_i!} \Bigr]= 0  \;\;,
\label{fortyfourd}\end{eqnarray}
where $k \! \geq \! 2$. The above is an interesting result where the partitions with an even number of elements are 
cancelled by those with an odd number of elements according to the frequencies of the elements. 

Eq.\ (\ref{fortyfourd}) is not the only instance where the sum over all partitions vanishes. For example, consider
the application of Theorem\ 1 to the function $f(z) =\exp(a\ln(1+z))$ or $f(z)=(1+z)^a$. Here, the coefficients of
the inner series are given by $p_0 \!=\!0$ and $p_k \!=\! (-1)^{k+1}/k$ for $k \!\geq\! 1$, while the coefficients
of the outer series are given by $q_k \!=\!1/k!$. Then from Theorem\ 1 we obtain
\begin{eqnarray}
D_k= (-1)^k  L_{P,k} \Bigl[ (-1)^{N_k} a^{N_k} \prod_{i=1}^k \frac{1}{i^{n_i}\, n_i!}\Bigr] \;\;.
\label{fortyfoure}\end{eqnarray}
In Ref.\ \cite{kow09a} it is shown that $f(z)$ represents the regularised value of the binomial series. That is,
for all values of $a$ and $z$, we have
\begin{eqnarray}
\sum_{k=0}^{\infty} \frac{ \Gamma(k-a)}{\Gamma(-a) \, k!} \; (-z)^k  \equiv (1+z)^a \;\;.
\label{fortyfoure1}\end{eqnarray}
For $|z| \!<\! 1$, the series is absolutely convergent and we can replace the equivalence symbol by an equals sign.
Since the $D_k$ are the coefficients of the power series expansion in $z$, they are equal to the coefficients
in the above result. Then we find that
\begin{eqnarray}
L_{P,k} \Bigl[ (-1)^{N_k} a^{N_k} \prod_{i=1}^k \frac{1}{i^{n_i}\, n_i!}\Bigr] = \frac{\Gamma(k-a)}{\Gamma(-a) \,k!} \;\;.
\label{fortyfoure1a}\end{eqnarray} 

The results in Theorem\ 1 are actually more general than the above result. As a consequence, for $p_0 \!=\! 0$ we 
arrive at 
\begin{eqnarray}
D_k= L_{P,k} \Bigl[ q_N \,a^{N_k} \, N_k! \, D_0^{-N_k} \prod_{i=1}^k \frac{p_i^{n_{i}}}{n_i!}\Bigr] \;\;,
\label{fortyfoure2}\end{eqnarray}
and
\begin{eqnarray}
E_k= L_{P,k} \Bigl[ (-1)^{N_k}  N_k! \, D_0^{-N_k} \prod_{i=1}^k \frac{D_i^{n_{i}}}{i^{n_i}\, n_i!}\Bigr] \;\;.
\label{fortyfoure3}\end{eqnarray}

When $a \!=\!1$, we have $f(z) \!=\! 1 \!+\! z$, which, in turn, means that $D_k \!=\!0$ for $k \!>\! 1$. Then we find 
that
\begin{eqnarray}
L_{P,k} \Bigl[ (-1)^{N_k} \prod_{i=1}^k \frac{1}{i^{n_i}\, n_i!}\Bigr] =0 \;\;,
\label{fortyfourf}\end{eqnarray}
for $k \!>\!1$. On the other hand, when $a \!=\! -1$, we have $f(z) \!=\! 1/(1+z)$, which represents the regularised value 
for the geometric series. Since the coefficients of the latter series are equal to $(-1)^k$, Eq.\ (\ref{fortyfoure})
reduces to    
\begin{eqnarray}
L_{P,k} \Bigl[ \prod_{i=1}^k \frac{1}{i^{n_i}\, n_i!}\Bigr] =1 \;\;.
\label{fortyfourg}\end{eqnarray}
More importantly, the above results can first be generalised by letting $a \!=\! l$, where $l$ is an arbitrary integer.
Then $f(z) \!=\! (1 \!+\! z)^l$, whose coefficients courtesy of the binomial theorem equal $\binom{l}{k}$ for 
$k \!\leq \! l$ and vanish for the remaining values of $l$. As a result, Eq.\ (\ref{fortyfoure}) yields
\begin{eqnarray}
L_{P,k} \Bigl[ (-l)^{N_k} \prod_{i=1}^k \frac{1}{i^{n_i}\, n_i!}\Bigr] = \begin{cases} 
0 \quad , & \quad  k > l \quad, \\ 
(-1)^k \binom{l}{k} \quad,  & \quad k \leq l \quad. 
\end{cases} 
\label{fortyfourh}\end{eqnarray}

If $-a$ is replaced by $\alpha$ in Eq.\ (\ref{fortyfoure1a}), then the equation reduces to 
\begin{eqnarray}
k! \;  L_{P,k} \Bigl[ \alpha^{N_k} \prod_{i=1}^k \frac{1}{i^{n_i}\, n_i!}\Bigr] = (\alpha)_k \;\;.
\label{fortyfouri}\end{eqnarray}
According to Chs.\ 24 and 18 of Refs.\ \cite{abr65} and \cite{spa87} respectively, the Pochhammer polynomials can be 
written as 
\begin{eqnarray}
 (\alpha)_k = (-1)^k \sum_{j=0}^k (-1)^j S^{(j)}_k \, \alpha^j \;\;,
\label{fortyfourj}\end{eqnarray}
where $S^{(j)}_k$ are known as the Stirling numbers of the first kind and satisfy
\begin{eqnarray}
S^{(j)}_{k+1} = S^{(j-1)}_k- k S^{(j)}_k \;\;.
\label{fortyfourk}\end{eqnarray}

We can proceed further with Eq.\ (\ref{fortyfouri}) by introducing a new operator that only considers a fixed number
of elements in the partitions. By setting this number equal to $j$, we can define the operator for a 
fixed number of elements as 
\begin{eqnarray}
L_{P,k}^{j}[\cdot] =  \sum_{\scriptstyle n_1, n_2,n_3, \dots,n_{k}=0 
\atop{\sum_{i=1}^{k}i n_i =k\;\;, \;\; \sum _{i=1}^k n_i=j}}^{k,[k/2], [k/3],\dots,1}  \;\;. 
\label{fortyfourl}\end{eqnarray}
Moreover, the above operator is related to the partition operator by
\begin{eqnarray}
L_{P,k}[\cdot] = \sum_{j=1}^{k} L_{P,k}^{j}[\cdot] \;\;.  
\label{fortyfourla}\end{eqnarray}
As mentioned in Sec.\ 2, the number of partitions of $k$ with exactly $j$ parts is denoted by 
$\begin{vmatrix} k \\ j\end{vmatrix}$. This means that
\begin{eqnarray}
L_{P,k}^{j}[1] =  \begin{vmatrix} k \\ j\end{vmatrix} \;\;,
\label{fortyfourm}\end{eqnarray}
while the recurrence relation given by Eq.\ (\ref{one}) becomes
\begin{eqnarray}
L_{P,k}^{j}[1] =  L_{P,k-1}^{j-1}[1] + L_{P,k-j}^{j}[1] \;\;. 
\label{fortyfouro}\end{eqnarray}
Furthermore, introducing the new operator into Eq.\ (\ref{fortyfouri}) with the rhs replaced with the aid of Eq.\
(\ref{fortyfourj}) yields
\begin{eqnarray}
L_{P,k}^{j}\Bigl[\prod_{i=1}^k \frac{1}{i^{n_i}\, n_i!} \Bigr] =  \frac{(-1)^{j+k}}{k!} \; S^{(j)}_k \;\;.
\label{fortyfourp}\end{eqnarray}
From the recurrence relation given by Eq.\ (\ref{fortyfourk}), we obtain
\begin{eqnarray}
(k+1)L_{P,k+1}^{j}\Bigl[ \prod_{i=1}^{k+1} \frac{1}{i^{n_i}\, n_i!}\Bigr] =  L_{P,k}^{j-1}\Bigl[ \prod_{i=1}^k
\frac{1}{i^{n_i}\,n_i!}\Bigr] +k  L_{P,k}^{j}\Bigl[ \prod_{i=1}^k \frac{1}{i^{n_i}\,i!}\Bigr] \;\;. 
\label{fortyfourp1}\end{eqnarray}

Specific results for the Stirling numbers of the first kind when $j$ is relatively small have been derived in Refs.\ 
\cite{kow10} and \cite{com74}. For example, these references give $S^{(2)}_k \!=\! (-1)^{k} \Gamma(k) H_1(k)$, where 
$H_1(k) =\sum_{j=1}^{k-1} 1/j$. Therefore, we find that
\begin{eqnarray}
L_{P,k}^{2}\Bigl[\prod_{i=1}^k\frac{1}{i^{n_i}\, n_i!} \Bigr] =  \frac{1}{k} \; H_1(k) \;\;.
\label{fortyfourq}\end{eqnarray}
Values for the Stirling numbers of the first kind are presented for $j$ close to $k$ in the appendix
of Ref.\ \cite{kow09}. For example, when $j \!=\! k \!-\! 1$, $S^{(k-1)}_k \!=\!- \binom{k}{2}$. This, in turn, 
leads to
\begin{eqnarray}
L_{P,k}^{k-1}\Bigl[\prod_{i=1}^k\frac{1}{i^{n_i}\, n_i!} \Bigr] =  \frac{(-1)^k}{k!} \; \binom{k}{2} \;\;.
\label{fortyfourr}\end{eqnarray}

Another fundamental result can be obtained by applying Theorem\ 1 to $\exp(-x)$. By writing the function as $1/\exp(x)$,
we see that the coefficients of the inner series, viz.\ $p_k$, are equal to $1/k!$ for $k \!\geq \! 1$, while  $p_0 \!=\! 0$.
Meanwhile, the outer series is given by the geometric series so that $q_k \!=\! (-1)^k$. The coefficients of the power
series for $\exp(-x)$ are $(-1)^k/k!$, which are also equal to the $D_k$. Therefore, according to Theorem\ 1 we have
\begin{eqnarray}
L_{P,k}^{k}\Bigl[N_k!\prod_{i=1}^k\frac{1}{(i!)^{n_i}\, n_i!} \Bigr] =  \frac{(-1)^k}{k!} \;\;.
\label{fortyfours}\end{eqnarray}
In the above result the $i \!=\! k$ term in the product is simply the result on the rhs. Therefore, it can be simplified
to
\begin{eqnarray}
L_{P,k}^{k}\Bigl[N_k!\prod_{i=1}^{k-1}\frac{1}{(i!)^{n_i}\, n_i!} \Bigr] =  \frac{2(-1)^k}{k!} \;\;.
\label{fortyfourt}\end{eqnarray}
Note in the above result that even though $n_k \!=\!0$, the constraint in the partition operator still applies to $k$.  

In Ref.\ \cite{kow11} we found that the cosecant numbers denoted by $c_k$ were the coefficients generated when the
partition method for a power series expansion was applied to $s \csc s$. This means that the method was basically
applied to
\begin{eqnarray}
s \csc s =\frac{1} {1 -s^2/3!+ s^4/5! -s^6/7! + \cdots} \;\;.
\label{fortyfive}\end{eqnarray}
By applying Theorem\ 1 to this example, in which $h(z) \!=\! z$, we see that $y \!=\! s^2$, $a \!=\!1$, 
$p_k = (-1)^{k+1}/(2k+1)!$ with $p_0 \!=\!0$ and $q_k \!=\! 1$ since the outer series corresponds to the geometric 
series. Therefore, from Theorem\ 1 we obtain 
\begin{eqnarray}
s\csc s \equiv \sum_{k=0}^{\infty} c_k s^{2k} \;\;,
\label{fortysix}\end{eqnarray}
where according to Eq.\ (\ref{five}),
\begin{eqnarray}
(-1)^k\, c_k=  L_{P,k}\left[ N_k! \prod_{i=1}^{k} \left(\frac{-1}{(2i+1)!}\right)^{n_i} \;\frac{1}{n_i!}  \right] \;\;,
\label{fortyseven}\end{eqnarray}
and $N_k \!=\! \sum_{i=1}^k n_i$. It should be noted that $\prod_{i=1}^{k} (-1)^{i n_i +n_i} = (-1)^{k+N_k}$, although we shall retain 
the phase factor of $(-1)^{n_i}$ in order to observe a remarkable correspondence arising from the inversion of Equivalence\ 
(\ref{fortysix}). It was also found in Ref.\ \cite{kow11} that the radius of absolute convergence for the power series expansion 
in Equivalence\ (\ref{fortysix}) was $\pi$, while the cosecant numbers were seen to be rapidly decreasing positive fractions 
given by Eq.\ (\ref{thirtythreea4}). Hence, Eq.\ (\ref{fortyseven}) represents a means of determining even integer values of 
the Riemann zeta function.

Now we invert the power series expansion in Equivalence\ (\ref{fortysix}) and apply Theorem\ 1 again, but 
in the case  we have  $p_k \!=\! -c_k$ with $p_0 \!=\!0$. The resulting power series expansion was found to yield the standard 
Taylor/Maclaurin power series for $\sin(s)/s$ or $\sum_{k=0}^{\infty} (-1)^k s^{2k}/(2k+1)!$, which is convergent for all values of 
$s$, despite the fact that Equivalence\ (\ref{fortysix}) has a radius of absolute convergence equal to $\pi$. This 
confirms the earlier remark concerning the fact that the resulting power series expansion obtained via Theorem\ 1 
can turn out to be convergent even though the inner or the outer series may in fact be divergent. Furthermore, 
from Eq.\ (\ref{five}) we obtain
\begin{eqnarray}
\frac{(-1)^k }{(2k+1)!} =  L_{P,k}\left[ N_k! \prod_{i=1}^{k} (-c_i)^{n_i} \;\frac{1}{n_i!}  \right] \;\;.
\label{fortyeight}\end{eqnarray} 

It should also be mentioned that in Ref.\ \cite{kow11} numerous recurrence relations were derived for the cosecant
numbers. One of these is
\begin{eqnarray}
\sum_{j=0}^{k-1} \frac{(-1)^{k-j-1}}{(2k-2j+1)!}\; c_j = c_k \;\;,
\label{fortyeighta}\end{eqnarray}
where $c_0 \!=\! 1$.  If we introduce Eq.\ (\ref{fortyseven}) into the Eq.\ (\ref{fortyeighta}), then we obtain the interesting 
result of
\begin{align}
& \sum_{j=0}^{k-1} \frac{1}{(2k-2j+1)!}  \; L_{P,j}\left[ (-1)^{N_j-1} \,N_j! \prod_{i=1}^{j} \left(\frac{1}{(2i+1)!}\right)^{n_i} 
\;\frac{1}{n_i!} \right] 
\nonumber\\
& = \;\; L_{P,k}\left[ (-1)^{N_k} \, N_k! \prod_{i=1}^{k} \left(\frac{1}{(2i+1)!}\right)^{n_i} \;\frac{1}{n_i!}  \right] \;\;.
\label{fortyeightb}\end{align}
Similar results to the above can be obtained by considering the other recurrence relations. Furthermore, in the same reference
it was found that just as the Bernoulli numbers give rise to Bernoulli polynomials, the cosecant numbers give rise to their
own polynomials, which are related to the former. In particular, the value at unity was found to be given by
\begin{eqnarray}
c_k(1) = (-1)^k\,2^{2k} L_{P,2k} \left[ (-1)^{N_{2k}} \,N_{2k}! \prod_{i=1}^{2k} \left( \frac{1}{(i+1)!}\right)^{\! n_i}
\frac{1}{n_i!} \right] \;\;,
\label{fortyeightc}\end{eqnarray}  
where $N_{2k} = \sum_{i=1}^{2k} n_i$ in the above result. Consequently, we see the reason for the introduction of the 
$k$-subscript to $N$ in Theorem\ 1. The value of the cosecant polynomials at unity was found to equal
\begin{eqnarray}
c_k(1) = \frac{c_k}{2^{1-2k}-1} \;\;.
\label{fortyeightd}\end{eqnarray}
Hence, by introducing Eq.\ (\ref{fortyseven}) into the above result and then equating it to Eq.\ (\ref{fortyeightc}), we
arrive at 
\begin{align}
&  L_{P,2k}\left[ (-1)^{N_{2k}-1} \,N_{2k}! \prod_{i=1}^{2k} \left(\frac{1}{(i+1)!}\right)^{n_i} 
\;\frac{1}{n_i!} \right] 
\nonumber\\
& = \;\; \frac{1}{(2^{2k}-2)} \; L_{P,k}\left[ (-1)^{N_k} \, N_k! \prod_{i=1}^{k} \left(\frac{1}{(2i+1)!}\right)^{n_i} \;\frac{1}{n_i!}  
\right] \;\;.
\label{fortyeighte}\end{align}
Consequently, we do not need to consider all the partitions up to $2k$ to determine the sum on the lhs, which is a 
significant reduction in computational effort. 

In Ref.\ \cite{kow11} another infinite set of related numbers denoted by $d_k$ and known as the secant numbers were obtained 
when the partition method for a power series expansion was applied to $\sec s$. Specifically, the method was applied to
\begin{eqnarray}
\sec s =\frac{1} {1 -s^2/2!+ s^4/4! -s^6/6! + \cdots} \;\;.
\label{fortynine}\end{eqnarray}
By applying Theorem\ 1 to the above result we have $y \!=\! s^2$, $p_0 \!=\!0$, $a \!=\! 1$ and $q_k \!=\! 1$ as before,
but on this occasion, $p_k \!=\!(-1)^{k+1}/(2k)!$. The resulting power series expansion, which can be expressed as
\begin{eqnarray}
\sec s \equiv \sum_{k=0}^{\infty} d_k\, s^{2k} \;\;,
\label{fortyninea}\end{eqnarray} 
where the coefficients $d_k$ from Eq. (\ref{five}) are given by
\begin{eqnarray}
(-1)^k d_k=  L_{P,k} \left[  N_k! \prod_{i=1}^{k} \left( -\frac{1}{(2i)!} \right)^{\! n_i} \frac{1}{n_i!}  \right] \;\;,
\label{fifty}\end{eqnarray}
and $N_k$ is the same sum over the frequencies as before. Equivalence\ (\ref{fortynine}) was found to possess a narrower 
radius of absolute convergence compared with the power series expansion for cosecant, viz. $\pi/2$ as opposed to $\pi$, 
while the $d_k$ or secant numbers were found to be not as rapidly decreasing fractions as their cosecant counterparts. In 
addition, instead of being related to the Riemann zeta function as the cosecant numbers are, they were found to be related 
to the Hurwitz zeta function by
\begin{eqnarray}
d_k = \frac{2^{2k+2}}{\pi^{2k+1}} \Bigl( \zeta(2k+1,1/4) -\zeta(2k+1,3/4) \Bigr) \;\;.
\label{fiftyone}\end{eqnarray}
The bracketed expression can also be written as $\sum_{j=1}^{\infty} (-1)^{j+1}/(2j-1)^{2k+1}$. To 
invert the analysis, either we can apply Theorem\ 1 to the power series expansion in Equivalence\ (\ref{fortyninea}) or we can
go directly to Equivalence\ (\ref{six}). In the latter case we replace $D_i$ by $d_i$ in Eq.\ (\ref{seven}), while the lhs of 
Equivalence\ (\ref{six}) equals $\cos s$, which we replace by its power series expansion. Then we find that the coefficients 
$E_k$ in Eq.\ (\ref{seven}) equal $(-1)^{k}/(2k)!$. Hence, we arrive at
\begin{eqnarray}
\frac{(-1)^k}{(2k)!} =  L_{P,k} \left[  N_k! \prod_{i=1}^{k} (-d_i)^{n_i} \;\frac{1}{n_i!}  \right] \;\;.
\label{fiftytwo}\end{eqnarray}

The secant numbers were also found to obey recurrence relations, although not as many as their cosecant counterparts. 
Nevertheless, an analogue of Eq.\ (\ref{fortyeighta}) was obtained, which is given by 
\begin{eqnarray}
\sum_{j=0}^{k-1} \frac{(-1)^{k-j-1}}{(2k-2j)!}\; d_j = d_k \;\;,
\label{fiftytwo1}\end{eqnarray}
with $d_0 \!=\! 1$. Introducing Eq.\ (\ref{fiftytwo}) into the above result yields
\begin{eqnarray}
& \sum_{j=0}^{k-1} \frac{1}{(2k-2j)!}  \; L_{P,j}\left[ (-1)^{N_j-1} \,N_j! \prod_{i=1}^{j} \left(\frac{1}{(2i)!}\right)^{n_i} 
\;\frac{1}{n_i!} \right] 
\nonumber\\
& = \;\; L_{P,k}\left[ (-1)^{N_k} \, N_k! \prod_{i=1}^{k} \left(\frac{1}{(2i)!}\right)^{n_i} \;\frac{1}{n_i!}  \right] \;\;.
\label{fiftytwo2}\end{eqnarray}
This result is virtually identical to Eq.\ (\ref{fortyeightb}) except that there are no ``+1's" in the denominators of
the above equation. 

More sophisticated results involving both the secant and cosecant numbers can also be derived. From No.\ 1.518(2) of Ref.\
\cite{gra94} we have
\begin{eqnarray}
\ln \sec(\pi z) = \sum_{k=1}^{\infty} \frac{(2^{2k}-1)}{k\, (2k)!} \;2^{2k-1} |B_{2k}| (\pi z)^{2k} \;\;,
\label{fiftytwoa}\end{eqnarray}
which is absolutely convergent for $|z| \!<\! 1/2$. In Ref.\ \cite{kow11} it is shown that the cosecant numbers are 
related to the Bernoulli numbers by
\begin{eqnarray}
c_k =\frac{(-1)^{k+1}}{(2k)!}\; \left( 2^{2k}-2\right) B_{2k} \;\;.
\label{fiftytwob}\end{eqnarray}     
Moreover, we can write the lhs of Eq.\ (\ref{fiftytwoa}) as 
\begin{eqnarray}
\ln \sec (\pi z) = \ln \Bigl( 1+ \sum_{k=1}^{\infty} d_k (\pi z)^{2k} \Bigr) \;\;.
\label{fiftytwoc}\end{eqnarray}
We now apply Theorem\  1 to the rhs of the above result. This means that we expand the logarithm in terms of its
Maclaurin series expansion, in which case $q_k \!=\! (-1)^{k+1}/k$, while $y \!=\! z^2$ and $p_k \!=\! d_k \pi^{2k}$.
We then equate the resulting power series expansion to like powers of $z^2$ or $y$ in Eq.\ (\ref{fiftytwoa}). By 
substituting Eq.\ (\ref{fiftytwob}) we replace the Bernoulli numbers by the secant numbers, thereby obtaining
\begin{eqnarray}
L_{P,k} \Bigl[ (-1)^{N+1} \, (N-1)! \prod_{i=1}^k \frac{d_i^{n_i}}{n_i!} \Bigr]= \frac{1}{2k} \, \frac{(2^{2k}-1)}
{(1-2^{1-2k})}\; c_k \;\;.
\label{fiftytwod}\end{eqnarray} 
Once again, we see the partition operator acting on another strange argument to yield an interesting finite quantity for
all values of $k$.

In Ref.\ \cite{kow10} the partition method for a power series expansion is applied to the reciprocal of the logarithmic
function $\ln(1 \!+\! z)$. There a power series expansion is obtained in terms of special coefficients $A_k$, which
are referred to as the reciprocal logarithm numbers. On p.\ 138 of Ref.\ \cite{ape008} these numbers are referred to as 
the Gregory or Cauchy numbers when their modulus is taken. That is, the following result is obtained: 
\begin{eqnarray}
\frac{1}{\ln(1+z)} \equiv \sum_{k=0}^{\infty} A_k \,z^{k-1} \;\;.
\label{fiftythree}\end{eqnarray}
The reciprocal logarithm numbers are found to be oscillating fractions, which are more slowly converging to zero than 
either the cosecant or secant numbers. Moreover, they are given by
\begin{eqnarray}
A_k = \frac{(-1)^k}{k!} \int_0^1 dt \; \frac{\Gamma(k+t-1)}{\Gamma(t-1)} \;\;.
\label{fiftythreea}\end{eqnarray}  
Expressing $\ln (1 \!+\! z)$ in terms of its Maclaurin series as in the previous example, which is absolutely 
convergent only for $|z| \!<\! 1$, we are in a position to apply Theorem\ 1 to $z/\ln(1 \!+\! z)$. In this case, 
$h(z) \!=\! 1/z$ and $f(z) \!=\! \ln(1 \!+\! z)$. Then the coefficients of the inner series, viz.\ $p_k$, are equal 
to $(-1)^{k+1}/(k+1)$ for $k \!>\! 0$, while for $k \!=\!0$, $p_0 \!=\! 0$. The resulting denominator can be regarded 
as the regularised value of the geometric series, which means that the $q_k$ are equal to $(-1)^k$ as in the preceding 
examples. Again, $a \!=\! 1$. Hence, according to Eq.\ (\ref{five}), the reciprocal logarithm numbers can be written as
\begin{eqnarray}
(-1)^k \, A_k =  L_{P,k} \left[  N! \prod_{i=1}^{k} \left(-\frac{1}{i+1} \right)^{n_i} \;\frac{1}{n_i!}  \right] \;\;.
\label{fiftyfour}\end{eqnarray}
The inverse of this result is obtained by putting $D_0 \!=\! 1$ and $D_k  \!=\! A_k$ in Eq.\ (\ref{seven}), while 
the $E_k$ equal the coefficients in the Maclaurin series for $\ln(1 \!+\! z)$, i.e.\ $(-1)^{k+1}/(k\!+\!1)$. Then we find 
that
\begin{eqnarray}
\frac{(-1)^{k+1}}{k+1} =  L_{P,k} \left[  N! \prod_{i=1}^{k} \frac{(-A_i)^{n_i}}{n_i!}  \right] \;\;.
\label{fiftyfive}\end{eqnarray}

As an aside, in Ref.\ \cite{kow10}, Euler's constant is derived in terms of an infinite series involving the reciprocal
logarithm numbers, where it is also referred to as Hurst's formula. Since the publication of Ref.\ \cite{kow10}, it has been
revealed that the formula was independently discovered by Kluyver \cite{klu24}. By using Eq.\ (\ref{fiftyfour}), we can express
Euler's constant as
\begin{eqnarray}
\gamma= - \sum_{k=1}^{\infty} \frac{1}{k} \; L_{P,k} \left[  (-1)^N \, N! \prod_{i=1}^{k} \left(\frac{1}{i+1} \right)^{n_i} 
\;\frac{1}{n_i!}  \right] \;\;.
\label{fiftyfivea}\end{eqnarray}
Alternatively, we can introduce Eq.\ (\ref{fiftyfive}) into Hurst's formula, which yields
\begin{eqnarray}
\gamma= - \sum_{k=1}^{\infty} A_k \; L_{P,k} \left[  (-1)^N \, N! \prod_{i=1}^{k} \frac{A_i^{n_i}}{n_i!}  \right] \;\;.
\label{fiftyfiveb}\end{eqnarray}
Moreover, Euler's constant is not the only result found in Ref.\ \cite{kow10} that can be expressed as an infinite sum over the
reciprocal logarithm numbers. For example, $\ln 2$ can be expressed as a similar sum to Hurst's formula. Therefore,
with the aid of Eq.\ (\ref{fiftyfour}) we find that
\begin{eqnarray}
\ln 2= \sum_{k=1}^{\infty} \frac{1}{k+1} \; L_{P,k} \left[ (-1)^N \,  N! \prod_{i=1}^{k} \left(\frac{1}{i+1} \right)^{n_i} 
\;\frac{1}{n_i!}  \right] \;\;.
\label{fiftyfivec}\end{eqnarray}

With the aid of Corollary\ 1 to Theorem\ 1 we can generalise the preceding examples to where the generating functions are raised 
to an arbitrary power $\rho$.  For example, the quotient in Eq.\ (\ref{fortyfourb}) raised to an arbitrary power $\rho$ becomes
\begin{eqnarray}
g(f(z))^{\rho} = \Bigl(\frac{1-bz}{1-2 bz}\Bigr)^{\rho} \equiv \sum_{k=0}^{\infty} (bz)^{k} \sum_{j=0}^k
\frac{\Gamma(j-\rho)}{\Gamma(-\rho)\, j!} \; \frac{2^{k-j}\,\Gamma(k-j+\rho)}{\Gamma(\rho)\, (k-j)!}\;\;,
\label{fiftysix}\end{eqnarray}
where we have used the regularised value for the binomial series in Lemma\ 1. On the other hand, according to Equivalence\
(\ref{thirtyeight}), the lhs of the above result can be written as
\begin{eqnarray}
g(f(z))^{\rho} \equiv  1+ \frac{1}{2} \sum_{k=0}^{\infty} (2bz)^{k} \;\;.
\label{fiftyseven}\end{eqnarray}
Hence, $D_0 \!=\! 1$ and $D_k \!=\! 2^{k-1}b^k$ for $k \!\geq \!1$, while $y \!=\! z$. From Eq.\ (\ref{thirtyfive}) we
obtain
\begin{eqnarray}
D_k(\rho) = (2b)^k L_{P,k} \Bigl[(-1/2)^N \,(-\rho)_N \prod_{i=1}^k \frac{1}{n_i!} \Bigr] \;\;.
\label{fiftyeight}\end{eqnarray}
Equating like powers of $z$ on the rhs's of the preceding equivalences yields
\begin{eqnarray}
L_{P,k} \Bigl[(-1/2)^N (-\rho)_N \prod_{i=1}^k \frac{1}{n_i!} \Bigr] = \sum_{j=0}^k \frac{\Gamma(j-\rho)}{\Gamma(-\rho)\, j!} 
\; \frac{\Gamma(k-j+\rho)}{2^j \,\Gamma(\rho)\, (k-j)!} \;\;.
\label{fiftynine}\end{eqnarray}
 
If Equivalence\ (\ref{fortysix}) is taken to the arbitrary power of $\rho$, then the $D_k$ of Equivalence\ (\ref{thirtyeight})
become the cosecant numbers, i.e.\ $D_i \!=\! c_i$. If we denote $D_k(\rho)$ by $c_{\rho,k}$, then we find that according to
Eq.\ (\ref{thirtyfive}), these generalised cosecant numbers are given by
\begin{eqnarray}
c_{\rho,k} = L_{P,k} \Bigl[(-1)^N (-\rho)_N \prod_{i=1}^k \frac{c_i^{n_i}}{n_i!} \Bigr] \;\;.
\label{sixty}\end{eqnarray}
Alternatively, the original equation can be expressed as 
\begin{eqnarray}  
s^{\rho} \csc^{\rho}s = \Bigl( 1 +\sum_{k=1}^{\infty} \frac{(-1)^k x^{2k+1}}{(2k+1)!} \Bigr)^{-\rho} \;\;. 
\label{sixtya}\end{eqnarray}
Now the $D_k$ in Equivalence\ (\ref{thirtyeight}) are equal to $(-1)^k/(2k+1)!$, while $\rho$ has changed sign. Thus, the
generalised cosecant numbers can be written as
\begin{eqnarray}
c_{\rho,k} = (-1)^k L_{P,k} \Bigl[(-1)^N \,(\rho)_N \prod_{i=1}^k \Bigl(\frac{1}{(2i+1)!}\Bigr)^{n_i} \frac{1}{n_i!} \Bigr] \;\;.
\label{sixtyone}\end{eqnarray}
This form for the generalised cosecant numbers has been recently been employed in evaluating a finite sum of inverse powers of 
cosines given as $S_{m,v}=(\pi/2m)^{2v} \sum_{k=1}^{m-1} \cos(k\pi/m)^{-2v}$ in Ref.\ \cite{kow11a}. 

In a similar manner we can generalise the secant numbers. By taking the $\rho$-th power of Equivalence\ (\ref{fortyninea}), we see that
the $D_k$ in Equivalence\ (\ref{thirtyeight}) are equal to  $d_i$. By denoting the generalised secant numbers as $d_{\rho,k}$, we
find via Eq.\ (\ref{thirtyfive}) that
\begin{eqnarray}
d_{\rho,k} = L_{P,k} \Bigl[(-1)^N (-\rho)_N \prod_{i=1}^k \frac{d_i^{n_i}}{n_i!} \Bigr] \;\;,
\label{sixtytwo}\end{eqnarray}
while taking the $\rho$-th power of Eq.\ (\ref{fortynine}) yields
\begin{eqnarray}
d_{\rho,k} = (-1)^k L_{P,k} \Bigl[(-1)^N \,(\rho)_N \prod_{i=1}^k \Bigl(\frac{1}{(2i)!}\Bigr)^{n_i} \frac{1}{n_i!} \Bigr] \;\;.
\label{sixtythree}\end{eqnarray}
The above result appears as Eq.\ (290) in Ref.\ \cite{kow11}.

To generalise the reciprocal logarithm numbers, we take the $\rho$-th power of Equivalence\ (\ref{fiftythreea}). Then the
$D_k$ in Equivalence\ (\ref{thirtyeight}) are equal to $A_k$.  Denoting the generalised reciprocal numbers by $A_k(-\rho)$,
we find via Eq.\ (\ref{thirtyfive}) that they are given by
\begin{eqnarray}
A_{k}(-\rho) = L_{P,k} \Bigl[(-1)^N \,(-\rho)_N \prod_{i=1}^k \frac{A_i^{n_i}}{n_i!} \Bigr] \;\;.
\label{sixtyfour}\end{eqnarray}
Alternatively, the generalised reciprocal logarithm numbers can be determined by taking the $\rho$-th power of the 
Maclaurin series for $\ln(1+z)$ , which is
\begin{eqnarray}
\ln(1+z) \equiv  \sum_{k=1}^{\infty} \frac{(-1)^{k+1}}{k} \; z^k \;\;.
\label{sixtyfive}\end{eqnarray} 
As explained in Ref.\ \cite{kow09}, the above result is absolutely convergent for $\Re\, z \!<\! -1$, in
which case the equivalence symbol can be replaced by an equals sign. With regard to Equivalence\ (\ref{thirtyeight}) 
the $D_k$ are now equal to $(-1)^{k+1}/k$, while $y \!=\! z$. Consequently, Eq.\ (\ref{thirtyfive}) yields
\begin{eqnarray}
A_{k}(\rho) = (-1)^k \,L_{P,k} \Bigl[(-1)^N \,(\rho)_N \prod_{i=1}^k \Bigl( \frac{1}{i+1}\Bigr)^{n_i} 
\frac{1}{n_i!} \Bigr] \;\;.
\label{sixtysix}\end{eqnarray}
Eq.\ (\ref{sixtysix}) appears as Eq.\ (118) in Ref.\ \cite{kow09}.

At the end of Sec.\ 2 an alternative algorithm was given for accessing the partitions via the tree diagram in Fig.\ \ref{figone}.
As a result, we can present a new formulation of the partition method for a power series expansion. Before doing so, however, we
need to amend the definition of the partition operator $L_{P,k}$. The amendment is necessary so that we can sum along each of the 
branches emanating from the seed number. For example, the entire sub-tree from \{1,5\} represents the tree diagram for partitions
summing to 5, but to obtain those summing to 6, we need to increment $n_1$ in all the partitions summing to 5 by one. This, of course,
affects all the contributions to the coefficients. In addition, when we consider the partitions emanating from \{2,4\}, we need
to ensure that no partitions with unity will appear, while there should be no ones or twos for the partitions emanating from \{3,3\}
and so on. Therefore, we define the restricted partition operator as follows:
\begin{eqnarray}
L_{RP,k,i}\Bigl[ \cdot \Bigr]= \sum_{\scriptstyle n_i=1,n_{i+1},\dots,n_{k-i}=0 
\atop{\sum_{j=i}^{k-i} j n_j=k}}^{1+[(k-i)/i],[(k-i)/(i+1)], \dots,1} \;\;.
\label{sixtyseven}\end{eqnarray}    

There are major differences between the above operator and the partition operator as defined by Eq.\ (\ref{fortyfour}). The first is
that the sum begins at $n_i$ rather than at $n_1$. This is due to the fact that we need to exclude elements less than $i$ when
summing along the branches emanating from the seed number. The next difference is that $n_i$ begins at unity rather than zero 
as for the other elements, which accounts for the fact already a one element of $i$ has occurred in moving from the seed
number. Furthermore, as a result of separating the element $i$, the maximum element in the resulting partitions becomes 
$k \!-\! i$. Hence, the elements range from $i$ to $k-i$, which not only affects the number of summations, but also their upper 
limits. In addition, whilst the number of summations is restricted in the constraint, the value remains invariant, viz.\ $k$.
  
We are now in a position to implement the algorithm described at the end of Sec.\ 2. Basically, this entails expressing a result 
like Eq.\ (\ref{five}) in terms of the partition operator on the lhs and the sum of restricted partition operators on the rhs.
Therefore, we arrive at
\begin{align}
&L_{P,k}\Bigl[ q_{N_{1,k}} \, a^{N_{1,k}} \, N_{1,k}! \; \prod_{i=1}^k \frac{p_i^{n_i}}{n_i!} \Bigr] = q_1 \, a \, p_k 
\nonumber\\
& +  \sum_{j=1}^{[k/2]} L_{RP,k,j} \Bigl[ q_{N_{j,k-j}}\, a^{N_{j,k-j}}\, N_{j,k-j}! \; \prod_{i=j}^{k-j} \frac{p_i^{n_i}}{n_i!} 
\Bigr] \;\;,
\label{sixtyeight}\end{align}
where $N_{i,k-i} = \sum_{j=i}^{k-i} n_j$ and $N_{1,k} \!=\! N$. For $p_0 \!\neq\!0$ in Theorem\ 1, we use Eq.\ (\ref{four}) instead, 
which amounts to replacing $q_1 a p_k$ and $q_{N_{j,k-j}}\, a^{N_{j,k-j}}\, N_{j,k-j}!$ in the above result by $a \,F^{(k)}(ap_0)$ 
and $a^{N_{j,k-j}}\, F^{N_{j,k-j}} (ap_0)$ respectively.

To complete this section, we consider the situation where the quotient of the pseudo-composite functions in Theorem\ 1 yields
a function $r(y)$, which is infinitely differentiable. This results in the following corollary.\newline
{\bfseries Corollary 2 to Theorem\ 1}. If the functions $f(z)$ and $g(z)$ obey the same conditions as in Theorem\ 1 and
the quotient of the pseudo-composite functions $g_a \circ f$ and $h_a \circ f$ yields an infinitely differentiable 
function, $r(y)$, then for $p_0 \!\neq\!0$,
\begin{eqnarray}
r^{(k)}(0)= k! \;  L_{P,k}\left[ a^N F^{(N)} (a p_0)\prod_{i=1}^k \frac{p_i^{n_i}}{n_i!} \right] \;\;,
\label{sixtynine}\end{eqnarray} 
while for $p_0 \!=\! 0$, 
\begin{eqnarray}
r^{(k)}(0) = k! \; L_{P,k} \left[ q_N \, a^N N! \prod_{i=1}^k \frac{p_i^{n_i}}{n_i!} \right] \;\;. 
\label{seventy}\end{eqnarray}
Furthermore, if $r(0) \! \neq \! 0$, then inversion of the quotient yields 
\begin{eqnarray}
\left. \Bigl(\frac{1}{r(y)}\Bigr)^{(k)} \right|_{y=0}= k! \; \frac{E_k}{D_0} \;\;,
\label{seventyone}\end{eqnarray}
where the $E_k$ are given by Eq.\ (\ref{seven}).

{\bfseries Proof}. From the proof of Theorem\ 1, we know that the ratio of $g_a \circ f$ over $h_a \circ f$ is equivalent
to the power series $\sum_{k=0}^{\infty} D_k y^k$, where the coefficients $D_k$ are given by either Eq. (\ref{four}) for
$p_0$ non-zero, or Eq.\ (\ref{five}) when $p_0 \!=\!0 $. Since the ratio of the pseudo-composite functions yields an infinitely
differentiable function $r(y)$ according to the corollary, the ratio can also be expressed as a Taylor/Maclaurin series
given by
\begin{eqnarray}
\frac{g_a \circ f}{h_a \circ f} \equiv \sum_{k=0}^{\infty} r^{(k)}(0) \; \frac{y^k}{k!} \;\;,
\label{seventytwo}\end{eqnarray} 
where the superscript $(k)$ now denotes the $k$ times differentiation of $r(y)$ w.r.t.\ $y$. From Theorem\ 1 we also know
that the above quotient can be expressed as 
\begin{eqnarray}
\frac{g_a \circ f}{h_a \circ f} \equiv  \sum_{k=0}^{\infty} D_k \, y^k  \;\;,
\label{seventythree}\end{eqnarray}
where the coefficients $D_k$ are given by either Eq.\ (\ref{four}) when $p_0 \! \neq \! 0$ or Eq. (\ref{five}) for
$p_0 \!=\!0$. Since the regularised value is unique as described in Refs.\ \cite{kow09} and \cite{kow09a}-\cite{kow95}, 
the rhs's of the two preceding equivalence statements are equal to one another. Since $y$ is arbitrary in the resulting
equation, we can equate like powers or the coefficients of both power series, which yields Eq.\ (\ref{sixtynine}) or
(\ref{seventy}) depending upon the the value of $p_0$. 

If the quotient of the pseudo-composite functions is inverted, then it will equal $1/r(y)$. If $r(0) \!\neq \!0$, then 
the function $1/r(y)$ can also be expressed as a Taylor/Maclaurin series since $r(y)$ is infinitely differentiable. Therefore,
we have 
\begin{eqnarray}
\left. \frac{h_a \circ f}{g_a \circ f} \equiv  \sum_{k=0}^{\infty} \left( \frac{d^k}{dy^k} \frac{1}{r(y)} \right)
\right|_{y=0}  \!\!\! \frac{y^k}{k!}  \;\;.
\label{seventyfour}\end{eqnarray}
From Theorem\ 1 we know that the above quotient also represents the regularised value of a power series in $y$ whose 
coefficients are equal to $E_k/D_0$, while the $E_k$ are given by Eq.\ (\ref{seven}). Moreover, since $r(0) \!\neq \! 0$,
$D_0$ does not vanish. Since both power series have the same regularised value, they are equal to one another for the 
same reason as in the first part of the proof. Again, as $y$ is arbitrary, we can equate like powers of $y$, which results 
in Eq.\ (\ref{seventyone}). This completes the proof of the corollary.

So far, we have described the partition method for a power series expansion in terms of a novel discrete operator, which has 
been referred to as the partition operator and is denoted by $L_{P,k}[\cdot]$. At this stage the operator has been used to 
derive general results for the coefficients of the power series expansions obtained via Theorem\ 1. However, it is not always 
possible to derive general results dependent only upon $k$ for the coefficients $D_k$ and $E_k$ in Theorem\ 1. Often the 
coefficients become polynomials dependent upon another variable or even a function. For these cases we need to develop a 
program that can automatically evaluate specific coefficients, which is addressed in the following section.    

\section{Programming the Partition Method for a Power Series Expansion}
At the end of the previous section we indicated the need for developing a computer program to enable the evaluation of
the coefficients via the partition method for a power series expansion. As indicated earlier, such an approach can be 
developed by employing the BRCP algorithm of Sec.\ 2, but before doing so here, some remarks are necessary. Because the 
partition method for a power series expansion relies on evaluating the contribution due to each partition and the number 
of partitions $p(k)$ according to the Hardy-Rademacher-Ramanujan formula is $O(\exp(\pi \sqrt{2k/3})/4k\sqrt{3})$ 
\cite{yam07,knu005,wik11}, any computer program based on partitions as its input will ultimately become very slow. In 
fact, since all the partitions summing to the order of each power are involved, such a program represents a 
brute-force approach to deriving power series expansions. Nevertheless, determining power series expansions for orders 
up to 40 ($p(50)=37\,338$) should be achievable with most number-crunching computers around today. So, at least for 
intermediate values of the order $k$, programming the partition method is still of great benefit, particularly for 
intractable functions where it represents the only method we have of deriving a power series expansion.

With regard to very high orders it should be noted that the partition method does not actually use the partitions themselves. 
What the method requires is each element appearing in each partition and their frequency, which is referred to as the
multiplicity representation in Sec.\ 2. This information can be stored in external arrays which can be called upon when 
one wishes to determine the series expansion for different situations. Therefore, there is no need to repeat the process 
of generating the partitions when dealing with different problems. In addition, the contributions due to many of the 
partitions will often be negligible even by today's computing standards. In those cases the calculation of the coefficients 
can be simplified yielding extremely accurate approximations. In other cases it is possible to sum the contributions in 
classes or groups, thereby avoiding the necessity of processing each partition separately. This issue will be addressed
later in this work. Finally, by developing a programming approach to the partition method, we will be in a position to 
consider different problems in the theory of partitions such as the evaluation of partitions with specific elements 
including those with discrete elements, doubly-restricted partitions and the transposes of partitions. The significant 
issue here is that such problems only require small changes to the BRCP algorithm, whereas separate programs are required 
when other codes such as those presented in Sec.\ 2 are used to generate partitions. E.g., in order to determine the 
partitions with a fixed number of elements in them, Knuth presents another algorithm based on the 18-th century 
dissertation by C.F. Hindenburg on p.\ 38 of Ref.\ \cite{knu005}. As we shall see in Sec.\ 5, this problem can be solved by 
inserting a few lines into the BRCP algorithm.     

In the previous section the partition method for a power series expansion was described in terms of the partition operator,
$L_{P,k}[\cdot]$. There many general results involving this operator were given without requiring to evaluate the sum over 
partitions for specific values of $k$. Consequently, the partition operator can be regarded as an intricate abstract 
operator when compared with more well-known operators such as the differential operator. However, whilst taking the 
derivative of a function is also viewed as an abstract operation, we at least have an understanding of the process 
because we can always calculate the limit of Newton's difference quotient provided, of course, it exists. As a result of 
this understanding, general shorthand rules such as $dx^k/dx \!=\! kx^{k-1}$ have evolved. Yet, the opposite situation 
applies to the partition operator--- we have a few general results in the previous section and in Refs.\ \cite{kow10}-
\cite{kow11}, but we do not even have a means of applying the operator outside of those cases to evaluate the first few 
coefficients of the expansions given in Theorem\ 1. Therefore, we require an approach that will allow us to apply the 
partition operator for any value of $k$ to any situation that obeys the conditions in Theorem\ 1, even if it is no 
longer feasible to evaluate the coefficients for very large values of $k$. 

Now that we have indicated why it is necessary to program the partition method for a power series expansion, we 
turn to the issue of the programming languages required for the task. The first point to be noted is that if we choose a
standard high-level programming language like C/C++ or Fortran, then our results for the coefficients will inevitably 
become decimal numbers when they could be rational. Moreover, they will invariably be rounded off or worse still,
may only equal zero if significantly smaller than the precision allowed by the computing system. In addition, the coefficients 
need not be numerical as exemplified by the examples appearing after Corollary\ 1 to Theorem\ 1. All this means is that 
we require a mathematical software package such as Mathematica to retain either the rationality of the coefficients
or when applicable, their symbolic form. However, programming the BRCP algorithm in Mathematica with its bi-variate recursion 
is also formidable. Instead, the issue can be overcome by using the material in Sec.\ 2. Therefore, the best option is 
to combine the strengths of both C/C++ and Mathematica. Basically, this means that the initial program is to be written
in C/C++ so that the coefficients can be printed out in a symbolic form. Then these forms can be introduced into Mathematica 
where we can use either the integer arithmetic routines to evaluate the coefficients, thereby avoiding the round-off that 
occurs with floating point numbers or its symbolic routines to reduce all the terms generated by the C/C++ code into simple 
mathematical statements such as polynomials. 
    
The appendix presents the C/C++ program called ${\bf partmeth}$, which outputs in symbolic form the coefficients
$D_k$ and $E_k$ given in Theorem\ 1. Here, we are only concerned with the case of $p_0 \!=\!0$ or Eq.\
(\ref{four}) for the $D_k$, while the $E_k$ are given by Eq.\ (\ref{seven}). The case of $p_0 \neq 0$ or
Eq.\ (\ref{three}) is left as an exercise for the reader. If we compare the code with the final code in Sec.\ 2, then
we see that the overall structure remains the same. That is, there is a main section with the same two function 
prototypes ${\bf termgen}$ and ${\bf idx}$. In fact, ${\bf idx}$ or the BRCP algorithm has not been altered at all, 
but ${\bf termgen}$ and the main routine have been changed to produce the symbolic forms for the coefficients
in the partition method for a power series expansion. Besides evaluating the execution time, ${\bf main}$ carries
out the calculation of the coefficients in one for loop, which is limited by the variable ${\it dim}$, representing 
the maximum value of $k$ or the coefficient of the highest order term which the user must input. Within this
for loop there are two calls to ${\bf idx}$, one of which applies to the calculation of the $D_k$ and the other
to the calculation of the $E_k$. Therefore, it is ${\bf termgen}$ that is doing the heavy work in the programme.
In fact, we shall see in the next section that by modifying ${\bf termgen}$ we can determine many of the properties 
of partitions, which often require separate programs.

Within ${\bf termgen}$ we see that the $D_k$ and $E_k$, which are represented by the variables ${\it DS[k,n]}$ and
${\it ES[k,n]}$ respectively, are evaluated depending upon the value of the variable ${\it inv_{-}case}$. If it 
equals zero, then the $D_k$ are evaluated, while if it equals unity, then the $E_k$ are evaluated. In evaluating 
the latter there is also an extra complication due to the phase factor of $(-1)^N$ in Eq.\ (\ref{seven}). Consequently, 
for this case ${\bf termgen}$ must determine the number of distinct elements in each partition. When ${\it dim } \!=\!4$,
${\bf partmeth}$ prints out the first four values of the $E_k$ and $D_k$ in symbolic form. E.g., the $k \!=\! 4$ 
values that it prints out are : \newline 
\indent DS[4,n$_{-}$]:= p[4,n] q[1] a + p[1,n] p[3,n] q[2] a$^{\wedge}$(2) 2! \newline
\indent + p[1,n]$^{\wedge}$(2) p[2,n] q[3] a$^{\wedge}$(3) 3!/2! + p[1,n]$^{\wedge}$(4) q[4] a$^{\wedge}(4)$ \newline 
\indent + p[2,n]$^{\wedge}$(2) q[2] a$^{\wedge}$(2) \newline

ES[4,n$_{-}$]:= -DS[0,0]$^{\wedge}$(-2) DS[4,n] + DS[0,0]$^{\wedge}$(-3) DS[1,n] DS[3,n] 2! \newline 
\indent - DS[0,0]$^{\wedge}$(-4) DS[1,n]$^{\wedge}$(2) DS[2,n] 3!/2! + DS[0,0]$^{\wedge}$(-5) DS[1,n]$^{\wedge}$(4) \newline 
\indent + DS[0,0]$^{\wedge}$(-3) DS[2,n]$^{\wedge}$(2) \newline 
From these results we see that each coefficient is composed of five distinct terms corresponding to the fact 
that the number partitions summing to 4, i.e.\ $p(4)$, is equal to 5. These results allow for the situation where
the ${\it p[k,n]}$ may be dependent upon another variable, viz.\ ${\it n}$, even though it may not be necessary.       

The first code presented in the appendix is suitable for values up to and around $k \!=\! 20$. In fact, all the values 
of $D_k$ and $E_k$ for $k \leq 20$ are computed within one CPU second. For $k \!\geq\! 20$, however, the expressions 
become unwieldy and thus, it is better to evaluate them separately so that each can be introduced directly into 
Mathematica. This amounts to removing the for loop in ${\bf main}$ and computing only for the value of $k$ or 
${\it dim}$. The code for computing the $D_k$, which is called ${\bf mathpm}$, appears immediately after ${\bf partmeth}$ 
in the appendix. Furthermore, in order that the coefficients can be introduced directly in Mathematica, only
three terms appear on each line of output, while a plus sign now appears as the last character on each line except,
of course, on the final line.

Let us now consider the evaluation of $D_{30}$ via ${\bf mathpm}$. Since $p(30) \!=\! 5604$, this is the number
of distinct terms when ${\bf mathpm}$ prints out ${\rm DS}[30,{\rm n}]$. Even though the output file for 
${\rm DS}[30,{\rm n]}$ is very large, it can still be imported into Mathematica. If we now set $p_k \!=\!(-1)^k/(2k+1)!$, 
$q_k \!=\!1 (-1)^k$ and $a \!=\!1$, which represent the inner and outer series for the cosecant numbers, then it takes 
0.15 CPU sec to evaluate $c[30]$ or $c_{30}$ in integer form on the same Sony VAIO laptop mentioned in Sec.\ 2. In this 
instance the numerator is given by a 60 digit integer, while the denominator is given by a 90 digit integer. In decimal 
notation the value of $c_{30}$ is approximately $2.965 \times 10^{-30}$. If we use Eq.\ (\ref{fortyseven}) to evaluate 
$c_{30}$, then we find that it takes almost zero CPU sec to evaluate the same result. On the other hand, if we set 
$p_k \!=\! (-1)^k/(2k)!$, i.e.\ the situation for the secant numbers $d_k$, then we find that it takes 0.14 CPU sec to 
evaluate $d_{30}$ in integer form on the Sony VAIO laptop. In this case the numerator and denominator are respectively 
67 and 78 digit numbers, while in decimal form $d_{30}$ is approximately equal to $2.176 \times 10^{-12}$. Unfortunately, 
if Eq.\ (\ref{fiftyone}) is implemented in Mathematica, then we only obtain approximate values in decimal form for the 
secant numbers. Hence, we need to implement a recurrence relation such as Eq.\ (\ref{fiftytwo1}) in order to obtain them 
in integer form. When this is done, it is found that Mathematica takes 6548 CPU sec to compute $d_{30}$. 

If we set $p_k \!=\! (-1)^{k}/(k+1)$, which represents the situation for the reciprocal logarithm numbers $A_k$, 
then we find that it takes only 0.1 CPU sec to determine $A_{30}$ in integer form. In this instance the numerator 
and denominator are 35 and 38 digit numbers yielding an approximate decimal value of $1.474 \times 10^{-3}$, the 
slowest converging of the numbers considered so far. The reciprocal logarithm numbers can be evaluated by either 
relating them to the Stirling numbers of the first kind \cite{kow10} via
\begin{eqnarray}
A_k = \frac{(-1)^k}{k!} \sum _{j=0}^{k-1} \frac{S_k^{(j)}}{j+1} \;\;.
\label{seventysix}\end{eqnarray}
or by the recurrence relation of
\begin{eqnarray}
A_k = \sum_{j=0}^{k-1} \frac{(-1)^{k-j+1}}{(k-j+1)}\; A_j \;\;.
\label{seventyseven}\end{eqnarray}   
If the first form is implemented in Mathematica, then it takes 0.1 CPU sec to compute $A_{30}$, while with the second
form it takes 5719 CPU sec. Therefore, we see that the evaluating the coefficients of power series expansions via the 
partition method for a power series expansion can be vastly superior to using recurrence relations and is almost on a 
par with the cases where intrinsic forms have already been implemented within a mathematical software package. 
Furthermore, by altering the relations for the coefficients of the inner and outer series in addition to $a$, we obtain 
results for other mathematical quantities with different power series expansions.

As discussed previously, the coefficients of the inner and outer series do not need to be numbers as has been the 
case so far. If we set $q_k \!=\! (\rho)-k/k!$, $p_k \!=\! (-1)^{k+1}/(2k+1)!$ and $a \!=\!1$, 
then {\it DS[30,n]} yields the generalised cosecant number $c_{\rho,30}$ as given by Eqs.\ (\ref{sixty}) and (\ref{sixtyone}).  
Then we find that it only takes 0.36 CPU sec to compute the resulting polynomial, which is thirtieth order in $\rho$. By 
invoking the Simplify routine in Mathematica the polynomial can be arranged in increasing order within another 0.36 CPU 
sec. Before this calculation was performed, the results produced by {\it DS[5,n]} and {\it DS[8,n]} had been found to 
agree with the generalised cosecant numbers, $c_{\rho,5}$ and $c_{\rho,8}$, obtained in Ref.\ \cite{kow11a}.
 
As stated earlier, the calculation of the coefficients via the forms generated by either of the first two 
programs in the appendix can be continued beyond the thirtieth order, but eventually problems arise due
to the combinatorial explosion occurring in the number of partitions. It has already been stated that 
there are $190\,569\,272$ partitions summing to 100, which means that this number of terms will be present
in {\it DS[100,n]}. If ${\bf mathpm}$ is run for this case, then it takes around 600 CPU sec to compute {\it DS[100,n]}. 
Whilst this is not an overly long time of computation compared with the earlier results obtained via recurrence
relations, it produces a file whose size is over 16 GB. Files of this size are going to present a problem when imported
into mathematical software packages. For example, it appears that Mathematica\ 8.0.1 is only able to import files
with 2Gb of data. One method of circumventing this problem would be to divide the file into smaller files so 
that could be handled by the different processors on a supercomputer. Then the results generated by each processor
could combined to yield the final answer.

Another method of overcoming this problem is to introduce the values for $p_k$, $q_k$ and $a$ first and then
evaluate a specific number or limit of terms via ${\bf mathpm}$. Once the limit point is reached, the values
where this occurs would need to be stored. In terms of Fig.\ \ref{figone} this amounts to storing the values of
both arguments in ${\bf idx}$. Then the partial value of the coefficient could be evaluated and stored, while
all the terms outputted in running ${\bf mathpm}$ can be either deleted or overwritten in a re-run of ${\bf mathpm}$.
In the re-run of ${\bf mathpm}$ the code would not print out any terms until ${\bf idx}$ reaches the values 
of the arguments of ${\bf idx}$ stored from the first run. Then the code would either continue to print out 
the next limit of terms stopping at two new values of ${\bf idx}$ or would terminate on reaching the central 
partition. Then the terms stored in the second run could be evaluated and combined with the result obtained from the
first run. If the central partition has not been reached in the second run, then the process can be continued until
the central partition is eventually reached. Of course, the disadvantage in this approach is that we have lost
the ability to evaluate a new coefficient by altering the $p_k$, $q_k$ and $a$ as we were able to do with
{\it DS[30,n]} above. This second method of solving the problem of very large data files produced by running 
${\bf mathpm}$is contingent on whether we can stop and re-start the program at specific points in the tree diagrams 
for the partitions. Hence, we need to be able to adapt the BRCP algorithm so that specific partitions can be 
determined, which represents the topic of the following section.    
   
\section{Specific Types of Partitions}
In this section we aim to investigate how the BRCP algorithm presented in Sec.\ 2 can be modified to determine
specific types of partitions. By specific partitions we mean such general problems in partitions as the determination 
of: (1) partitions with a fixed number of elements, (2) doubly-restricted partitions, (3) discrete partitions, 
(4) conjugate partitions and (5) partitions with specific elements in them. Solving such problems invariably 
means creating a different algorithm or program for each problem as can be seen in Refs.\ \cite{yam07}, \cite{knu005} 
and \cite{fen80}. However, we shall see here that such problems and others can be solved with relatively 
minor modifications to the BRCP algorithm presented in Sec.\ 2, once again highlighting its versatility. As we
shall be modifying ${\bf partgen}$ in Sec.\ 2 when solving these problems, we shall generate the partitions in 
the compact multiplicity representation, although it should be mentioned that the various algorithms presented in 
this section would only require minor modification to ${\bf termgen}$ to generate partitions in the standard 
representation.  

Of all the problems mentioned in the previous paragraph perhaps the simplest one to consider is the determination 
of those partitions with a specific element or elements in them. E.g., suppose we wish to determine all those 
partitions summing to 15 with the element of \{5\} in them? From the material presented in Sec.\ 2, it is obvious 
that the total number of partitions for this problem must equal the number of partitions summing to 10 or $p(10)$. 
Moreover, the partitions generated by the new code should yield the same partitions obtained by running the various 
codes presented in Sec.\ 2 except that each partition generated by the new code will have an additional element of 
\{5\} in them. We shall see that this is indeed the case, although the order in which the partitions are generated 
is different from those discussed in Sec.\ 2.     

In order to modify {\bf partgen} in Sec.\ 2 so that it generates partitions with a specific element
in them, we first need to alter the {\bf main} prototype of the program. This is required to enable
users to input the specific values of the elements that they wish to appear in the partitions printed out
by the new program. Consequently, {\bf main} becomes

\begin{lstlisting}{}
int main( int argc, char *argv[] )
{
int i;
if(argc!=3) printf("usage: specpart <partition sum> 
	<sp_val> \n");
else{
        tot=atoi(argv[1]);
        sp_val=atoi(argv[2]);
        part=(int *) malloc(tot*sizeof(int));
        if(part == NULL) printf("unable to allocate 
		array\n");
        else{
                for(i=0;i<tot;i++) part[i]=0;
                idx(tot,1);
                free(part);
            }
     }
printf("\n");
return(0);
}
\end{lstlisting}
Here we see that the new program called {\bf specpart} has a global variable called ${\it sp_{-}val}$, which
represents the specific element that is to appear at least once in each partition printed out by the program.

As discussed in Sec.\ 2, {\bf idx(tot,1)} scans over all the partitions summing to {\it tot}, while {\bf termgen},
which is called in {\bf idx}, is responsible for generating or printing out partitions. Thus, in order to determine those
partitions with a specific element or elements in them, we need to modify {\bf termgen}. That is, we still need to
scan over all partitions by calling {\bf idx(tot,1)}. In fact, aside from making minor modifications to {\bf main}, we 
shall find that to solve all the problems mentioned in the introduction to this section, we need only modify {\bf termgen}. 

In Sec.\ 2 the ${\bf termgen}$ function prototype in {\bf partgen} was responsible for printing out all the partitions 
summing to {\it tot} in the multiplicity representation. This was achieved by processing the array {\it part}, which 
stored the frequencies of the elements in each partition. That is, the variable {\it freq} was used to represent the 
frequency of the element $i+1$ in the partition with $i$ ranging from 0 to ${\it tot -1}$. Now that we wish to determine 
those partitions with a specific element or elements in them, we need to restrict the partitions that are printed out 
by ${\bf termgen}$. This is accomplished simply by introducing a local variable called ${\it freq_{-}spval}$, which 
evaluates the frequency of ${\it sp_{-}val}$ in each partition. If this value is non-zero, then we know that there 
is at least one occurrence of ${\it sp_{-}val}$ in the partition and the partition is then printed out in the same 
manner as in Sec.\ 2. If ${\it freq_{-}spval}$ is zero, then the partition is ignored. Therefore, the {\bf termgen} 
function prototype for {\bf specpart} becomes
\begin{lstlisting}{}
void termgen()
{
int freq,i,freq_spval;

freq_spval=part[sp_val-1];
if(freq_spval){
        printf("%ld: ",term++);
        for (i=0;i<tot;i++){
               freq=part[i];
               if(freq) printf("%i(%i) ",freq,i+1);
                           }
        printf("\n");
              }
}
\end{lstlisting}   

In Sec.\ 2 the partitions summing to 5, which amounted to 7, were printed out by running ${\bf partgen}$. The output 
produced by running ${\bf specpart}$ with ${\it tot}$ and ${\it sp_{-}val}$ set equal to 11 and 6, respectively, is:\newline
1: 5(1) 1(6) \newline
2: 3(1) 1(2) 1(6) \newline
3: 2(1) 1(3) 1(6) \newline
4: 1(1) 2(2) 1(6) \newline 
5: 1(1) 1(4) 1(6) \newline
6: 1(2) 1(3) 1(6) \newline
7: 1(5) 1(6) \newline
Hence, we observe that the number of partitions is once again 7 or $p(5)$. If we remove one six from each partition, 
then we obtain the same partitions as those generated in Sec.\ 2 by ${\bf partgen}$ except that the order in which
they appear is now different. If we define the specific element partition operator $L_{SEP,k,j}[ \cdot]$ as
\begin{eqnarray}
L_{SEP,k,j}\Bigl[ \cdot \Bigr]= \sum_{\scriptstyle n_1,\dots,n_{j-1}=0,n_{j}=1,n_{j+1},\dots,n_{k}=0 
\atop{\sum_{i=1,i\neq j}^{k} i n_i=k-j}}^{k,\dots,[k/(j-1)],[k/j],[k/(j+1)] \dots,1} \;\;,
\label{five-one}\end{eqnarray}
then it follows that 
\begin{eqnarray}
L_{SEP,k,j}\Bigl[ 1 \Bigr] = L_{P,k-j}\Bigl[ 1\Bigr]= p(k-j) \;\;.
\label{five-two}\end{eqnarray}

As a result of the preceding analysis, it becomes a simple matter to consider partitions with more than 
one specific element in them. All we need to do is introduce more values in ${\bf main}$ and then create
local variables like ${\it freq_{-}spval}$ to represent the number of occurrences or frequencies for each 
of these values. As before, each frequency is set equal to ${\it part[sp_{-}val-1]}$ in ${\bf termgen}$. Next, 
we modify the if statement to include all the frequencies. For example, if we wish to determine all the 
partitions with two specific elements, viz.\ ${\it sp_{-}val}$ and ${\it sp_{-}val2}$, then the if statement 
becomes
\begin{lstlisting}{}
if (freq_spval && freq_spval2){  etc.         
\end{lstlisting}
On the other hand, we may want those partitions with at least one occurrence of either ${\it sp_{-}val}$ or
${\it sp_{-}val2}$. In this instance all we need to do is replace the logical AND operator by the logical 
OR operator or \&\& in the if statement.

By making a few minor modifications, mostly to {\bf termgen}, we have been able to solve several 
problems involving specific partitions. Now we consider evaluating the partitions with a fixed 
number of elements in them, which has already been considered when we modified the partition operator 
into the form given by Eq.\ (\ref{fortyfourl}). Previously, it was stated that to solve this problem,
often a completely different approach is employed. For example, Knuth presents a different algorithm 
on p.\ 38 of Ref.\ \cite{knu005} compared with the algorithm used to generate partitions in reverse 
lexicographic form. In the case of the BRCP algorithm, however, the number of elements in a partition is 
determined by the number of branches along its path prior to termination in the tree diagram. E.g., the 
partitions with only two elements in them in Fig.\ \ref{figone} are obtained by searching for the terminating 
tuples appearing in the third column, viz.\ (0,5), (0,4) and (0,3). In other words, these are the terminating
tuples two branches away from the seed number situated in the first column. When we include 
the elements from the second column, they yield the partitions \{1,5\}, \{2,4\} and \{3,3\}. Since
partitions with a fixed number of elements can be determined from the tree diagram, it means that only 
minor modifications to the program in Sec.\ 2 are again required in order to generate these partitions.

As in the previous example, we need to modify {\bf main} so that the user can input the number of parts, 
which will be represented by the global variable {\it numparts}. Then we can proceed to the modification 
of {\bf termgen}, which is presented below. There we see that a new local variable called ${\it sumparts}$ 
has been introduced. This variable determines the number of elements in each partition. When this number equals 
${\it numparts}$, the partition is printed out. Otherwise, the partition is discarded. 

\begin{lstlisting}{}
void termgen()
{
int freq,i,sumparts=0;
/* sumparts is the number of parts or elements in a 
   partition */


for(i=0;i<tot;i++){
      sumparts= sumparts+part[i];
                  }
if(sumparts == numparts){
    printf("%ld: ",term++);
    for(i=0;i<tot;i++){
                freq=part[i];
                if(freq) printf("%i(%i) ",freq,i+1);
                      }
     printf("\n");
                        }
}
\end{lstlisting}

When the code is run for the case where the partitions sum to 10 and possess 5 elements, i.e.\ for ${\it tot} \!=\! 10$
and ${\it numparts} \!=\! 5$, the following output is produced: \newline
1: 4(1) 1(6) \newline
2: 3(1) 1(2) 1(5) \newline
3: 3(1) 1(3) 1(4) \newline
4: 2(1) 2(2) 1(4) \newline
5: 2(1) 1(2) 2(3) \newline
6: 1(1) 3(2) 1(3) \newline
7: 5(2) \newline
Therefore, we see that the number of partitions in this instance is 7, which can also be expressed as either 
$\begin{vmatrix} 10 \\ 5 \end{vmatrix}$ according to  Sec.\ 2 or as $L^{5}_{P,10}[1]$ according to Eq.\ 
(\ref{fortyfourm}).

It should also be mentioned that if the condition,   
\begin{lstlisting}{} 
if(sumparts==numparts){ etc., 
\end{lstlisting}
in the above program is replaced by 
\begin{lstlisting}{} 
if(sumparts<=numparts){ etc., 
\end{lstlisting} 
then the resulting code generates all those partitions summing to ${\it tot}$ with at most ${\it numparts}$ 
elements. E.g., the number of partitions summing to 10 with at most 5 elements, which according to Sec.\ 2  
is equal to 30. As stated in Sec.\ 2, this is also equivalent to $\begin{vmatrix} 15 \\ 5 \end{vmatrix}$.
However, the partitions generated by the code given above for ${\it tot}$ and ${\it numparts}$ equal to 15 
and 5 respectively, are different from those generated by the version of the program with the modified 
${\bf if}$ statement.

Doubly-restricted partitions are those partitions where all the elements are greater than a particular
value and less than another value. Since this is a combination of two separate conditions, first we need
to be able to modify ${\bf partgen}$ so that it generates the partitions where the elements are either greater 
than or lower than a specified value. Therefore, let us consider the situation where all the elements in the
partitions are less than or equal to a value, which will be represented by the global variable 
${\it largest_{-}elt}$. As in the previous examples, this value must be introduced into ${\bf main}$. Then
we proceed to the modification of ${\bf termgen}$.

To modify ${\bf termgen}$ so that only partitions with elements less than or equal to ${\it largest_{-}elt}$ 
are generated, we need to introduce an extra for loop. This extra loop is required so that if a 
partition is encountered where an element is greater than ${\it largest_{-}elt}$, it is discarded via a 
${\bf goto}$ statement as demonstrated by the modified version of ${\bf termgen}$ given below. Although ${\bf goto}$ 
statements are generally frowned upon by programmers, it is being used here to abandon processing in a 
nested structure of two loops. In fact, the code behaves much like ${\bf partgen}$ in Sec.\ 2
when all the elements in the partitions are less than or equal to ${\it largest_{-}elt}$. However, when an 
element is greater than ${\it largest_{-}elt}$, the ${\bf goto}$ statement discards the partition by 
diverting to the ${\bf end}$ statement label. 

\begin{lstlisting}{}
void termgen()
{
int f,i;
for(i=largest_elt;i<dim;i++){
      f=part[i];
      if(f)
          goto end;
                              }
printf("%ld: ",term++);
for(i=0;i<dim;i++){
      f=part[i];
      if(f) printf("%i(%i) ",f,i+1);
                   }
printf("\n");
end:;
}
\end{lstlisting}
When the above code is run for partitions summing to 14 in which the elements are greater than 3, the following 
output is produced: \newline
1: 1(14) \newline
2: 1(4) 1(10) \newline
3: 2(4) 1(6)  \newline
4: 1(4) 2(5)  \newline
5: 1(5) 1(9)  \newline
6: 1(6) 1(8) \newline
7: 2(7)\newline 
Hence, we see that there only seven partitions summing to 14, in which all the elements are greater than 3. 

As a result of the above code, it is now a simple matter to consider the case where all the elements
are greater than or equal to another value specified by the user. In this instance we simply replace 
${\it largest_{-}elt}$ by another global variable called ${\it smallest_{-} elt}$ and alter the 
condition in the first for loop of the previous version of ${\bf termgen}$. That is, the
first for loop in the preceding version of ${\bf termgen}$ simply becomes 
\begin{lstlisting}{}
for (i=0;i<smallest_elt-1;i++){   etc.
\end{lstlisting} 
Note also that if the upper value in the for loop had been ${\it smallest_{-}elt}$ instead of 
${\it smallest_{-}elt}-1$, then all the elements generated by the code would only be greater than the 
value specified by the user.

For doubly-restricted partitions, where the elements are greater or equal to value and less than or 
equal to another (larger) value, all we need to do is incorporate two for loops that divert
to the ${\bf end}$ statement label. For example, ${\bf termgen}$ for this situation would become
\begin{lstlisting}{}
void termgen()
{
int f,i;
for (i=0;i<smallestpart-1;i++){
      f=part[i];
      if (f)
           goto end;
                              }
for(i=largestpart;i<dim;i++){
      f=part[i];
      if(f)
          goto end;
                            }
      printf("%ld: ",term++);
      for(i=0;i<dim;i++ ){
          f=part[i];
          if(f) printf("%i(%i) ", f,i+1);
                          }
printf("\n");
end:;
}
\end{lstlisting}

When the above code is run for partitions summing to 13 with the elements greater than or equal to 3 and less than 
or equal to 9, the following output is produced: \newline
1: 2(3) 1(7) \newline
2: 3(3) 1(4) \newline
3: 1(3) 1(4) 1(6) \newline
4: 1(3) 2(5) \newline
5: 1(4) 1(9) \newline
6: 2(4) 1(5) \newline
7: 1(5) 1(8) \newline
8: 1(6) 1(7) \newline
Thus, we see that there are 8 partitions with all elements lying in the interval [3,9].

In Ch.\ 3 of Ref.\ \cite{and03} Andrews defines restricted partitions differently from the earlier definition 
given by Eq.\ (\ref{sixtyseven}). There, the partitions represent those in which the elements are less than a 
value, say ${\it el_{-}max}$, while the number of elements is less than or equal to another value, which we take 
again to be ${\it numparts}$ as in the preceding examples. Although there is now a condition pertaining
to the number of parts, evaluating these partitions is again similar to the doubly-restricted case studied above. 
First, we must introduce ${\it el_{-}max}$ into ${\bf main}$ in addition to ${\it numparts}$. Then we need to insert
an extra for loop into ${\bf termgen}$ so that it can make use of the different condition. The new loop appears first 
since if it is true, we immediately by-pass any action to process the current partition. Therefore, this modified 
version of 
${\bf termgen}$ becomes
\begin{lstlisting}{} 
void termgen()
{
int freq,i,sumparts=0;
for(i=0;i<tot;i++){
       if(i> el_max-1 && part[i]>0) goto end;
                  }
/*(1) sumparts is the number of elements in a partition 
  (2) all elements are now less than or equal to el_max 
*/

for(i=0;i<tot;i++){
      sumparts= sumparts+part[i];
                  }
if(sumparts <= numparts){
    printf("%ld: ",term++);
    for(i=0;i<tot;i++){
                freq=part[i];
                if(freq) printf("%i(%i) ",freq,i+1);
                      }
     printf("\n");
                        }
end: ;
}
\end{lstlisting}
When the code is run with ${\it tot}$, ${\it numparts}$ and ${\it el_{-}max}$ set equal to 10, 3 and 5 respectively, 
the following output is produced: \newline
1: 1(1) 1(4) 1(5) \newline
2: 1(2) 1(3) 1(5) \newline
3: 1(2) 2(4) \newline
4: 2(3) 1(4) \newline
5: 2(5) \newline
Hence, we see that there are 5 partitions summing to 10 with at most 3 elements and all elements less than or equal to 5.

The number of partitions summing to $k$ with at most $M$ parts and each element less than or equal to $N$ is represented 
as $p_G(N,M,k)$ in Ref.\ \cite{and03}. The subscript G has been introduced here so that the reader is not confused with
similar notation in the next section. From the above example we have $p_G(5,3,10) \!=\! 5$. If $k \!>\! MN$, then 
$p_G(N,M,k)$ vanishes, while $p_G(N,M,NM) \!=\! 1$. These numbers also appear as the coefficients in the generating function 
for Gaussian polynomials, which are given by
\begin{eqnarray}
G(N,M;q) = \prod_{i=1}^N \frac{(1-q^{M+i}\,)}{(1-q^i\,)} = 1+ \sum_{k=1}^{NM} p_G(N,M,k) \,q^k\;\;.
\label{five-three}\end{eqnarray} 
Hence, Gaussian polynomials are polynomials in $q$ of degree $NM$. Moreover, to avoid confusion with the restricted
partitions studied earlier, we shall refer to the above partitions as Gaussian partitions. As a consequence, we define 
the Gaussian partition operator $L_{GP,k,N,M}[ \cdot]$ as
\begin{eqnarray}
L_{GP,k,N,M}\Bigl[ \cdot \Bigr]= \sum_{\scriptstyle n_1,n_2\dots,n_{N}=0 
\atop{\sum_{i=1}^{N} i n_i=k,  \quad \sum_{i=1}^{N} n_i \leq M}}^{M,Min\{[k/2],M\},\dots,Min\{[k/i],M\}} \;\;.
\label{five-four}\end{eqnarray}
As a result, we have $L_{GP,k,N,M}[1] \!=\! p_G(N,M,k)$. 

We now turn our attention to a more complicated example--- the problem of determining discrete or distinct partitions.
By discrete partitions, we mean those partitions in which the elements appear at most once, if at all. Like all
the preceding examples, they too represent a subset or class of the set of integer partitions. As we shall see in 
the next section, such partitions figure prominently in the theory of partitions. Because of this, we define the 
discrete partition operator, $L_{DP,k}[ \cdot]$, by  
\begin{eqnarray}
L_{DP,k}\Bigl[ \cdot \Bigr]= \sum_{\scriptstyle n_1,n_2\dots,n_{k}=0 
\atop{\sum_{i=1}^{k} i n_i=k}}^{1,1,\dots,1} \;\;.
\label{five-five}\end{eqnarray}
Therefore, the only difference between this operator and the partition operator introduced in Sec.\ 3 is that 
the upper limits of the summations are restricted to unity, whereas for the latter the upper limits were set to 
$[k/i]$ for each element $i$.

\begin{figure}
\begin{center}
\includegraphics{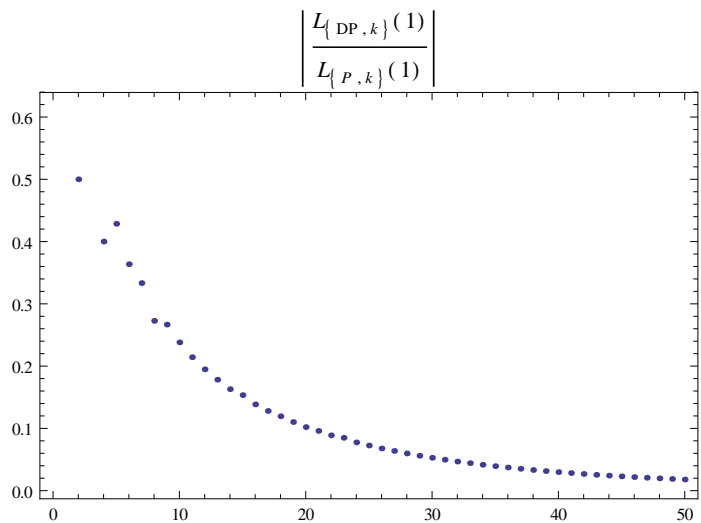}
\end{center}
\caption{The number of discrete partitions to the total number of partitions versus $k$}
\label{figthree}
\end{figure}

A computer program that generates discrete partitions is presented in its entirety in the appendix. As stated
previously, the number of partitions summing to 100 or $p(100)$ equals $190\,569\,272$, which is the reason why
it is cumbersome to evaluate the 100-th coefficient via the partition method for a power series expansion. On the 
other hand, if we run the program in the appendix called ${\bf dispart}$ to determine the discrete partitions 
summing to 100, then we find that it only takes 5 CPU seconds to generate the $444\,793$ partitions, which 
represent 0.2 per cent of $p(100)$. Moreover, Fig.\ \ref{figthree} displays the ratio of the number of discrete 
partitions or $L_{DP,k}[1]$ to the number of standard partitions versus $k$ for $k \! \leq \! 50$. From the figure 
we see that this ratio decreases monotonically for $k \! \geq \! 10$, reaching a value of $ 0.0179\cdots$ 
when $k \!=\! 50$.

The final example we shall consider in this section is the evaluation of the conjugate partition as a partition
is generated. According to Ref.\ \cite{knu005}, the conjugate partition is obtained by transposing the rows and
columns of the corresponding Ferrers diagram for the original partition. In a Ferrers diagram a partition is
represented by an array of dots in which the first element, say $a_1$, is allocated a row of $a_1$ dots, while 
the next element ($a_2$) is represented by another row of $a_2$ dots immediately below the first row of dots. The
process of allocating rows of dots for each element in a partition is continued up to the final element in the
partition. The conjugate partition or $\alpha^{T}$ of the partition $\{a_1,a_2,\dots,a_k\}$ is obtained by
transposing the rows and columns of the Ferrers diagram for the partition. For example, the Ferrers diagram
for the partition $\{1,1,2,3,3,4\}$ has two rows with one dot, a row with two dots followed by two rows
with 3 dots and finally a row with 4 dots. The transpose is obtained by counting the dots in the vertical columns
of the Ferrers diagram. Hence, $\alpha^T$ for our example is found to be $\{6,4,3,1\}$. Note also that the conjugate
partition does not necessarily possess the same number of elements as the original partition.   
    
A code that determines the conjugate partition for each partition generated by ${\bf partgen}$ is presented in 
the appendix. As with all the other examples in this section, it is ${\bf termgen}$ that has been modified. The 
major complication in this code compared with the others is that one is required to create and allocate space for a 
two-dimensional array ${\it ferrers}$ of size ${\it tot} \times {\it tot}$ in addition to creating a one-dimensional 
array of pointers called ${\it rptr}$ to it. In both cases the arrays are of integer type. In reality they should be 
of type char since a true Ferrers diagram consists of dots. That is, instead of allocating dots the program 
allocates ones in creating a Ferrers diagram. Nevertheless, by adding the ones vertically, one obtains the conjugate 
partition. For the purists it should be a simple matter to alter the ${\bf termgen}$ to count dots rather than ones.  

When the program called ${\bf transp}$ for the partitions summing to 4, the following output is produced: \newline
Partition 1 is: 1(4)  and its conjugate is: 4(1)  \newline
Partition 2 is: 1(1) 1(3)  and its conjugate is: 1(2) 2(1) \newline 
Partition 3 is: 2(1) 1(2)  and its conjugate is: 1(3) 1(1) \newline
Partition 4 is: 4(1)  and its conjugate is: 1(4) \newline
Partition 5 is: 2(2)  and its conjugate is: 2(2) \newline
As can be seen from the output, the conjugate partitions are printed out in a different order to those produced
by the BRCP algorithm. For example, the conjugate partition to Partition 2 or \{1,3\} is printed out as 1(2) 2(1) or 
\{2,1,1\}, whereas the form produced by the BRCP algorithm, viz.\ Partition 3, is 2(1) 1(2) or \{1,1,2\}. Throughout this 
work we have not been concerned with the order of the elements in the partitions. When the order is important, partitions
become compositions \cite{and03}. For example, there are three compositions of the partition \{1,1,2\} are: (112), (121) and
(211). We shall not be concerned with compositions here.

There are some interesting properties concerning conjugate partitions. For example, there is at least one partition 
whose conjugate is itself, although it can turn out to be a different composition. Such self-conjugate partitions 
arise when the partitions sum to a square of an integer $k$ since the partition $k(k)$ in the multiplicity representation
is a self-conjugate. On the other hand, if the partitions sum to $k(k+1)/2$, then the partition \{1,2,\dots,k\} is a
self-conjugate, although again the composition is different. When partitions sum to $2k$ for $k \!>\!2$, the 
partition given by $(k-2)(1) 1(2) 1(k)$ is also a self-conjugate, while if they sum to $2k \!+\!1$, where $k\geq 7$, 
then the partition given by $(k-6)(1) 1(2) 1(3) 1(k-2)$ is another self-conjugate. In fact, according to Ref.\ \cite{wik11},
the number of self-conjugate partitions is the same as the number of partitions with distinct odd elements. 

There are other problems that can be solved by modifying ${\bf termgen}$ given in Sec.\ 2. One such problem is 
determining those partitions with only odd elements in them. As a consequence, we can also define an odd-element  
partition (OEP) operator as 
\begin{eqnarray}
L_{OEP,2k+1}[\cdot] = \sum_{\scriptstyle n_1,n_3, \dots,n_{2k+1}=0 \atop{\sum_{i=1}^{k}(2i+1) n_{2i+1}=2k+1}}^{2k+1,[(2k+1)/3], \dots,1}  
\;\;,
\label{five-six}\end{eqnarray} 
which is only valid for partitions summing $(2k\!+\!1)$, while for those summing to $2k$, the operator becomes
\begin{eqnarray}
L_{OEP,2k}[\cdot] = \sum_{\scriptstyle n_1,n_3, \dots,n_{2k}=0 \atop{\sum_{i=1}^{k}(2i-1) n_{2i-1}=2k}}^{2k,[2k/3], \dots,1}  
\;\;.
\label{five-sixa}\end{eqnarray}
To obtain partitions with only odd elements in them, all we need to do is insert the following for loop at the beginning
of ${\bf termgen}$ just after the type declarations:
\begin{lstlisting}{}
for (i=0;i<tot;i++){
             freq=part[i];
             if (i % 2 && freq > 0) goto end;
                   }
printf("%ld: ",term++);      .
\end{lstlisting}

It has already been mentioned that the number of partitions with discrete odd elements equals the number of self-conjugate 
partitions. To obtain the former, the upper limits in the definitions for $L_{OEP,2k+1}[\cdot]$ and $L_{OEP,2k}[\cdot]$ must be
set equal to unity. In addition, if we wish to generate partitions with discrete odd elements in them, all we need to do is 
introduce the above for loop at the beginning of ${\bf termgen}$ in the program ${\bf dispart}$ that is presented in 
the appendix. With regard to self-conjugate partitions, ${\bf transp}$ in the appendix would need to be modified so that the 
original partition is stored in a temporary array before it undergoes conjugation. Then a test would need to be introduced 
to see if both partitions are identical to one another. If they are, then the partition is printed out. Otherwise, it is 
discarded. This problem is left as an exercise for the reader. 
    
In a similar manner we can define an even-element partition (EEP) operator that only applies to those partitions summing to $2k$
with even elements in them. This is defined as 
\begin{eqnarray}
L_{EEP,2k}[\cdot] = \sum_{\scriptstyle n_2,n_4, \dots,n_{2k}=0 \atop{\sum_{i=1}^{k}i n_{2i}=k}}^{k,[k/2], \dots,1}  
\;\;.
\label{five-seven}\end{eqnarray} 
The number of partitions generated by this operator is equal to the number of partitions summing to $k$. Hence, we find 
that
\begin{eqnarray}
L_{EE,2k} [1] = L_{P,k}[1] \;\;.  
\label{five-sevena}\end{eqnarray} 
In Ref.\ \cite{kow11} the cosecant numbers or $c_k$ are first derived in terms of the partitions summing to $2k$ with even
elements in them before the form in terms of the partition operator given by Eq.\ (\ref{fortyseven}) is derived. Consequently,
we arrive at the following identity:
\begin{align}
& L_{P,k}  \left[ (-1)^{N_k} N_k! \prod_{i=1}^{k}  \left(\frac{1}{(2i+1)!}\right)^{n_i} \;\frac{1}{n_i!}  \right]  
\nonumber \\ 
&=  L_{EEP,2k}\left[ (-1)^{N^{*}_{k}} N^{*}_{k}! \prod_{i=1}^{k} \left(\frac{1}{(2i+1)!}\right)^{n_{2i}} \;\frac{1}{n_{2i}!}  \right] \;\;,
\label{five-eight}\end{align} 
where $N^{*}_{k} = \sum_{i=1}^k n_{2i}$.

In this section various programs have been presented for obtaining specific types/classes of partitions, which can be regarded
as subsets of the total number of partitions summing to a particular value. Consequently, the various operators given in this 
section represent restricted forms of the partition operator $L_{P,k}[\cdot]$. In the previous section we showed how the 
partition method for a power series expansion could be developed into a computer program scanning the entire set 
of partitions. From the material presented in this section, it should, therefore, be possible to determine the 
contributions that the specific partitions contribute to the partition method for a power series expansion. In particular, 
in the next section we shall see that the partition function $p(k)$ can be obtained via the partition method for a power 
series expansion involving a restricted set of the partitions summing to $k$, although this restricted set is more difficult 
to evaluate than many of the examples considered in this section. We shall not only present the program that can generate 
this restricted set of partitions, but also present the program that outputs the partition function in symbolic form 
so that it can be handled by Mathematica. As for the more important programs discussed in this section, the new program will
appear in the appendix. Finally, it should be noted that by possessing the capacity to adapt the partition method for a 
power series expansion to handle specific partitions, we may be able to determine which partitions make the largest 
contribution to the coefficients in Theorem\ 1. As a result, accurate approximations to the coefficients can be 
considered, which may avoid the combinatorial explosion that occurs for large orders as described at the beginning of the 
previous section. 

\section{Generating Function for ${\bf P(k)}$}
An important topic in the theory of partitions are the generating functions whose power series expansions possess coefficients 
that are dependent on the properties of partitions. One of the greatest achievements in this context is the derivation 
of the asymptotic formula for the number of partitions $p(k)$. The first step that led to this formula was the derivation of
a remarkable formula for the product $P(z)$ in Equivalence\ (\ref{fortyfoura}) by Dedekind. As explained in Ref.\ \cite{knu005}, 
this can be derived by the application of standard analytic techniques, namely Poisson's summation formula, to the logarithm 
of $P(z)$. Then by studying the behaviour of Dedekind's formula for $\ln P(\exp(-t)$ with $\Re\, t \!>\! 0$, Hardy and Ramanujan 
\cite{har18} were able with amazing insight \cite{and03} to deduce the asymptotic behaviour of the partition function $p(k)$ for 
large $k$. Eventually, the asymptotic behaviour of $p(k)$ was completely evaluated by Rademacher \cite{rad37}, which culminated 
in the now famous Hardy-Rademacher-Ramanujan formula mentioned in the introduction to Sec.\ 4.

Although we shall not reach such lofty heights, we shall, nevertheless, turn our attention in the remainder of this work to 
how the partition method for a power series expansion can be applied in the analysis of the various generating functions that 
occur in the theory of partitions and their extensions or generalisations. We begin in this section by applying the partition 
method for a power series expansion to $P(z)$, but before we can embark upon this task, we need to determine when the power series
expansion or generating function is convergent or in another words, for what values of $z$ Equivalence\ (\ref{fortyfoura}) becomes 
an equation. According to Knuth \cite{knu005}, it was Euler who noticed that the coefficient of $z^n$ in the infinite product of 
$(1+z+z^2+z^3+ \cdots+z^j+ \cdots) (1+z^2+z^4+\cdots+z^{2k} + \cdots) (1+z^3+z^6+ \cdots+ z^{3k} + \cdots) \cdots$ is the number 
of non-negative solutions to $k+2k+ 3k+ \cdots = n$ or the partition function $p(n)$ and that $1+z^m+ z^{2m} + \cdots$ equals 
$1/(1-z^m)$. As a result, he arrived at Equivalence\ (\ref{fortyfoura}) except that the equivalence symbol was replaced by an 
equals sign, which is not entirely correct as can be seen by the following theorem.

{\bf Theorem\ 2}. The equivalence statement relating $P(z)$ or $\prod_{m=1}^{\infty} 1/(1-z^m)$ to the generating function given
by the power series expansion with coefficients equal to $p(k)$, viz.\ Equivalence\ (\ref{fortyfoura}), is absolutely convergent 
for $|z| \!<\! 1$, in which case the equivalence symbol can be replaced by an equals sign. On the other hand, it is divergent for 
all other values of $z$. Then $P(z)$ represents the regularised value of the series on the rhs. For $|z| \!=\! 1$, the generating
function is singular.

{\bf Proof}. The reason why an equivalence symbol has had to be introduced is due to Euler's second observation concerned with
the geometric series. Replacing the series by its limit value of $1/(1-z^m)$ is strictly not valid for all values of $z$ as 
described in Refs.\ \cite{kow09}, \cite{kow09a}, \cite{kow001} and \cite{kow11b}. There it is shown that the standard geometric 
series, i.e.\ $\sum_{k=0}^{\infty} z^k$, is divergent for $\Re\, z \!>\! 1$, absolutely convergent for $|z| \!<\! 1$ and 
conditionally convergent for $|z| \!>\! 1$ and $\Re\, z \!<\! 1$. Regardless of the type of convergence, the limit value of the 
series is found to equal $1/(1-z)$. For $\Re\, z \!>\! 1$, however, summing the series yields an infinity. In this case, if the 
infinity is removed after the series is summed, which is the essence of regularisation, then the remaining finite part is found 
to equal $1/(1-z)$ again. Hence, for $z \!\geq\! 1$, the regularised value of the geometric series is equal to $1/(1-z)$. 
Moreover, along the line $\Re\, z \!=\! 1$, the limit of the series is undefined or indeterminate, while at the point $z \!=\! 1$, 
where the line is tangent to the unit disk of absolute convergence, it is singular. This type of behaviour is expected since 
$\Re\, z \!=\! 1$, represents the border between the convergence to the left and divergence to the right. Since the limit value is 
the same on both sides after regularisation, we set the regularised value to $1/(1-z)$ along $\Re \, z \!=\! 1$. 

With regard to Equivalence\ (\ref{fortyfoura}) we are dealing with an infinite product of geometric series involving 
different powers of $z$ in the limit value. Nevertheless, we can use the above knowledge of the geometric series to determine
where the series on the rhs of Equivalence\ (\ref{fortyfoura}) is convergent and where it is divergent. For the product on the 
lhs to equal the series on the rhs of Equivalence\ (\ref{fortyfoura}), we must have for all positive integer values of $l$, 
$\Re\, z^l < 1$. For $l \!=\! 1$ we end up with the standard geometric series, but for $l \!=\! 2$, the series will now only 
be convergent for $\Re \, z^2 <1$ or $-1 < \Re\, z < 1$. Thus, the range of values of $z$ has changed, which means that the 
convergence of the series of the series on the rhs of Equivalence\ (\ref{fortyfoura}) will be affected by each value of $l$ 
or each series in the product. 

Let us examine the third series in the product, whose limit is $1/(1-z^3)$. In order to analyse this version of the geometric 
series, we write the limit value as
\begin{eqnarray} 
\frac{1}{1-z^3}= \frac{1}{(1-z) (1-z e^{2i \pi/3}) (1-z^{-2i \pi/3})} \;\;.
\label{six-one}\end{eqnarray}
Decomposing the rhs into partial fractions, we see that this version of geometric series is actually the sum of three
geometric series, each with a different limit. The first yields the standard geometric series discussed above. The second 
series has a limit of $1/(1-z \exp(2i\pi/3))$. In this case we replace $z$ by $z\exp(2i\pi/3)$ and continue with the same 
analysis. Then the second series is convergent for $\Re\, (z \exp(2i\pi/3)) <1$ or $y < (2-x)/\sqrt{3}$ when $z \!=\! x+ iy$, 
while it is divergent for $y > (2-x)/\sqrt{3}$. That is, the line $\Re\, z \! = \! 1$ separating the regions or planes of 
convergence and divergence has been rotated by $2\pi/3$ in a clockwise direction.  The ``left side" of the line representing 
where the series is convergent is now given by $y < (2-x)/\sqrt{3}$. On the other hand, it is divergent for $y > (2-x)/\sqrt{3}$ 
in which case the limit becomes the regularised value of the second series. The third series, whose limit is 
$1/(1-z \exp(-2i \pi/3)$, represents the opposite of the previous series. That is, the singularity at $z \!=\!1$ in the standard 
geometric series has now been rotated by $2\pi/3$ in an anti-clockwise direction. Hence, the third series is convergent for 
$\Re\, (z \exp(-2i\pi/3)) \!<\! 1$ or $y > -(x+2)/\sqrt{3}$ when $z \!=\! x+iy$. This represents the ``left" side, while the 
``right" side or where it is divergent is given by $y< -(x+2)/\sqrt{3}$. For these values of $z$ the limit represents the 
regularised value of the series. 

It is the intersection of the ``left" sides for the three series that represents the region of the complex plane for which
the third series in the product or $\sum_{k=0}^{\infty} z^{3k}$, is convergent. Outside this region the series is divergent. 
The intersection of the tangent lines yields an equilateral triangle, where the mid-points of the edges coincide with the 
three singular points of the component geometric series on the circle $|z| \!=\! 1$. Moreover, the unit disk of absolute 
convergence is circumscribed by this triangle. Those parts of the triangle not in the unit disk of absolute convergence 
represent the regions of the complex plane for which the third series in the product is conditionally convergent. In total 
they are significantly less than the corresponding region of the complex plane for the second series in the product, 
which we found is given by the region outside the unit disk of absolute convergence in the plane $-1 < \Re\,z < 1$.   

If we consider the fourth series in the product, i.e.\ $\sum_{k=0}^{\infty} z^{4k}$, then decomposing its limit into partial 
fractions yields four distinct geometric series with the singularities situated again on the unit circle, but at $\pm 1$ 
and $\exp( \pm i \pi/2)$. If we draw tangent lines through each of these singularities, then we find the common region 
or the intersection of their ``left" sides is now a square circumscribing the unit disk of absolute convergence. The 
singularities in the component geometric series appear at the mid-points of the square's sides. Outside the square, the 
fourth series in the product is divergent. Then the regularised value of the series is obtained by combining the limits 
of the component series. The regions inside the square, but outside the unit disk, represent the values of $z$ for 
which the fourth series in the product is conditionally convergent. As expected, these regions in total are less than 
either of the regions of conditional convergence for the second and third series in the product.

If we continue this analysis to the $l$-th series in the product, i.e.\ for $\sum_{k=0}^{\infty} z^{lk}$, then we find 
that the intersection of tangent lines yields an $l$-sided polygon that circumscribes the unit of disk of absolute 
convergence. For the values of $z$ outside the polygon the series will be divergent, while for those values of $z$ within 
the polygon, but outside of the unit disk, the series will be conditionally convergent. Furthermore, as higher values of 
$l$ are considered, the number of tangent lines not only increases, but also the total region of conditional convergence 
contracts. In the limit as $l \to \infty$ we will be left with the unit disk as the sole region where the series is 
(absolutely) convergent, while outside the disk the series is divergent. Since the product in Equivalence\ (\ref{fortyfoura}) 
includes all values of $l$, the $l \!=\! \infty$ limit by virtue of the fact that it possesses the smallest region of 
convergence in the complex plane determines the values of $z$, where the series on the rhs of Equivalence\ (\ref{fortyfoura}) 
is convergent. This means that the series or generating function, which we shall call from here on the partition number 
series, is equal to the product on the lhs only for $z$ situated in the unit disk. Outside the unit disk, the lhs represents 
the regularised value of the series and an equivalence symbol must be used instead of an equals sign. Finally, the circle 
$|z| \!=\! 1$ represents a ring of singularity separating the divergent values of the partition number series from the 
absolutely convergent values. In the case of the standard geometric series, there is only one point where the absolutely 
convergent region is separated from the divergent region, namely the singularity situated at $z \!=\!1$. Elsewhere, the line 
$\Re \, z \!=\! 1$ separates conditionally convergent values from divergent values. As mentioned above, the limit for the 
geometric series is indeterminate along the line, but is assigned a regularised value of $1/(1-z)$. However, for $z \!=\! 1$ 
the regularised value yields infinity. Hence, it can be seen that there is a difference between separating absolutely convergent 
values from divergent values and separating conditionally convergent values from divergent values. This completes the proof 
of the theorem.     

For $|z| \!<\! 1$, we can invert Equivalence\ (\ref{fortyfoura}), thereby obtaining
\begin{eqnarray}
\frac{1}{P(z)} = \prod_{k=1}^{\infty} \left( 1-z^k \right) = \frac{1} {1+ \sum_{k=1}^{\infty} p(k) z^k} \;\;. 
\label{six-two}\end{eqnarray}
Since the above product produces a power series that is valid for all values of $z$, we can write it as 
\begin{eqnarray}
(z;z)_{\infty}=\prod_{k=1}^{\infty} \left( 1-z^k \right) = 1+ \sum_{k=1}^{\infty} q(k) \, z^k \;\;.
\label{six-four}\end{eqnarray}
The leftmost expression is a special case of the q-Pochhammer symbol \cite{wik11a}, which is defined as
\begin{eqnarray}
(a;z)_n= \prod_{k=0}^{n-1} \left( 1 -a z^k \right) \;\;.
\label{six-foura}\end{eqnarray}

According to Knuth \cite{knu005}, it was Euler, who was the first to discover that much cancellation occurs when multiplying
the various terms in the infinite product given in Eq.\ (\ref{six-foura}). Specifically, he found that
\begin{eqnarray}
\prod_{m=1}^{\infty} \left( 1-z^m \right) = 1+ \sum_{k=1}^{\infty} (-1)^k \left( z^{(3k^2-k)/2}
+z^{(3k^2+k)/2}  \right) \;\;. 
\label{six-fourb}\end{eqnarray}
Therefore, comparing the above result with rhs of Eq.\ (\ref{six-four}) we that the $q(k)$ are frequently equal to zero and
when non-zero are either equal to 1 or -1. The values of $k$ for which the $q(k)$ vanish are known today as the pentagonal 
numbers \cite{wei011a}. They are themselves a particular case of a broader class of numbers known as the figurate or figural 
numbers \cite{wei011b,wik11b}. By applying Theorem\ 1 to Eq.\ (\ref{six-two}), we see that the coefficients of the outer series, 
viz.\ $q_k$, are equal to $(-1)^k$, while those for the inner series or the $p_k$ are equal to the partition function $p(k)$ 
for $k \geq 1$ and zero for $k \!=\! 0$. Therefore, with the aid of Eq.\ (\ref{five}) we arrive at
\begin{eqnarray}
q(k)= L_{P,k} \left[(-1)^{N_k}\, N_k!\prod_{i=1}^k \frac{p(i)^{n_i}}{n_i!} \right] \;\;.
\label{six-five}\end{eqnarray}
As a consequence of the fact that the $q(k)$ are non-zero when $k \!=\! (3j^2 \pm j)/2$, Eq.\ (\ref{six-five}) can also be written as 
\begin{eqnarray}
L_{P,k} \left[(-1)^{N_k}\, N_k!\prod_{i=1}^k \frac{p(i)^{n_i}}{n_i!} \right] = \begin{cases} (-1)^j, & \;\; k=(3j^2 \pm j)/2,  \\
0, & \;\; {\rm otherwise} \;\;. \end{cases}
\label{six-six}\end{eqnarray}

Although they do not give the actual number of discrete partitions, we shall refer to the $q(k)$ as the discrete 
partition numbers. Shortly, we shall see how these coefficients are related to the number of discrete partitions.
They do, however, have an interesting connection with the partition function, which follows when both power series on 
the rhs's of Eqs.\ (\ref{six-two}) and (\ref{six-four}) are multiplied by one another. Then we find that  
\begin{eqnarray}
\sum_{k=0}^{\infty} z^k \sum_{j=0}^k p(j) \, q(k-j) =1 \;\;,
\label{six-eight}\end{eqnarray}
where $p(0) \!=\! q(0) \!=\! 1$. Again, since $z$ is fairly arbitrary, like powers of $z$ can be equated on both sides of 
the above equation. For $k \geq 1$, we obtain the following recurrence relation:
\begin{eqnarray}
\sum_{j=0}^k p(j) \, q(k-j) =0 \;\;.
\label{six-nine}\end{eqnarray}
This is simply Euler's recurrence relation for the partition function, which is given by Eq.\ (20) in Ref.\ \cite{knu005}. 
Occasionally, it is referred to as MacMahon's recurrence relation as in Eq.\ (20) of Ref.\ \cite{wei011c}. 
Because most of the discrete partition numbers or $q(k)$ vanish, it means that only a few of the previous values of the 
partition function are required to evaluate the latest value of the partition function. 

Eq.\ (\ref{six-five}) is an interesting result, but perhaps, not very practical for determining the discrete partition
numbers, when it is realised that the partition function or $p(k)$ grows exponentially. We can, however, use the 
partition method for a power series expansion to derive another result for the discrete partition numbers. First, 
assuming that $|z| \!<\! 1$, we write $1/P(z)$ as
\begin{eqnarray}
\frac{1}{P(z)} = \exp \Bigl( \sum_{m=1}^{\infty} \ln \bigl( 1- z^m \, \bigr) \Bigr) =\exp \Bigl( - \sum_{m=1}^{\infty}
\sum_{j=1}^{\infty} \bigl( z^m \bigr)^{j}/j \Bigr) \;\;.
\label{six-ten}\end{eqnarray} 
The Taylor/Maclaurin series for the logarithm has been introduced in this result, since it is absolutely convergent 
for $|z| \!<\! 1$ \cite{kow10}. The double sum can be expressed as a single sum by realising that the coefficients of 
$z$ can be expressed as a sum over divisors or factors of the power \cite{wei011}. Then we arrive at
\begin{eqnarray}
P(z)^{-1} = 1 + \sum_{k=1}^{\infty} \frac{(-1)^k}{k!}\Bigl( z+\frac{3z^2}{2}+ \frac{4z^3}{3}+\frac{7z^4}{4}
+ \cdots + \gamma_j z^j+ \cdots \Bigr)^k \, ,
\label{six-eleven}\end{eqnarray} 
where $\gamma_j= \sum_{d|j} d/j$ and $d$ represents a divisor of $j$. That is, the sum is only over the divisors of 
$j$. Some values of the $\gamma_j$ are: $\gamma_1 \!=\! 1$, $\gamma_2 \!=\! 3/2$, $\gamma_3 \!=\!4/3$, 
$\gamma_4 \!=\! 7/4$, and $\gamma_5 \!=\! 6/5$. More explicitly, we find that $\gamma_6 = 1/6+1/3+1/2+1=2$, while for 
the case of $j \!=\! \rho^m$, where $\rho$ is a prime number, the sum yields
\begin{eqnarray}
\gamma_{\rho^m} = \frac{(1- 1/\rho^{m+1})}{(1-1/\rho)}  \;\;.
\label{six-twelve}\end{eqnarray}
Eq.\ (\ref{six-twelve}) can be derived simply by using the limit for the geometric series. Furthermore, Eq.\
(\ref{six-eleven}) also means that
\begin{eqnarray}
P(z) = 1 +  \sum_{k=1}^{\infty} \frac{1}{k!}\Bigl( z+\frac{3z^2}{2}+ \frac{4z^3}{3}+\frac{7z^4}{4}
+ \cdots + \gamma_j z^j+ \cdots \Bigr)^k \, .
\label{six-twelvea}\end{eqnarray} 
Therefore, we have obtained alternative representations for the generating functions of both $P(z)$ and its inverse.

Now we are in a position to apply Theorem\ 1 to Eq.\ (\ref{six-eleven}), whereupon we see that the coefficients 
of the inner series $p_k$ equal $ \gamma_k$ for $k \! \geq \! 1$, while $p_0 \!=\! 0$. On the other hand, the 
coefficients of the outer series $q_k$ are equal to $(-1)^k/k!$. Hence, from Eq.\ (\ref{five}) the discrete
partition numbers are given by
\begin{eqnarray}
q(k) = L_{P,k} \Bigl[ (-1)^{N_k} \prod_{i=1}^k \frac{\gamma_i^{n_i}}{n_i!} \Bigr] \;\;.
\label{six-thirteen}\end{eqnarray}
Hence, we have derived an alternative version of Eq.\ (\ref{six-five}). Moreover, whenever $k \!=\! (3j^2\pm j)/2$,
for $j$ an integer, this result can also be written as
\begin{eqnarray}
L_{P,(3j^2\pm j)/2} \Bigl[ (-1)^{N_{(3j^2 \pm j)/2}} \prod_{i=1}^{(3k^2 \pm j)/2} \frac{\gamma_i^{n_i}}{n_i!} \Bigr] 
=(-1)^j\;\;.
\label{six-thirteena}\end{eqnarray}
For all other values of $k$, the sum over all partitions in Eq.\ (\ref{six-thirteen}) vanishes. Therefore, to 
evaluate $q(6)$, we require all the $\gamma_j$ ranging from $j \!=\!1$ to 6, which have already been given above. 
Summing over the eleven partitions summing to 6 in Eq.\ (\ref{six-thirteen}) yields
\begin{align}
q(6) & =-2 + 6/5 + 21/8 - 7/8 + 16/18 - 2 + 4/18 
\nonumber\\
& - \;\; 27/48 +  9/16 - 3/48 + 1/6!=0 \;\;,
\label{six-fourteen}\end{align} 
which is indeed the value of this discrete partition number. By using these results, the reader can readily verify that 
$q(0) \!=\! 1$, $q(1) \!=\! -1$, $q(2) \!=\!-1$, $q(5) \!=\! 1$ and $q(3) \!= \! q(4) \!=\! 0$. Moreover, since $D_0$ in 
Theorem\ 1 is non-zero in this case due to the fact that $q(0) \!=\!1$, we can use Eq.\ (\ref{seven}) to determine the 
coefficients of the inverted power series expansion or the generating function for $P(k)$. Hence, we find that
\begin{eqnarray}
p(k)= L_{P,k} \Bigl[ (-1)^N N!\prod_{i=1}^k \frac{q(i)^{n_i}}{n_i!} \Bigr] \;\;.
\label{six-fifteen}\end{eqnarray}
This result, which represents the inverse of Eq.\ (\ref{six-five}), incorporates much redundancy since the $q(i)$ are 
only non-zero when for specific values of $i$. Consequently, both the sum over the partitions and the product are only 
non-zero for those values of $i$, which are of the form of $(3j^2 \!-\!j)/2$ or $(3j^2 \!+\! j)/2$, where $j$ is an integer 
ranging from 1 to $j_m=[(1+\sqrt{1+24k})/6]$. In addition, when the $q(i)$ are non-zero, they are only equal to unity 
in magnitude. Therefore, it is the factor of $N!$ that is responsible for the exponential increase in the partition 
function as $k$ increases, although this factor will often be countered by the $1/n_i!$ terms in the denominator of 
the product. For example, if we wish to determine $p(6)$, there will not be any sums over $n_3$, $n_4$ and $n_6$ in 
the above equation since $q(3)$, $q(4)$ and $q(6)$ vanish, while in the product, $n_3!$, $n_4!$ and $n_6!$ will equal 0! 
or unity. Furthermore, since $q(1)=q(2)=-1$ and $q(5)=1$, Eq.\ (\ref{six-fifteen}) yields 
\begin{align}
p(6) & = 1 + (-1)^5 (-1)^4 (-1)\, 5!/4! + (-1)^4 (-1)^2 (-1)^2 4!/(2! \cdot 2!)
\nonumber\\
& + (-1)^3 (-1)^3 \, 3!/3! + (-1)^2 (-1) 2! = 11  \;\;.
\label{six-sixteen}\end{align}
Hence, the contributions from the five relevant partitions are positive except for the last one, which represents the
contribution due to the partition \{1,5\}.

On the other hand, if we apply Theorem\ 1 to Eq.\ (\ref{six-twelvea}), then the only difference to the previous 
evaluation of the discrete partition numbers is that the coefficients of the outer series $q_k$ are now equal to $1/k!$. 
That is, the coefficients of the outer series are still equal to $\gamma_k$. Hence, we arrive at 
\begin{eqnarray}
p(k) = L_{P,k} \Bigl[ 1 \Bigr] = L_{P,k} \Bigl[ \prod_{i=1}^k \frac{\gamma_i^{n_i}}{n_i!} \Bigr] \;\;.
\label{six-sixteena}\end{eqnarray}
Therefore, we have an entirely different means of evaluating the partition function with the sum of the reciprocals
of the divisors of each element being the assigned values in the partition method for a power series expansion. 
Consequently, Eq.\ (\ref{six-nine}) becomes
\begin{eqnarray}
\sum_{j=0}^k  L_{P,k} \Bigl[ \prod_{i=1}^j \frac{\gamma_i^{n_i}}{n_i!} \Bigl] \,
L_{P,k} \Bigl[ (-1)^{N_{k-j}} \prod_{i=1}^{k-j} \frac{\gamma_i^{n_i}}{n_i!} \Bigl] =0 \;\;.
\label{six-sixteenb}\end{eqnarray}

In deriving these new results for the partition function $p(k)$ we have encountered a more complex example involving a specific
class of partitions than any of those studied in the previous section. In this case we only require partitions
whose elements are pentagonal numbers or of the form of $(3j \pm 1)j/2$, where $j$ is any integer lying between zero and 
$j_m$. In view of the importance of the Eq.\ (\ref{six-fifteen}), let us consider introducing the modifications to the 
program ${\bf partgen}$ presented in Sec.\ 2. As in the examples of the previous section, most of these modifications 
will be made to the function prototype ${\bf termgen}$.

The first modification is that we need to make is introduce the math library with the other header files. This is 
needed so that the floor function can be called to evaluate the maximum value of $j$ or $i_m$ when the value of 
$k$ or ${\it tot}$ is typed in by the user. Since $j_m$ is called in both the ${\bf main}$ and ${\bf termgen}$ 
prototypes, it must be declared as a global variable. Once these modifications are carried out, we can concentrate on 
the changes that are required for ${\bf termgen}$, which is displayed below:
\begin{lstlisting}{}
void termgen()
{
int freq,i,j,jval;

for (i=0;i<tot;i++){
              jval=0;
              freq=part[i];
              if (freq >0){
                   for (j=1;j<=j_m;j++){
                      if (i == ((3*j-1)*j-2)/2) jval=j;
                      if (i == ((3*j+1)*j-2)/2) jval=j;
                                       }
                   if ((jval==0) && (freq>0)) goto end;
                          }
                   }

printf("%ld: ",term++);
for(i=0;i<tot;i++){
            freq=part[i];
            if(freq) printf("%i(%i) ",freq,i+1);
                  }
printf("\n");
end:        ;
}

\end{lstlisting}

\begin{figure}
\begin{center}
\includegraphics{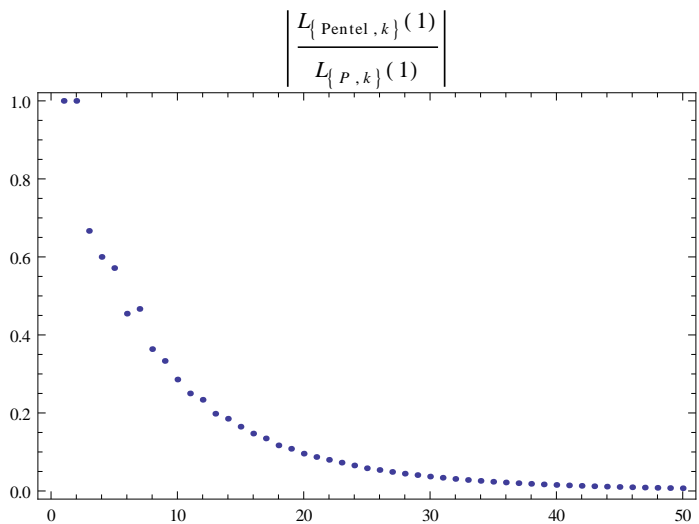}
\end{center}
\caption{The ratio of the number of partitions whose elements are pentagonal numbers to the total number of partitions versus $k$}
\label{figthreeb}
\end{figure}

Basically, this version of ${\bf termgen}$ is similar to those that appeared in the previous section. That is, 
before a partition is printed out, testing is done in the function prototype to see if each partition belongs 
or conforms to the particular set of partitions, which the user desires. In this case our aim is to print out 
only those partitions whose elements can be written in the pentagonal number forms of $(3j^2-j)/2$ and 
$(3j^2+j)/2$, where $j$ is an integer. This is accomplished by introducing another for loop in the function prototype, 
which evaluates the value of $j$ called ${\it jval}$, for the appropriate elements. If an element is not of the 
required form, then ${\it jval}$ remains zero. Otherwise, it is non-zero. If ${\it jval}$ is zero for an element, 
then the partition is examined to see if there any occurrences of the element by checking the variable ${\it freq}$. 
If it is greater than zero, then the entire partition is discarded by the goto statement. This test is 
carried out on all elements in the partition. The procedure is then applied to all partitions summing to ${\it tot}$. 
Only those partitions whose elements are of the required form are printed out by the code. As an example, when 
$k \!=\!6$, the output for this program called ${\bf partfn}$ is:\newline
1: 1(1) 1(5) \newline
2: 4(1) 1(2) \newline
3: 6(1) \newline
4: 2(1) 2(2) \newline
5: 3(2) \newline
The above partitions represent the five that contributed to the calculation of $p(6)$ given above. Interestingly,
when the code is run for partitions summing to 100, it only prints out 42205 partitions which represents 
about 0.02 percent of the total number of partitions or $p(100)$. Moroever, the number of partitions summing to
$k$ whose elements are pentagonal numbers or $L_{Pentel,k}[1]$ starts off greater than the number of 
discrete partitions for the same value of $k$, but then drops away when $k \!> \! 17$. Fig.\ \ref{figthreeb}
displays the ratio of $L_{Pentel,k}[1]$ to the number of partitions summing to $k$ for $k \leq 50$. Once again, 
it is monotonically decreasing for $k$ beyond a certain value, which in this case is about 8. 

It was also stated at the end of the previous section that the partition method for a power series expansion 
can be adapted to handle situations where only a subset of the partitions summing to a particular value are 
required as in the above example. We shall demonstrate this by modifying the code in Sec.\ 4 so that the 
partition function $p(k)$ can be evaluated via Eq.\ (\ref{six-fifteen}). The resulting code called ${\bf pfn}$ 
is presented in its entirety in the appendix. Basically, it expresses the partition function in a symbolic form 
where the final values can be evaluated by introducing the output into Mathematica. As expected, ${\bf termgen}$ 
has the same for loop presented above. The interesting feature about the code is that the for loop
appears twice at the beginning of ${\bf termgen}$. This is necessary because the first step prints out 
$p(k):=$ before considering the single element partition of \{$k$\} in the first branch of the if statement.
Clearly, the single element partition may not be of the required form--- hence, the first appearance of the
for loop. All the other partitions are processed by the second branch of the if statement, which means that 
we require the for loop in this part of the program in order to determine which are of the required form.

The code is also capable of expressing the final symbolic form for the partition function in two different forms. 
In the first form a printf statement prints the final form for $p(k)$ in terms of the $q(i)$ as they appear 
in Eq.\ (\ref{six-fifteen}). This statement has been commented out in favour of the version of the code appearing 
in the appendix. In the second version of the code the value of each $q(i)$ is evaluated. That is, the code prints 
out -1 to the power of ${\it jval}$ for the elements in a suitable or candidate partition. In order to accomplish this, 
the for loop mentioned in the previous paragraph has had to be introduced into the latter part of ${\bf termgen}$ again. 
Finally, the quantity $q[N] a^N$ appearing in the output of the code discussed in Sec.\ 4 has been replaced 
$(-1)^{\wedge}(N)$ in ${\bf pfn}$. Therefore, if the version of ${\bf pfn}$ in the appendix is run for ${\it tot}=6$, 
then the following output is produced:\newline   
p[6]:= ((-1)$^{\wedge}$(1)) ((-1)$^{\wedge}$(2)) (-1)$^{\wedge}$(2) 2! + ((-1)$^{\wedge}$(1))$^{\wedge}$(4) 
((-1)$^{\wedge}$(1)) (-1)$^{\wedge}$(5) 5!/4! + 
((-1)$^{\wedge}$(1))$^{\wedge}$(6) (-1)$^{\wedge}$(6) + ((-1)$^{\wedge}$(1))$^{\wedge}$(2) ((-1)$^{\wedge}$(1))$^{\wedge}$(2) 
(-1)$^{\wedge}$(4) 4!/(2! 2!) + ((-1)$^{\wedge}$(1))$^{\wedge}$(3) (-1)$^{\wedge}$(3) \newline 
Time taken to compute the coefficient is 0.000000 seconds. \newline
This output can be imported into Mathematica whereupon it gives the correct value of 11 for $p(6)$. When the code is 
run for ${\it tot}$ equal to 100, it takes 34 seconds to compute the symbolic form for $p(100)$ or $L_{p,100}[1]$ and only 
0.27 seconds to produce the value of $190\,569\,292$ in Mathematica on the same Sony VAIO laptop mentioned in previous 
sections. However, whilst this is not too bad, the fastest method of obtaining $p(k)$ or $L_{p,k}[1]$ from the various 
programs considered in this work is to comment out ${\bf termgen}$ and introduce the statement 
\begin{lstlisting}{}
	 term++; 
\end{lstlisting} 
immediately below in ${\bf partgen}$ of Sec.\ 2. Then ${\it term}$ needs to be initialised to zero and a printf 
statement introduced into ${\bf main}$. When the resulting program is run for ${\it tot}$ equal to 100, it only 
takes 3 seconds to compute $p(100)$. Nevertheless, both methods for computing the partition function will slow down 
dramatically as ${\it tot}$ continues to increase due to combinatorial explosion. This is where either Eq.\ (\ref{six-nine})
or even the Hardy-Ramanujan-Rademacher formula should be used. In fact, we may write
\begin{eqnarray}
\lim_{k \to \infty} L_{P,k} \Bigl[ 1 \Bigr] \to
\frac{1}{4\sqrt{3} \,k} \; \exp( \sqrt{2k/3}\, \pi) \;\;,
\label{six-seventeen}\end{eqnarray} 
and 
\begin{eqnarray}
\lim_{k \to \infty} L_{P,k} \Bigl[ (-1)^N N!\prod_{i=1}^k \frac{q(i)^{n_i}}{n_i!} \Bigr] \to
\frac{1}{4\sqrt{3} \,k} \; \exp( \sqrt{2k/3}\, \pi) \;\;.
\label{six-eighteen}\end{eqnarray} 
In addition, the arrow symbols in the above results can be replaced by equals signs.

\section{Generalisation of the Inverse of ${\bf P(z)}$}
In the previous section we demonstrated how the partition method for a power series expansion can be used to derive 
different forms for the generating functions of the product $P(z)$ and its inverse or reciprocal. In the case of $P(z)$ 
two different forms for the coefficients of the generating function or the partition function $p(k)$ were obtained in 
terms of the partition operator acting on different arguments. The first given by Eq.\ (\ref{six-fifteen}) involved 
the discrete partition numbers $q(k)$, which represent the coefficients of the generating function for the inverse 
of $P(z)$, while the second given by Eq.\ (\ref{six-sixteena}) involved a sum over the inverses of the divisors for 
each element $i$ in a partition. For the inverse of $P(z)$ two different forms for the discrete partition numbers 
were also obtained, but now the partition operator was found to act on either the partition function as in the first 
form given by Eq.\ (\ref{six-five}) or in the second form given by Eq.\ (\ref{six-thirteen}) the same sum over the 
inverses of the divisors with an extra phase factor of $(-1)^{N_k}$, where $N_k$ represents the number of elements 
in each partition summing to $k$. This extra phase factor is responsible for ensuring that the discrete partition 
numbers equal either $\pm 1$ when $k$ is a pentagonal number or zero, otherwise. Because the phase factor does not 
appear in Eqs.\ (\ref{six-fifteen}) and (\ref{six-sixteena}), the partition function experiences exponential growth 
as can be seen from the Hardy-Ramanujan-Rademacher formula given above. Nevertheless, it was possible to derive all 
these results because the coefficients of the powers of $z$ in the product $P(z)$ are simple, namely equal to -1. 
In this section we aim to generalise $P(z)$ to the situation where the coefficients of $z^k$ are now equal to $C_k$. 
We begin by studying the inverse of $P(z)$. As a consequence, we arrive at the following theorem.

{\bfseries Theorem\ 3}. The infinite product defined by 
\begin{eqnarray}
H(z) = \prod_{i=1}^{\infty} \left( 1+ C_i z^i \right) \;\;,
\label{six-nineteen}\end{eqnarray}
can be written as a power series expansion or generating function of the form: 
\begin{eqnarray}
H(z) = 1+ \sum_{k=1}^{\infty} h_k z^k \;\;,
\label{six-twenty}\end{eqnarray}
where the coefficients $h_k$ can be expressed in terms of the discrete partition operator defined by Eq.\ (\ref{five-five}) 
and are given by  
\begin{eqnarray}
h_k =L_{DP,k} \Bigl[ \prod_{i=1}^k C_i^{n_i} \Bigr] \;\;.
\label{six-twentyone}\end{eqnarray}

{\bf Proof}. If we multiply out the lowest order terms in $z$ in the product given in Eq.\ 
(\ref{six-nineteen}), then we obtain a power series expansion for $H(z)$ in a similar manner to the proof
of Theorem\ 1. The zeroth order term in the resulting power series is unity, while the first order term 
becomes $C_1 z$. We then find that the second order term in the power series becomes $C_2 z^2$. Once we go 
beyond second order, however, the coefficients become more complex to evaluate due to an ever-increasing 
number of terms appearing in them. For example, the third, fourth and fifth order coefficients are respectively 
equal to $C_3+C_1 C_2$, $C_4+ C_1 C_3$, and $C_5+C_1 C_4+C_2 C_3$. In fact, on close inspection we find that the 
coefficient of $z^k$ in the power series depends on the number of discrete partitions summing to $k$. For example,
the third and fourth order coefficients are composed of only two terms because there are only two discrete
partitions for $k \!=\! 3$ and $k \!=\! 4$, while the fifth order term is composed of three terms due to the three
discrete partitions summing to 5, viz.\ \{5\}, \{1,4\} and \{2,3\}. Therefore, instead of summing over all the partitions 
summing to $k$ as we did in Theorem\ 1, we need only sum over the subset comprising the discrete partitions summing to
$k$, which we have already seen is significantly less than $p(k)$. Hence, the sum over all partitions in Eq.\ 
(\ref{five}) simplifies drastically with all the frequencies lying between zero and unity, not between zero and 
$[k/i]$ for each element $i$. Furthermore, since there is no outer series in the expansion of the product, we can 
drop $q^N$ and $a^N$ in Eq.\ (\ref{five}). In addition, there is no multinomial factor. Consequently, we find 
that according to Eq.\ (\ref{five}) the coefficients $h_k$ in the power series for $H(z)$ are given by
\begin{eqnarray}
h_k = \sum_{\scriptstyle n_1, n_2,n_3, \dots,n_k=0 \atop{\sum_{i=1}^k in_i =k}}^{1,1,1,\dots,1} 
\prod_{i=1}^k C_i^{n_i} \;\;.
\label{six-twentyonea}\end{eqnarray}
The sum over partitions in the above result is simply the discrete partition operator as defined by Eq.\ 
(\ref{five-five}). When this is introduced into Eq.\ (\ref{six-twentyonea}), we arrive at Eq.\ 
(\ref{six-twentyone}). This completes the proof of the theorem. 

It should also be mentioned that we can invert the rhs of Eq.\ (\ref{six-twenty}) and apply Theorem\ 1 again. Then
we obtain a power series expansion for $1/H(z)$. As a result, we find that
\begin{eqnarray}
H_k = L_{P,k} \left[ (-1)^{N_k} \,N_k! \prod_{i=1}^k \frac{h_i^{n_i}}{n_i!} \right] \;\;,
\label{six-twentyoneb}\end{eqnarray} 
while from Eq.\ (\ref{seven}), we obtain
\begin{eqnarray}
h_k = L_{DP,k} \left[ \prod_{i=1}^k C_i^{N_i}\right]=L_{P,k} \left[ (-1)^{N_k} \,N_k! \prod_{i=1}^k \frac{H_i^{n_i}}{n_i!} \right] \;\;.
\label{six-twentyonec}\end{eqnarray} 

Since Theorem\ 3 is quite general, it means conversely that any power series expansion can be expressed as an 
infinite product of the form given by the rhs of Eq.\ (\ref{six-nineteen}). For example, as discussed on
p.\ 111 of Ref.\ \cite{knu005} the geometric series can be represented by the following product:
\begin{eqnarray}
\sum_{k=0}^{\infty} z^k = \prod_{i=0}^{\infty} \left( 1+ z^{2^i} \right) \;\;.
\label{six-twentytwo}\end{eqnarray}
In this example $C_{2^i} \!=\! 1$, while the other $C_i$ simply vanish. This means that the product of
$\prod_{i=1}^k C_i^{n_i}$ in the $h_k$ is either zero or unity. In addition, since the coefficients of the 
geometric series are equal to unity, we have $L_{DP,k}[\prod_{i=1}^k C_i^{n_i}] \!=\! 1$. Since all the terms
in the discrete partition operator are either equal to zero and unity, this means that there can only be one 
discrete partition where all the elements are of the form $2^j$ with each value of $j$ lying between zero and 
$[\log k/\log 2]$. For $k \!=\! 5$ and $k \!=\! 7$, these discrete partitions are respectively \{1,4\} and 
\{1,2,4\}, while for $k \!=\! 11$ it is \{1,2,8\}. 

Another example of a well-known power series that can be expressed as an infinite product is the exponential 
function, which can be expressed as
\begin{align}
e^y & = (1+y)(1+y^2/2)(1-y^3/3)(1+3y^4/8)(1-y^5/5)(1+13 y^6/72) 
\nonumber\\
&\times \;\; (1-y^7/7) (1+27y^8/128) (1-8y^9/81) (1+91 y^{10}/800) \cdots \;\; .
\label{six-twentytwoa}\end{align}
From Eq.\ (\ref{six-twentytwoa}) we find that $C_1 \!=\! 1$, while the other $C_i$ can be determined by the 
following recurrence relation involving the divisors of $i$:
\begin{eqnarray}
\sum_{d|i} \frac{(-1)^{d+1}}{d} \; C^d_{i/d} =0 \;\;.
\label{six-twentytwob}\end{eqnarray}
This means that whenever $i$ is a prime number greater than 2, say $p$, we find that $C_p=-1/p$. Eq.\ 
(\ref{six-twentytwob}) is left as an exercise for the reader.

We can derive other interesting results by taking the logarithm of both Eqs.\ (\ref{six-nineteen}) and (\ref{six-twenty}), 
which yields
\begin{eqnarray}
\ln \Bigl( 1+ \sum_{k=1}^{\infty} h_k \, z^k \Bigr)= \sum_{i=1}^{\infty} \ln \Bigl( 1+C_i z^i \Bigr) \;\;. 
\label{six-twentytwoc}\end{eqnarray}
By introducing the Taylor/Maclaurin series expansion for logarithm into the rhs, we obtain
\begin{eqnarray}
\ln \Bigl( 1+ \sum_{k=1}^{\infty} h_k \, z^k \Bigr)= \sum_{i=0}^{\infty} \sum_{j=1}^{\infty} 
\frac{(-1)^{j+1}}{j}\; C_i^j z^{ij} \;\;, 
\label{six-twentytwod}\end{eqnarray}
where we have now assumed that $|z| \!<\! 1$. We now apply Theorem\ 1 to the left hand side of the above result.
The inner series coefficients are given by $p_k \!=\! h_k$, while the outer series coefficients are given by
$q_k \!=\! (-1)^{k+1}/k$. By equating like powers of $z$ in the resulting expression, we arrive at
\begin{eqnarray}
L_{P,k} \Bigl[ (-1)^{N_k }(N_k-1)! \prod_{i=1}^k \frac{h_i^{n_i}}{n_i !} \Bigr] = \sum_{i|k} \frac{(-1)^{i}}{i}\;
C_{k/i}^i \;\;.
\label{six-twentytwoe}\end{eqnarray}
For the particular case of the geometric series, viz.\  Eq.\ (\ref{six-twentytwo}), $h_k \!=\! 1$ and 
$C_{2^j} \!=\! 1$, while the other values of the $C_k$ vanish. Introducing these results into the above equation 
produces
\begin{eqnarray}
L_{P,2^{j}} \Bigl[ (-1)^{N_{2^{j}} }(N_{2^{j}}-1)! \prod_{i=1}^{2^{j}} \frac{1}{n_i !} \Bigr] = \Bigl( -
\frac{1}{2} \Bigr)^j \;\;.
\label{six-twentytwof}\end{eqnarray}
 
If we put the $C_i$ equal to unity in Eq.\ (\ref{six-nineteen}), then we find that
\begin{eqnarray}
\prod_{k=1}^{\infty} \left( 1 + z^k \right) =1+ \sum_{k=1}^{\infty} L_{DP,k} \Bigl[ 1 \Bigr] \, z^k \;\;.
\label{six-twentythree}\end{eqnarray}
In another words, the coefficients of the power series expansion for the product are given by the 
number of discrete partitions summing to $k$. Therefore, from Eq.\ (\ref{six-twentyone}) we have
\begin{eqnarray}
h_k = L_{DP,k} \Bigl[ 1 \Bigr] \;\;.
\label{six-twentythree1}\end{eqnarray}
On the other hand, if we put the $C_i = -1$, then we obtain 
\begin{eqnarray}
q(k)= L_{DP,k} \Bigl[ (-1)^{N_k} \Bigr] \;\;.
\label{six-twentyfour}\end{eqnarray}       
Prior to Eq.\ (\ref{six-five}) it was mentioned that the $q(k)$ or discrete partition numbers are equal 
to $(-1)^j$ whenever $k$ is a pentagonal number \cite{wei011a} or is equal to $j(3j\pm 1)/2$. Hence, Eq.\ 
(\ref{six-twentyfour}) can be expressed as
\begin{eqnarray}
L_{DP,k} \Bigl[ (-1)^{N_k} \Bigr] = \begin{cases} (-1)^j, \;\; & \;\; {\rm for} \; k=j(3j\pm 1)/2,  \\
0, \;\; & \;\; {\rm otherwise} \;\;. \end{cases}
\label{six-twentyfoura}\end{eqnarray}       
The above result tells us that the number of discrete partitions with an odd number of elements
is equal to the number of discrete partitions with an even number of elements when $k$ is not a 
pentagonal number. It also tells us that the discrete partition numbers $q(i)$ come in pairs of either 
1 or -1, the former corresponding to a pentagonal number derived from an even number, i.e. by setting 
$j$ equal an even number and the latter to a pentagonal number derived by an odd number. A value of 
1 in the above result or for $j$ an even number, means that the number of discrete partitions with an 
even number of elements is one more than the number of discrete partitions with an odd number of elements. 
A value of -1 or $j$ odd, represents the opposite situation. A proof of this result based on Ferrers
diagrams appears in Ch.\ 1 of Ref.\ \cite{and03}.   

Now if we put $C_i =-\omega$ in Eq.\ (\ref{six-twenty}), then the power series for the ensuing product 
becomes
\begin{eqnarray}
\prod_{k=1}^{\infty} \left( 1- \omega z^k \right) = \frac{(-\omega;z)_{\infty}}{(1-\omega)}= 1+ 
\sum_{k=1}^{\infty} q(k,-\omega) z^k \;\;,
\label{six-twentyfive}\end{eqnarray}
where $(-\omega;z)_{\infty}$ is again the q-Pochhammer symbol presented in Eq.\ (\ref{six-foura}) and 
\begin{eqnarray}
q(k,-\omega) = L_{DP,k} \Bigl[ (-\omega)^{N_k}  \Bigr] \;\;.
\label{six-twentysix}\end{eqnarray}
Thus we see that $q(k,-1) \!=\! q(k)$ or the discrete partition numbers. Note also that the power of $\omega$ gives 
the number of elements in each discrete partition. By summing over all partitions each power of $\omega$, say 
$\omega^n$, in $q(k,\omega)$ gives the total number of discrete partitions with $n$ elements in them. A similar 
situation will arise when we study the inverse of Eq.\ (\ref{six-twentyfive}). Moreover, by putting $C_i \!=\! 
\omega$ in Eq.\ (\ref{six-twenty}), we obtain
\begin{eqnarray}
Q(z,\omega)= \prod_{k=1}^{\infty} \left( 1+ \omega z^k \right) = \frac{(\omega;z)_{\infty}}{(1+\omega)}= 1+ 
\sum_{k=1}^{\infty} q(k,\omega) z^k \;\;,
\label{six-twentysixa}\end{eqnarray}
where 
\begin{eqnarray}
q(k,\omega) = L_{DP,k} \Bigl[ \omega^{N_k}  \Bigr] \;\;.
\label{six-twentysixb}\end{eqnarray}
The above polynomials will be referred to as the discrete partition polynomials. From Eq.\ (\ref{six-twentythree})
we have already seen that the $\omega \!=\!1 $ case gives the coefficients that represent the total number of 
discrete partitions summing to $k$, viz. $L_{DP,k}[1]$, while from the above the discrete partition numbers are 
given by $q(k)= L_{DP,k}[(-1)^{N_k}]$. Hence, we see that the difference between the discrete partition numbers 
and the number of discrete partitions for a particular value of $k$ is that for the latter the discrete partition 
operator acts on unity when scanning the discrete partitions summing to $k$, while for the former it acts on $(-1)^{N_k}$, 
where $N_k$ represents the number of elements in each discrete partition. That is, the difference is caused again by a
phase difference in the number of elements in the discrete partitions summing to $k$.

\begin{table}
\begin{center}
\begin{tabular}{|c|l|l|} \hline
$k$ &  $q(k,\omega)$ & $p(k,\omega)$ \\ \hline
$0$ & $ 1$ & $1$ \\  
$1$ &  $  \omega $ & $ \omega$  \\ 
$2$ & $  \omega $ & $ \omega^2 + \omega$ \\ 
$3$ & $\omega^2 + \omega  $ & $\omega^3 + \omega^2 +\omega$ \\ 
$4$ & $\omega^2 + \omega $ & $\omega^4 +\omega^3+2 \omega^2+\omega$ \\ 
$5$ & $2 \omega^2 +\omega $ & $\omega^5 +\omega^4+2\omega^3+ 2\omega^2+\omega$ \\ 
$6$ & $\omega^3 + 2 \omega^2 + \omega $ & $ \omega^6 +\omega^5+2\omega^4+ 3 \omega^3+3\omega^2+\omega$ \\ 
$7$ & $ \omega^3 + 3 \omega^2 +\omega $ & $ \omega^7 +\omega^6+2 \omega^5+3 \omega^4+4 \omega^3+3 \omega^2+ \omega$ \\ 
$8$ & $2 \omega^3 + 3 \omega^2 +   \omega $ & $ \omega^8+\omega^7+2 \omega^6+3 \omega^5+5 \omega^4+
5 \omega^3+4 \omega^2+\omega$\\
$9 $ & $3 \omega^3+ 4 \omega^2 +\omega $ & $\omega^9 +\omega^8 +2 \omega^7 
+3 \omega^6 +5 \omega^5 + 6\omega^4+7\omega^3+4 \omega^2$ \\
$$ & $$ & $+\omega$ \\
$10$ & $\omega^4+ 4 \omega^3+ 4 \omega^2+ \omega$ & $ \omega^{10}+\omega^9 +2\omega^8+3 \omega^7+5\omega^6+7 \omega^5
+9 \omega^4+8 \omega^3$\\
$$ & $$ & $ +5 \omega^2+ \omega$ \\ \hline
\end{tabular}
\end{center}
\caption{Discrete partition polynomials $q(k,\omega)$ and partition function polynomials $p(k,\omega)$.}
\label{table2a}
\end{table}

The discrete partition polynomials up to $k \!=\! 10$ are displayed in the second column of Table\ \ref{table2a}. 
As expected, for $\omega \!=\! -1$ these results reduce to the discrete partition numbers $q(k)$, while for 
$\omega \!=\!1$ they yield the number of discrete partitions. From the table it can be seen that the discrete partition 
polynomials are polynomials of degree $n$, where $n(n+1)/2 \!\leq\! k \!<\! (n+1)(n+2)/2$ since the partition 
with the most discrete elements when $k \!=\! n(n+1)/2$ is \{1,2,3,\dots,n\}. The lowest order term in $\omega$ 
corresponds to the only one-element partition, viz. \{k\}. If we run the program ${\bf dispart}$, which is displayed 
in the appendix, then we find that there are 4 discrete partitions summing to 6, which are \{6\}, \{1,5\}, \{2,4\} 
and the self-conjugate \{1,2,3\}. By referring to Fig.\ \ref{figone}, we see that the first partition is just one 
branch from the seed number, the second and third partitions are two branches away and the third is three branches 
away. Hence, we arrive at $q(6,\omega) \!=\! \omega +2 \omega^2 +\omega^3$, where the magnitude of the 
coefficients of $\omega^k$ represent the number of distinct partitions with $k$ elements. That is, if the number of 
distinct elements summing to $k$ with $i$ elements in them are given by $q_i(k)$, then the polynomials can be expressed 
as
\begin{eqnarray}
q(k,\omega) = \sum_{i=1}^{n} q_i(k) \omega^i \;\;,
\label{six-twentyseven}\end{eqnarray}     
where $n(n+1)/2 \!\leq\! k \!<\! (n+1)(n+2)/2$. Hence, if we let $q(k,1) \!=\! L_{DP,k} \bigl[ 1\bigr]$ or the 
number of distinct partitions summing to $k$. then we obtain the trivial equation of $q(k,1)=\sum_{i=1}^{n} q_i(k)$. 
Moreover, the lower bound on $k$ gives us a limit as to the maximum number of elements that can appear in a 
discrete partition, which is given by
\begin{eqnarray}
i_{max} = [(\sqrt{8k+1} -1)/2]\;\;.
\label{six-twentysevena}\end{eqnarray}

By fixing the number of elements to $i$ in the discrete partition operator so that it becomes 
\begin{eqnarray}
L_{DP,k,i}\Bigl[ \cdot \Bigr]= \sum_{\scriptstyle n_1,n_2\dots,n_{k}=0 
\atop{\sum_{j=1}^{k} j n_j=k,\; \;\;\sum_{j=1}^{k} \! n_j=i}}^{1,1,\dots,1} \;\;,
\label{six-twentyeight}\end{eqnarray}
we arrive at $L_{DP,k,i}[1]= q_i(k)$. This means that we need to input a second value into the program,
which represents the number of elements the user desires. If this is set equal to a global variable called
${\it numparts}$, then the only changes to be made to ${\bf dispart}$ are: (1) introduce a local 
variable ${\it sumpart}$. which adds all the values of ${\it freq}$ in the first for loop 
and (2) insert the following ${\it if}$ statement before anything is printed out:
\begin{lstlisting}  
	if(sumparts != numparts) goto end;
\end{lstlisting}
When these modifications are implemented and the resulting code run for the discrete partitions summing
to 100 with the number of elements, i.e.\ ${\it numparts}$, set equal to 5, one finds that there are
25,337 discrete partitions beginning with \{1,2,3,4,90\} and ending with \{18,19,20,21,22\}. That is, $q_5(100)$ 
or $L_{DP,100,5}[1]$ is equal to $25\,337$. According to Eq.\ (\ref{six-twentysevena}), the maximum number of 
elements in the discrete partitions summing to 100 is 13. When the code is run for ${\it numparts}$ set 
equal to 13, 30$(=q_{13}(100))$ partitions are printed out beginning with \{1,2,3,4,...,12,22\} and ending 
with \{1,2,3,4,6,7,...,13,14\}. Running the code for higher values of ${\it numparts}$ with $k \!=\! 100$ 
does not result in any partitions being printed out. Hence, $q_{i}(100) \!=\! 0$ for $i \!>\! 13$, which 
is consistent with Eq.\ (\ref{six-twentysevena}). 

From Eq.\ (\ref{six-twentyfive}) we have 
\begin{eqnarray}
\prod_{k=1}^{\infty} \left( 1- \omega^2 z^{2k} \right) = \frac{(-\omega;z^2)_{\infty}}{(1-\omega^2)}= 1+ 
\sum_{k=1}^{\infty} q(k,-\omega^2) z^{2k} \;\;.
\label{six-twentyeighta}\end{eqnarray}
The product on the lhs of this result can also be written as 
\begin{eqnarray}
\prod_{k=1}^{\infty} \left( 1- \omega^2 z^{2k} \right) = \prod_{k=1}^{\infty} \left( 1- \omega z^{k} \right) 
\left( 1 + \omega z^{k} \right) \;\;.
\label{six-twentyeighta1}\end{eqnarray}
Introducing the rhs of Eq.\ (\ref{six-twentyfive}) into the rhs of the above equation yields
\begin{eqnarray}
\prod_{k=1}^{\infty} \left( 1- \omega^2 z^{2k} \right) = \sum_{k=1}^{\infty} z^k \, \sum_{j=0}^k
q(j,\omega) \,q(k-j,-\omega) \;\;.
\label{six-twentyeighta2}\end{eqnarray}
The equals sign is only valid in Eqs.\ (\ref{six-twentyeighta}) and (\ref{six-twentyeighta2}) for $|\omega| \!<\! 1$ 
and $|z| \!<\! 1$. By equating like powers of $z$ on the rhs's of both these equations, we find that for $k$, an 
odd number equal to $2m \!+\!1$,
\begin{eqnarray}
\sum_{j=0}^{2m+1} q(j,\omega) \,q(2m+1-j,-\omega) = 0 \;\;.
\label{six-twentyeighta3}\end{eqnarray}
On the other hand, for $k \!=\! 2m$, we obtain
\begin{eqnarray}
\sum_{j=0}^{2m} q(j,\omega) \, q(2m-j,-\omega)= q(m,-\omega^2) \;\;.
\label{six-twentyeighta4}\end{eqnarray}
When $\omega \!=\!1$, Eq.\ (\ref{six-twentyeighta4}) reduces to
\begin{eqnarray}
\sum_{j=0}^{2m} L_{DP,j}\left[1\right] \, L_{DP,2m-j} \bigl[ (-1)^{N_{2m-j}} \bigr]= L_{DP,m}\left[(-1)^{N_{m}} \right] \;\;.
\label{six-twentyeighta5}\end{eqnarray}
From Eq.\ (\ref{six-twentyfoura}) we see that the lhs of the above equation is effectively a sum over the pentagonal numbers
less than $2m$, while the rhs is non-zero if and only if $m$ is a pentagonal number. Furthermore, if we let 
$\omega$ equal the complex number $i$, then Eq.\ (\ref{six-twentyeighta4}) gives the number of discrete partitions summing
to $m$ or the number of discrete partitions summing to $2m$ with even elements.

The foregoing analysis can also be extended by raising Eq.\ (\ref{six-twentyfive}) to an arbitrary power $\rho$ 
and applying Corollary\ 1 to Theorem\ 1. Then we obtain a power series expansion like Equivalence\
(\ref{thirtyfour}), but in powers of $z$ with the coefficients depending upon $\rho$. Hence, we arrive at
\begin{eqnarray}
\prod_{k=1}^{\infty} \bigl(1-\omega z^k \bigr)^{\rho} \equiv 1+ \sum_{k=1}^{\infty}
q(k,-\omega,\rho) \, z^k \;\;,
\label{six-twentyeightb}\end{eqnarray}
where the coefficients $q(k,\omega,\rho)$ can determined from Eq.\ (\ref{thirtyfive}) with $D_i \!=\! q(i,\omega)$
and are given by  
\begin{eqnarray}
q(k,\omega,\rho)= L_{P,k} \Bigl[ (-1)^{N_k} \, (-\rho)_{N_k} \prod_{i=1}^k \frac{q(i,\omega)^{n_i}}{n_i!} \Bigr] \;\;.
\label{six-twentyeightc}\end{eqnarray}
It should be emphasised that in the above results $\rho$ can be any value including a complex number. For integer values of 
$\rho$ greater than zero the equivalence symbol can be replaced by an equals sign. As we shall see shortly, when $\rho \!=\! -1$ 
and $\omega \!=\!1$ in Equivalence\ (\ref{six-twentyeightb}), the $q(k,-1,-1)$ equal the partition function or $p(k)$, while 
if $\rho$ is equal to a positive integer, say $j$, then Eq.\ (\ref{six-twentyeightc}) simplifies drastically due to the fact 
that for $k \!>\!j$, the factor $(-\rho)_{N_k}$ vanishes for $N_k \!>\!j$. That is, the partitions with more than $j$ 
elements in them do not contribute to the $q(k,\omega,\rho)$. Moreover, if $\omega \!=\! -1$, further redundancy occurs
in Eq.\ (\ref{six-twentyeightc}) since the $q(i,\omega)$ become the discrete partition numbers or $q(i)$, which we have 
seen are only non-zero when $i$ is a pentagonal number \cite{wei011a}. 

There is also another approach to developing a power series expansion to the generating function given on the lhs of
Equivalence\ (\ref{six-twentyeightb}). This too involves the partition operator, but rather than consider the 
generating function in Equivalence\ (\ref{six-twentyeightb}), we shall investigate the more general case of Eq.\ 
(\ref{six-nineteen}) raised to an arbitrary power $\rho$ as set out in the following theorem.

{\bf Corollary to Theorem\ 3}. The generalised version of the product in Equivalence\ (\ref{six-twentyeightb})
whereby the coefficients of $z^k$ are set equal to $C_k$ can be expressed as a power series or generating form given by 
\begin{eqnarray}
\prod_{k=1}^{\infty} \bigl(1+C_k z^k \bigr)^{\rho_k} \equiv 1+ \sum_{k=1}^{\infty} B_k({\boldsymbol \rho}) \,z^k \;\;,
\label{six-twentyeightc1}\end{eqnarray}  
while the coefficients in the series are found to be 
\begin{eqnarray}
B_k( \boldsymbol{\rho}) = L_{P,k} \Bigl[ (-1)^{N_k} \prod_{i=1}^{k} \frac{(-\rho_i)_{n_i}}{n_i!} \; C_i^{n_i}\Bigr] \;\;.
\label{six-twentyeightc2}\end{eqnarray}
In the above results $(\boldsymbol{\rho})$ denotes $(\rho_1,\rho_2,\dots,\rho_k)$. If $S \!=\! {\rm inf}\; |C_i|^{-1/i} 
\!>\! 0$, then for $|z| \!<\! S$ the equivalence symbol can be replaced by an equals sign. In addition, the coefficients 
$B_k({\boldsymbol \rho})$ satisfy the following relations:
\begin{eqnarray}
L_{P,k} \Bigl[ (-1)^{N_k} N_k! \prod_{i=1}^{k} \frac{B_i(\boldsymbol{\rho})^{n_i}}{n_i!}\Bigr] =
L_{P,k} \Bigl[ (-1)^{N_k} \prod_{i=1}^{k} \frac{(\rho_i)_{n_i}}{n_i!} \; C_i^{n_i}\Bigr] \;\;,
\label{six-twentyeightc2a}\end{eqnarray}
and 
\begin{eqnarray}
B_k(\boldsymbol{\mu}+\boldsymbol{\nu}) = \sum_{j=0}^k B_j(\boldsymbol{\mu}) \, B_{k-j}(\boldsymbol{\nu}) \;\;.
\label{six-twentyeightc2b}\end{eqnarray}

{\bf Remark}. The reader should observe that in previous cases involving a constant power of $\rho$ the Pochhammer 
symbol appeared outside the product in the partition operator as in Corollary\ 1 to Theorem\ 1. Thus, the 
$\rho$-dependence of the coefficients in the resulting power series was only affected by the total number of 
elements in each partition. In the above case each element $i$ is assigned a value that is dependent upon $\rho_i$
and consequently, $(\rho_i)_{n_i}$ appears inside the product being acted upon by the partition operator. 
 
{\bf Proof}. In order to prove the theorem, we shall use Lemma\ 1 again. Then the generating function can be
expressed as
\begin{align}
\prod_{k=1}^{\infty}  & \bigl(1+C_k z^k \bigr)^{\rho_k}  \equiv \Bigl(1 + \sum_{j=1}^{\infty} \frac{(-\rho_1)_j}{j!}\,
(-C_1 \, z)^j \Bigr) \Bigl(1 + \sum_{j=1}^{\infty} \frac{(-\rho_2)_j}{j!}\, \left(-C_2 \, z^2 \right)^j \Bigr) 
\nonumber\\
& \times \;\; 
\Bigl(1 + \sum_{j=1}^{\infty} \frac{(-\rho_3)_j}{j!}\, \left( -C_3 \, z^3 \right)^j \Bigr) \Bigl(1 + \sum_{j=1}^{\infty} 
\frac{(-\rho_4)_j}{j!}\, \left( -C_4 \,  z^4 \right)^j \Bigr)  \cdots \;\;.
\label{six-twentyeightc3}\end{align}  
According to Lemma\ 1, each binomial series in the above product is absolutely convergent in the disk given by $|z| \!<\! 
|C_i|^{-1/i}$. Each series is also conditionally convergent in a specific region of the complex plane. From the proof of 
Theorem\ 2, which deals with the $\rho_i \!=\! -1$ and $C_i \!=\!-1$ case, we found that these regions can cancel each other 
when all the series appear in a product as in the above generating function. As a result, we were left with the region in which 
all the series are only absolutely convergent. Therefore, if $S= {\rm inf}\, |C_i|^{-1/i} \!> \! 0$, then all the series in 
Equivalence\ (\ref{six-twentyeightc3}) are absolutely convergent whenever $|z| \!<\! S$ and the equivalence symbol can be 
replaced by an equals sign. Furthermore, when the $\rho_i$ all equal a positive integer, say $k$, the binomial series 
become polynomials of degree $k$. Then the equivalence symbol can be replaced by an equals sign for any value of $z$.

Expanding the product of all the series in Equivalence\ (\ref{six-twentyeightc3}) in powers of $z$ yields
\begin{align}
\prod_{k=1}^{\infty}  & \bigl(1+ C_k z^k \bigr)^{\rho_k}  \equiv  1 - (-\rho_1)_1\,C_1 \,z 
+ \left( \frac{(-\rho_1)_2}{2!} \, C^2_1 - (-\rho_2)_1\, C_2 \right) z^2
\nonumber\\
&+ \;\; \left( -\frac{(-\rho_1)_3}{3!}\, C_1^3+ (-\rho_1)_1 \, (-\rho_2)_1\, C_1\, C_2 - (-\rho_3)_1 \,C_3 \right) z^3 +  
\cdots \;\;.
\label{six-twentyeightc4}\end{align}  
From this result we see that there is one term appearing in the first order coefficient in $z$, while there are two terms
appear in the second order coefficient. The third order coefficient is composed of three terms. Had the fourth order
term been displayed, there would have been five terms in the coefficient. In fact, the number of terms in the $k$-th order 
coefficient is the number of partitions summing to $k$ or $p(k)$. Therefore, we need to develop a means of coding the partitions 
as was done for the partition method for a power series expansion in Theorem\ 1.

Since the first order term corresponds to \{1\}, we assign a value of $-C_1$ to each occurrence of a one in a
partition. The first term in the second order coefficient possesses a factor of $C_1^2$. Therefore, it must correspond
to \{1,1\}, while the other term must correspond to the other partition summing to 2, viz.\ \{2\}.  This term possesses 
a factor of $-C_2$. Hence, each occurrence of a two in a partition yields a factor of $-C_2$. If we continue this process 
indefinitely for the higher orders, then we find that each occurrence of an element $i$ yields a factor $-C_i$. 

This, however, is not the complete story. Accompanying the factor of $-C_1$ in the first order term is the factor of 
$(-\rho_1)_1$, while the first and second terms in the second order coefficient possess factors of $(-\rho_1)_2/2!$ and 
$(-\rho_2)_1$, respectively. That is, when there is one occurrence of an element $i$ in a partition, its assigned value must be 
multiplied by $(-\rho_i)_1$, but if there are two occurrences of an element in a partition, then the assigned value must be 
multiplied by $(-\rho_i)_2/2!$. Therefore, if an element $i$ occurs $n_i$ times in a partition, then it contributes a 
factor of $(-\rho_i)_{n_i} (-C_i z^i)/n_i!$ to the coefficient. This can be checked with the various terms comprising
the third order term. The total contribution made by a partition is then given by the product over all elements, viz.\ 
$\prod_{i=1}^{k} (-\rho_i)_{n_i}(-C_i)^{n_i}/n_i!$. Finally, the coefficient $B_k(\boldsymbol{\rho})$ is evaluated by summing 
over all partitions summing to $k$. Hence, we arrive at Eq.\ (\ref{six-twentyeightc2}).      

If we assume that $z \!<\! S$, which is actually not necessary, we can invert the equation form of Equivalence\ 
(\ref{six-twentyeightc1}), thereby obtaining
\begin{eqnarray}
\prod_{k=1}^{\infty} \bigl(1+C_k z^k \bigr)^{-\rho_k} = \frac{1}{1+ \sum_{k=1}^{\infty} B_k(\boldsymbol{\rho}) \,z^k } \;\;.
\label{six-twentyeightc5}\end{eqnarray}  
Because the rhs can be regarded as the regularised value of the geometric series, we can apply the method for a power 
series expansion where the coefficients of the inner series, viz.\ $p_k$ in Theorem\ 1, are equal to $-B_k(\boldsymbol{\rho})$,
while the coefficients of the outer series $q_k$ are once again equal to $(-1)^k$. Then we have 
\begin{eqnarray}
\prod_{k=1}^{\infty} \bigl(1+ C_k z^k \bigr)^{-\rho_k} \equiv 1+ \sum_{k=1}^{\infty} D_k \,z^k  \;\;,
\label{six-twentyeightc6}\end{eqnarray}  
where according to Eq.\ (\ref{five}), the coefficients of this expansion are given by
\begin{eqnarray}
D_k = L_{P,k} \Bigl[ (-1)^{N_k} \, N_{k}! \prod_{i=1}^k \frac{B_i(\boldsymbol{\rho})^{n_i}}{n_i!}  \Bigr] \;\;.
\label{six-twentyeightc7}\end{eqnarray}  
We also know from earlier in the proof that the above product can be expressed in terms of a power series of the form 
given on the rhs of Equivalence\ (\ref{six-twentyeightc1}) except that now the $\rho_k$ are replaced by $-\rho_k$. 
Because the resulting power series expansions possess the same regularised value, they are equal to one another. Moreover,
since $z$ is arbitrary, we can equate like powers of $z$, thereby yielding Eq.\ (\ref{six-twentyeightc2a}). 

The final identity is easily proved by noting that 
\begin{eqnarray} 
\prod_{i=1}^{\infty} \left( 1+ C_k z^k \right)^{\mu_k+\nu_k} = \prod_{i=1}^{\infty} \left( 1+ C_k z^k \right)^{\mu_k} 
\left( 1+ C_k z^k \right)^{\nu_k} \;\;. 
\label{six-twentyeightc8}\end{eqnarray}
The lhs represents the regularised value for the series on the rhs of Equivalence\ (\ref{six-twentyeightc1}) with
coefficients $B_k(\boldsymbol{\mu}+\boldsymbol{\nu})$, while the rhs represents the regularised value of the product of 
two series, one with coefficients $B_k(\boldsymbol{\mu})$ and the other with coefficients $B_k(\boldsymbol{\nu})$. As the 
regularised value is the same in both situations, we have 
\begin{align} 
1+ \sum_{k=1}^{\infty} B_k(\boldsymbol{\mu}+\boldsymbol{\nu}) \,z^k & = \Bigl( 1+ \sum_{j=1}^{\infty} B_k(\boldsymbol{\mu}) 
\,z^k \Bigr) \Bigl( 1 + \sum_{k=1}^{\infty} B_k(\boldsymbol{\nu}) \,z^k \Bigr) 
\nonumber\\
& = \;\; 1+ \sum_{k=1}^{\infty} z^k \sum_{j=0}^k B_j(\boldsymbol{\mu}) B_{k-j}(\boldsymbol{\nu}) \;\;.
\label{six-twentyeightc9}\end{align}
Since $z$ is arbitrary, we can equate like powers yet again. Therefore, we obtain Eq.\ (\ref{six-twentyeightc2b}), which
completes the proof. 

In order to make the foregoing material clearer, we now consider a few examples. We have already generalised the
discrete partition polynomials $q(k,\omega)$ by introducing the power of $\rho$ into their associated product as 
demonstrated by Equivalence\ (\ref{six-twentyeightb}).  It was found that the new polynomials $q(k,\omega,\rho)$ 
could be expressed in terms of the partition operator acting on the discrete partition polynomials given by
Eq.\ (\ref{six-twentyeightc}). Now we apply the corollary to Theorem\ 3 to the product with the $C_k$ and $\rho_k$ set 
equal to $\omega$ and $\rho$ respectively. Consequently, we find that
\begin{eqnarray}
q(k,\omega,\rho)= L_{P,k} \Bigl[ (- \omega)^{N_k} \prod_{i=1}^{k} \frac{(-\rho)_{n_i}}{n_i!} \Bigr] \;\;.
\label{six-twentyeightc10}\end{eqnarray}
So, let us evaluate $q(3,\omega,\rho)$ via Eqs.\ (\ref{six-twentyeightc}) and (\ref{six-twentyeightc10}). In the case of
Eq.\ (\ref{six-twentyeightc}) we require $q(1,\omega)$, $q(2,\omega)$ and $q(3,\omega)$, which have already been evaluated. 
Therefore, we obtain
\begin{eqnarray}
q(3,\omega,\rho)= -\frac{(-\rho)_3}{3!} \, \omega^3 + (-\rho)_2 \, \omega^2- (-\rho)_1 \, (\omega^2 +\omega) \;\;.
\label{six-twentyeightc11}\end{eqnarray}
After a little algebra, the above becomes
\begin{eqnarray}
q(3,\omega,\rho)= \Bigl( \frac{\rho^3}{6} -\frac{\rho^2}{2} + \frac{\rho}{3} \Bigr) \, \omega^3 +\rho^2 \omega^2+
\rho \omega  \;\;.
\label{six-twentyeightc12}\end{eqnarray}
By setting $k \!=\!3 $ in Eq.\ (\ref{six-twentyeightc10}), we obtain
\begin{eqnarray}
q(3,\omega,\rho)= -\frac{(-\rho)_3}{3!} \,\omega^3 +(-\rho)_1 \, (-\rho)_1\,\omega^2 - (-\rho)_1 \, \omega  \;\;.
\label{six-twentyeightc13}\end{eqnarray}
Again after a little algebra, we end up with Eq.\ (\ref{six-twentyeightc12}), although we can see that Eqs.\
(\ref{six-twentyeightc11}) and (\ref{six-twentyeightc13}) are composed of different quantities. This demonstrates
that while Eqs.\ (\ref{six-twentyeightc}) and (\ref{six-twentyeightc10}) both possess the same sums over partitions,
they are in fact different representations for $q(k,\omega,\rho)$. 

The next example, which is a result attributed to Euler, appears on p.\ 56 of Ref.\ \cite{knu005}. This is 
\begin{eqnarray}
\prod_{k=1}^{\infty} \bigl( 1 - z^k \bigr)^{3} = 1 -3z+ 5z^3-7z^6 + \cdots= \sum_{k=0}^{\infty}
(-1)^k (2k+1) z^{\binom{k+1}{2}} \;\;.
\label{six-twentyeightd}\end{eqnarray}
From Equivalence\ (\ref{six-twentyeightb}), we have
\begin{eqnarray}
\prod_{k=1}^{\infty} \bigl(1- z^k \bigr)^{3} = 1+ \sum_{k=1}^{\infty}
q(k,-1,3) \, z^k \;\;.
\label{six-twentyeightd1}\end{eqnarray}
Note that the equivalence symbol has been replaced by an equals sign because $\rho \!=\! 3$ in this case. Euler's result
shows that only when the power of $z$ is another type of figurate number called a triangular number \cite{wei011b,wik11b}, are
the coefficients of the power series or generating function non-zero. If we equate like powers of $z$ in both power series 
expansions given above, then we arrive at
\begin{eqnarray}
L_{P,k} \Bigl[ (-1)^{N_k} \, (-3)_{N_k} \prod_{i=1}^k \frac{q(i)^{n_i}}{n_i!} \Bigr] =\begin{cases}
(-1)^j (2j+1), & {\rm if} \;\; k= \binom{j+1}{2}, \\
0, & {\rm otherwise}. \end{cases}
\label{six-twentyeighte}\end{eqnarray}
On the other hand, putting $C_i \!=\! -1$ and $\rho_i \!=\! 3$ in Eq.\ (\ref{six-twentyeightc2}) yields
\begin{eqnarray}
q(k,-1,3) = L_{P,k} \Bigl[ \prod_{i=1}^k \frac{(-3)_{n_i}}{n_i!} \Bigr] \;\;.
\label{six-twentyeighte1}\end{eqnarray}
In the above result $(-3)_{n_i}$ is only non-zero for $n_i$ equal to 1, 2 and 3, in which case it equals -3, 6 and -6, respectively.
This is a different type of restricted partition from those we have encountered previously since it means that partitions 
in which an element appears more than three times are excluded. To generate such partitions all that needs to be done is 
to scan each partition twice by introducing another for loop in ${\bf termgen}$ of the partition generating program
presented towards the end of Sec.\ 2. If in the first scan the $n_i$ or the variable ${\it freq}$ is greater than 3, 
then a goto statement is required so that the program avoids the next for loop, which is responsible for printing out 
the specific partitions. When such a code is constructed, one will find that out of the total of 627 partitions summing 
to 20, there are only 320 partitions in which all the elements occur at most three times. As a result of Eq.\ 
(\ref{six-twentyeighte1}), we also arrive at 
\begin{eqnarray}
L_{P,k} \Bigl[ \prod_{i=1}^k \frac{(-3)_{n_i}}{n_i!} \Bigr] =\begin{cases}
(-1)^j (2j+1), & {\rm if} \;\; k= \binom{j+1}{2}, \\
0, & {\rm otherwise}. \end{cases}
\label{six-twentyeighte2}\end{eqnarray}
Moreover, it should be noted that if $\rho \!=\! 1$ such that $n_i \!\leq 1$ for all elements $i$ in the partition, then the
sum over all partitions in Eq.\ (\ref{six-twentyeightc2}) reduces to the discrete partition operator, viz.\ $L_{DP,k}[ \cdot]$,
irrespective of the values for the coefficients $C_i$. 

If we put $\rho \!=\! 1$, take the cube power of the series on the rhs of Eq.\ (\ref{six-twentyeightb}) and equate like 
powers of the resulting series with those on the rhs of Eq.\ (\ref{six-twentyeightd1}), then  we obtain 
\begin{eqnarray}
q(k,-1,3) = \sum_{j_1=0}^{k} \sum_{j_2=0}^{j_1}q(k-j_1) q(j_1-j_2) q(j_2) \;\;,
\label{six-twentyeightf}\end{eqnarray}
while multiplying the $\rho \!=\!1$ and $\rho \!=\!2$ versions of the series on the rhs of Eq.\ (\ref{six-twentyeightb})
yields
\begin{eqnarray}
q(k,-1,3) = \sum_{j=0}^{k} q(j) q(k-j,-1,2) = \begin{cases} (-1)^i (2i+1) , & \;\; k= \binom{i+1}{2}, \\
0, & \;\; {\rm otherwise}. \end{cases}
\label{six-twentyeightg}\end{eqnarray}
In the above result $q(2,\omega,2)$ can be evaluated by putting $\rho \!=\! 2$ in Eq.\ (\ref{six-twentyeightc}). They can also be
determined by equating like powers of $z$ when taking the square of the series on the rhs of Equivalence\ (\ref{six-twentyeightb}) 
with the $\rho \!=\! 2$ series. This gives
\begin{eqnarray}
q(k,\omega,2) =\sum_{i=0}^k q(i,\omega) q(k-i,\omega) \;\; . 
\label{six-twentyeightg1}\end{eqnarray}
For $\omega \!=\! -1$, Eq.\ (\ref{six-twentyeightg1}) reduces to
\begin{eqnarray}
q(k,-1,2) =\sum_{i=0}^k q(i) q(k-i) \;\; . 
\label{six-twentyeightg2}\end{eqnarray}
Table\ \ref{table2b} displays both $q(k,\omega,2)$ and $q(k,\omega,3)$ up to $k \!=\! 10$. 
 
\begin{table}
\begin{center}
\begin{tabular}{|c|l|l|} \hline
$k$ &  $q(k,\omega,2)$ & $q(k,\omega,3)$ \\ \hline
$0$ & $ 1$ & $1$ \\  
$1$ &  $ 2 \omega $ & $3 \omega$  \\ 
$2$ & $ 2 \omega + \omega^2$ & $ 3 \omega + 3 \omega^2$ \\ 
$3$ & $2\omega + 4 \omega^2  $ & $3\omega + 9\omega^2 +\omega^3$ \\ 
$4$ & $2 \omega + 5 \omega^2 +2 \omega^3 $ & $3 \omega +12 \omega^2+9 \omega^3$ \\ 
$5$ & $2 \omega + 8 \omega^2 +4 \omega^3 $ & $3 \omega +18 \omega^2+18 \omega^3+ 3\omega^4$ \\ 
$6$ & $2 \omega + 9 \omega^2 + 10 \omega^3 +\omega^4 $ & $3 \omega +21 \omega^2+37 \omega^3+ 12 \omega^4$ \\ 
$7$ & $2 \omega +12  \omega^2 +14 \omega^3 +4 \omega^4 $ & $3 \omega+27 \omega^2+54 \omega^3+33 \omega^4+3 \omega^5$ \\ 
$8$ & $2 \omega +13 \omega^2 + 22  \omega^3 +9 \omega^4 $ & $-3 \omega+30 \omega^2+81 \omega^3+66 \omega^4+12 \omega^5$\\
$9 $ & $2 \omega+ 16 \omega^2 +30\omega^3 +16  \omega^4$ & $3\omega +36 \omega^2 +109 \omega^3 
+114 \omega^4+39 \omega^5$ \\
$$ & $2\omega^5$ & $ +\omega^6$ \\ 
$10$ & $2 \omega+ 17 \omega^2+40 \omega^3+30 \omega^4$ & $ 3 \omega+39 \omega^2 +144 \omega^3+
189 \omega^4+81 \omega^5$ \\ 
$$ & $4 \omega^5$ & $+9 \omega^6$ \\ \hline
\end{tabular}
\end{center}
\caption{Coefficients $q(k,\omega,2)$ and $q(k,\omega,3)$ in the power series expansions of the $\rho \! =\! 2$ and
$\rho \!=\! 3$ cases of Equivalence\ (\ref{six-twentyeightb}).}
\label{table2b}
\end{table}

On p.\ 23 of Ref.\ \cite{and03} it is stated that Gauss derived the following result:
\begin{eqnarray}
\sum_{k=0}^{\infty} z^{(k^2+k)/2} = \prod_{k=1}^{\infty} \frac{\bigl( 1-z^{2k} \bigr)}{\bigl(1-z^{2k-1}\bigr)}
=\prod_{k=1}^{\infty}  \bigl(1-z^k \bigr) \bigl(1+z^k \bigr)^2 \;\;.
\label{six-twentyeighth}\end{eqnarray}
Once again, the equals sign is only valid for $|z| \!<\! 1$. From this result we see that only the powers of $z$ equal 
to $(j^2+j)/2$, where $j$ is any non-negative integer, possess a non-zero coefficient. For any power $k$ we obtain 
contributions from all the partitions summing to $k$ with only odd elements in them and from the discrete partitions 
summing to $k$ consisting of only even elements. If the number of even elements in these discrete partitions is 
even, then the partition will yield a value of unity. Otherwise, it will yield a value of -1. In addition, the partition 
function possesses mixed partitions composed of discrete partitions with only even elements and standard partitions 
with only odd elements. The values contributed by the mixed partitions to the coefficients of the power series depend 
upon the number of elements in the discrete partitions. For example, if we consider $z^{10}$ or $k \!=\!4$ on the lhs 
of Eq.\ (\ref{six-twentyeighth}), then it will be composed of the contributions due to the partitions summing to 10 with 
only odd elements. There are 10 of these beginning with \{1,9\} and ending with \{5,5\}. Hence, these partitions contribute 
a value of 10 to the coefficient. On the other hand, there are only 3 discrete partitions with even numbers summing to 10. 
These are \{10\}, \{2,8\} and \{4,6\}. Since the last two possess an even number of elements, they each contribute a value 
of unity, while the single element partition gives a value of -1. Overall, the discrete partitions summing to 10 contribute 
a value of unity to the the coefficient, which now becomes 11 when the contribution from the standard partitions summing to
10 with only odd elements is included. However, when the discrete partition is \{2\}, we need to consider the partitions 
summing to 8 with odd elements. There are six of these, beginning with \{1,7\} and ending with \{3,5\}. Because there is 
only one element in the discrete partition, these mixed partitions contribute a value of -6 to the coefficient, which now
drops to 5. However, there are still more mixed partitions. We need to consider the discrete partition of \{4\}. In this 
case we need the partitions summing to 6 with odd elements. There are four of these, beginning with \{1,5\} and ending with 
\{3,3\}. In this instance the mixed partitions contribute a value of -4, yielding a value of 1 for the coefficient as 
indicated above. We do not need to consider the contributions where the discrete partitions sum to either 6 or 8 because 
in these cases there are only 2 partitions, one of which has two elements and one with one element. Hence, they cancel each 
other yielding a value of 0.   

If we introduce Eq.\ (\ref{six-twentyeightb}) into the Eq.\ (\ref{six-twentyeighth}), then we find that
\begin{eqnarray}
\sum_{k=0}^{\infty} z^{k^2+k/2} = 1+ \sum_{k=0}^{\infty} z^k \sum_{j=0}^{k} q(j) q(k-j,1,2) \;\;.
\label{six-twentyeighti}\end{eqnarray}
By equating like powers of $z$, we obtain for $i$, a positive integer,
\begin{eqnarray}
\sum_{j=0}^{k} q(j) q(k-j,1,2) =  \begin{cases} 1, & \;\; k=(i^2+i)/2, \\
0, & \;\;  {\rm otherwise} . \end{cases}
\label{six-twentyeightj}\end{eqnarray}
Therefore, when $k$ is not a triangular number \cite{wei011b}, by combining the above result with Eq.\ (\ref{six-twentyeightg}) 
we arrive at
\begin{eqnarray}
\sum_{j=0}^{k} q(j) \Bigl( q(k-j,1,2) \pm q(k-j,-1,2) \Bigr)=0 \;\;.
\label{six-twentyeightk}\end{eqnarray}
 
\section{Other Products}
By using the material of the previous section we are in a position to study more advanced products. We begin by introducing
the variable $\omega$ next to the power of $z$ in the denominator of $P(z)$. Based on the similar extension of the product 
yielding discrete partitions in the previous section, we expect to obtain polynomials as the coefficients of the 
resulting generating function. Therefore, we define the new product as
\begin{eqnarray}
P(z,\omega) = \prod_{k=1}^{\infty} \frac{1}{(1-\omega z^k)} \;\;.
\label{six-twentynine}\end{eqnarray}
Obviously, when $\omega \!=\! 1$, this reverts to the generating function of the partition function given
by Equivalence\ (\ref{fortyfoura}), which we have seen becomes an equation when $|z| \!<\! 1$. According
to p.\ 112 of Ref.\ \cite{knu005}, the product in the above result can be written alternatively as
\begin{eqnarray}
P(z,\omega) = (1- \omega) \sum_{k=0}^{\infty} \frac{\omega^k \; z^{k^2}}{(z;z)_{k} \; (\omega;z)_{k+1}} \;\;,
\label{six-thirty}\end{eqnarray} 
while inversion of the rhs of Eq.\ (\ref{six-twentyfive}) yields
\begin{eqnarray}
P(z,\omega) = \frac{1}{1+ \sum_{k=1}^{\infty} q(k,-\omega) \,z^k}  \;\;.
\label{six-thirtyone}\end{eqnarray}

We have seen that $P(z)$ or $\omega \!=\! 1$ in Eq.\ (\ref{six-twentynine}) yields a power series expansion whose
coefficients are given by the partition function or $p(k)$. This expansion is obtained by expanding each term in the 
generating function into the geometric series for each value of $i$. It is this value in the generating function,
which is responsible for yielding the specific elements in a partition, while the power of $z^i$ in each geometric series 
represents the frequency or number of occurrences of the element in the partition. For example, multiplying $(z^2)^3$
in the expansion of $1/(1 \!-\! z^2)$ by $(z^3)^4$ in the expansion of $1/(1 \!-\! z^3)$ means that the partition 
has three twos and four threes in it. By introducing $\omega$ into the generating function as indicated above, we see 
that the overall power of $\omega$ yields the total number of elements in a partition. In the example just mentioned 
we now obtain $(\omega z^2)^3$ multiplied by $(\omega z^3)^4$, which yields $\omega^7 z^{18}$. The power of 7 on 
$\omega$ represents the total number of twos and threes in the partition. Therefore, we expect the coefficient of
each power of $\omega$ in the coefficients of the resulting generating function to indicate the total number of partitions
where the number of elements equals the power of $\omega$. 

By regarding the rhs as the regularised value of the geometric series and the $q(k,-\omega)$ as the inner 
series whose coefficients are $p_k$, we can apply Theorem\ 1 to the last form for $P(z,\omega)$ with
$y \!=\!z$.  If we express the generating function as a power series expansion with coefficients, $p(k,\omega)$,
i.e.
\begin{eqnarray}
P(z,\omega) = \prod_{k=1}^{\infty} \frac{1}{(1-\omega z^k)} \equiv 1+\sum_{k=1}^{\infty} p(k,\omega)
\, z^k \;\;,
\label{six-thirtytwo}\end{eqnarray}
then according to Eq.\ (\ref{five}) the coefficients of the resulting power series expansion are given by
\begin{eqnarray}
p(k,\omega) = L_{P,k} \left[ (-1)^{N_k} \, N_k! \prod_{i=1}^k \frac{q(i,-\omega)^{n_i}}{n_i!} \right] \;\;. 
\label{six-thirtytwoa}\end{eqnarray}
Furthermore, if we put $\rho \!=\! -1$ and $C_i \!=\! -\omega$ in the corollary to Theorem\ 3, then we observe that 
the $p(k,\omega)$ become the coefficients $B_k(-1)$ given by Eq.\ (\ref{six-twentyeightc2}). Thus, we arrive at
\begin{eqnarray}
p(k,\omega)= L_{P,k} \Bigl[ \omega^{N_k} \Bigr] \;\;.
\label{six-thirtytwob}\end{eqnarray}
This tells us that the coefficients in the $p(k,\omega)$ will be the number of partitions summing to $k$ where
each power of $\omega$ corresponds to the number of elements in the partitions. That is,
\begin{eqnarray}
p(k,\omega) = \sum_{i=1}^k p_i(k) \, \omega^i \;\;,
\label{six-thirtythree}\end{eqnarray}
where from Sec.\ 2, $p_i(k) \!=\! \begin{vmatrix} k \\ i\end{vmatrix}$. It has already been stated the sub-partition 
numbers obey the recurrence relation given by Eq.\ (\ref{one}). In terms of the fixed number of elements operator defined 
by Eq.\ (\ref{fortyfourl}) we also find that
\begin{eqnarray}
p_i(k) = L_{P,k}^i \Bigl[ 1 \Bigr]  \;\;.
\label{six-thirtythreea}\end{eqnarray}

The first few partition function polynomials are found to be: $p(0,\omega) \!=\! 1$, $p(1,\omega) \!=\! 
\omega$, $p(2,\omega) \!=\! \omega^2 \!+\! \omega$, $p(3,\omega) \!=\! \omega^3 \!+\! \omega^2 \!+\! \omega$ and 
$p(4,\omega) \!=\! \omega^4 \!+\! \omega^3 \!+\! 2 \omega^2 \!+\! \omega$. In fact, those up to $k\!=\! 10$ are
displayed in the third column of Table\ \ref{table2a}. Since the coefficients of these polynomials represent 
the number of partitions in which the number of elements is given by the power of $\omega$, the highest order 
term of these polynomials is $k$, which arises from the $k$-element partition of \{1,1,\dots,$1_k$\}. The other 
partitions are unable to provide an $\omega^k$ term because the highest order of all other $q(k,\omega)$ is less 
than $k$. The partition \{1,1,\dots,$1_{k-1}$,2\} only produces an $\omega^{k-1}$ term as its highest order term 
because the highest order term of $q(2,\omega)$ is 1. Hence, deg $p(k,\omega) \!=\! k$. Conversely, the lowest order 
term in the $p(k,\omega)$ is the lowest order term in $q(k,\omega)$ stemming from the single element partition 
\{$k$\}, which is unity. Therefore, $p_k(k) \!=\! p_{k-1}(k) \!=\! p_1(k) \!=\! 1$, $p_{k-2}(k) \!=\! 2$ and 
$p(k,1) \!=\! \sum_{i=1}^{k} p_i(k) \!=\! p(k)$. Moreover, the total number of partitions with even elements is 
given by $\sum_{i=1}^{[k/2]} p_{2i}(k)$, while the total number of partitions with odd elements is equal to 
$\sum_{i=1}^{m} p_{2i-1}(k)$, where $m \!=\! k/2$ when $k$ is even and $m \!=\! [k/2] +1$ when $k$ is odd. 
On the other hand, the coefficient $p_2(k)$ can be evaluated by noting that it represents the product of the two 
lowest order terms in each partition, namely \{$j,k-j$\}, where $j$ ranges from 1 to $[k/2]$. According to Ref.\ 
\cite{wei011c}, the number of two-element partitions summing to $k$ is given by $p_2(k) \!=\! [k/2]$, while the 
number of three-element partitions is given by $p_3(k) \!=\! [k^2/12]$ for $k \!>\!3$. Moreover, a table of 
the sub-partition numbers for $k$ and $i$ ranging from 0 to 11 is presented on p.\ 46 of Ref.\ \cite{knu005}. All 
these results agree with those obtained via Eq.\ (\ref{six-twentyone}), confirming that the latter result does yield 
the number of partitions summing to $k$ with $i$ elements in them.

An interesting property of the partition number polynomials can be derived by setting $z \!=\! z^2$ and $\omega
\!=\! \omega^2$ in Equivalence\ (\ref{six-thirty}). Then we obtain
\begin{eqnarray}
P \left(z^2,\omega^2 \right) = \prod_{k=1}^{\infty} \frac{1}{(1-\omega^2 z^{2k})} \equiv 1+\sum_{k=1}^{\infty} 
p \left(k,\omega^2 \right) \, z^{2k} \;\;.
\label{six-thirtythreea1}\end{eqnarray}
The quantity on the lhs of the above equivalence can also be written 
\begin{eqnarray}
P \left(z^2,\omega^2 \right) = P(z,\omega) P(z,-\omega) \;\;.
\label{six-thirtythreea2}\end{eqnarray}
Introducing Equivalence\ (\ref{six-thirtytwo}) into the above equation yields
\begin{eqnarray}
P(z,\omega) P(z,-\omega) \equiv \Bigl( 1+\sum_{k=1}^{\infty} p(k,\omega) \, z^k \Bigr) 
\Bigl( 1+\sum_{k=1}^{\infty} p(k,-\omega) \, z^k \Bigr) \;\;.
\label{six-thirtythreea3}\end{eqnarray}
Since the series on the rhs's of Equivalences\ (\ref{six-thirtythreea1}) and (\ref{six-thirtythreea3}) are
derived from the identity given by Eq.\ (\ref{six-thirtythreea2}), they are equal to one another in accordance
with the concept of regularisation \cite{kow11,kow09a,kow001,kow002}. Because $z$ is arbitrary, once 
again we can equate like powers of $z$. Therefore, we obtain
\begin{eqnarray}
\sum_{j=0}^{2k+1} p(j,\omega) \, p(2k+1-j,-\omega) = 0 \;\;,
\label{six-thirtythreea4}\end{eqnarray}
and
\begin{eqnarray}
\sum_{j=0}^{2k} p(j,\omega) \, p(2k-j,-\omega) = p \bigl( 2k,\omega^2 \bigr) \;\;.
\label{six-thirtythreea5}\end{eqnarray}

As in the previous section we can generalise the foregoing analysis by replacing $-\omega$ in $P(z,\omega)$ 
by $C_i$, which effectively represents the inversion of Eq.\ (\ref{six-nineteen}). From Theorem\ 2 we know 
that
\begin{align}
H(z)^{-1} &=\prod_{i=1}^{\infty} \frac{1} {\bigl( 1+C_i\,z^i \bigr)} \equiv \bigl( 1 -C_1 z+ C_1^2 z^2 + \cdots \bigr) 
\bigl(1 -C_2 z^2+C_2^2 z^4+ \cdots \bigr)
\nonumber\\
& \times \;\; \bigl( 1-C_3z^3+ C_3^2 z^6+ \cdots \bigr) \bigl( 1-C_4 z^4 + C_4^2 z^8 + \cdots \bigr) \cdots  \;\;,
\label{six-thirtythreeb}\end{align}
where the equivalence symbol can be replaced by an equals sign provided $|C_i z^i| \!<\! 1$ for all $i$. Expanding the 
above yields
\begin{eqnarray}
H(z)^{-1} \equiv 1- C_1 z + \bigl( C_1^2-C_2) z^2+ \bigl(-C_1^3+ C_2 \, C_1-C_3 \bigr) z^3 + \cdots  \;\;. 
\label{six-thirtythreec}\end{eqnarray}
Each coefficient of $z^k$ in the above power series is composed of contributions that can be related to the partitions
summing to $k$ as was the case in the proof of Theorem\ 1. The major difference between the above situation and that in
the proof of Theorem\ 1 is that there is no multinomial factor associated with each contribution made by a partition as
we also found in the proof to the corollary to Theorem\ 3. Furthermore, each element $i$ in a partition is now assigned a 
value of $-C_i$ so that the above result becomes
\begin{eqnarray}
H(z)^{-1} \equiv \sum_{k=0}^{\infty} H_k \, z^k \;\;,
\label{six-thirtythreed}\end{eqnarray}  
where $H_0 \!=\! 1$, and 
\begin{eqnarray}
H_k = L_{P,k} \left[ (-1)^{N_k} \, \prod_{i=1}^k C_i^{n_i} \right] \;\;.
\label{six-thirtythreee}\end{eqnarray}
The above result represents the case when the $\rho_k$ in the corollary to Theorem\ 3 are set equal to -1.

If the $C_i$ are set equal to unity in Equivalence\ (\ref{six-thirtythreeb}), then we find that
\begin{eqnarray} 
\prod_{i=1}^{\infty} \frac{1}{\bigl(1+z^i \bigr)} \equiv 1+ \sum_{k=1}^{\infty} L_{P,k} \Bigl[ (-1)^{N_k} \Bigr]
\,z^k \;\;,
\label{six-thirtythree1}\end{eqnarray}
where the equivalence symbol can be replaced by an equals sign for $|z| \!<\! 1$. Therefore, the coefficients on
the power series expansion on the rhs represent the difference between the number of even- and odd-element partitions
summing to $k$. In addition, the above equivalence is analogous to putting the $C_i$ equal to unity in Eq.\ 
(\ref{six-nineteen}) and applying Theorem\ 1 to its inverted form. In this case the coefficients of the inner series
are given by $p_k \!=\! L_{DP,k}[1]$, the number of discrete partitions, while the coefficients of the outer series
are given by $q_k=(-1)^k$. Then the coefficients of the power series expansion on the rhs of Equivalence\ 
(\ref{six-thirtythree1}) can be expressed in terms of the number of discrete partitions as
\begin{eqnarray}
L_{P,k} \Bigl[ (-1)^{N_k} \Bigr]= L_{P,k} \Bigl[N_k! \prod_{i=1}^k L_{DP,i} \bigl[ 1 \bigr]^{n_i}/n_i! \Bigr] \;\;.
\label{six-thirtythree1a}\end{eqnarray}

If we multiply the product on the lhs of Equivalence\ (\ref{six-thirtythree1}) by $(1-z^i)$ in both the 
numerator and denominator, then we find that
\begin{eqnarray} 
\prod_{i=1}^{\infty} \frac{1}{\bigl(1+z^i \bigr)} =\prod_{i=0}^{\infty} \bigl( 1- z^{2i+1} \bigr) \;\;.
\label{six-thirtythree2}\end{eqnarray}
When the product on the rhs of Eq.\ (\ref{six-thirtythree2}) appeared in powers of $z^i$ rather
than $z^{2i+1}$, we saw that the resulting power series possessed coefficients which were equal to the number of 
discrete partitions summing to $k$. In the above result all the even powers are now missing. This means that the 
coefficients of the resulting power series will be the number of discrete partitions with only odd elements in them. 
That is,
\begin{eqnarray} 
\prod_{i=0}^{\infty} \bigl( 1- z^{2i+1} \bigr)  = 1+ \sum_{k=0}^{\infty} (-1)^k L_{ODP,k} \Bigl[1\Bigl] \, z^k \;\;,
\label{six-thirtythree3}\end{eqnarray}
where ODP denotes that only partitions with odd elements are to be considered in the sum over partitions. That is,
$n_{2i} \!=\! 0$ for all values of $i$. The phase factor of $(-1)^k$ in the series expansion arises from the fact that 
only an even number of odd discrete elements yields an even power of $z$, while only an odd number of discrete elements 
yields an odd power of $z$. That is,
\begin{eqnarray}
(-1)^k L_{ODP,k} \Bigl[ 1 \Bigr]= L_{DP,k}\Bigl[ \prod_{i=1}^{[k/2]} (-1)^{n_{2i-1}} \Bigr] \;\;.
\label{six-thirtythree3a}\end{eqnarray}
In a similar fashion we arrive at
\begin{eqnarray} 
\prod_{i=0}^{\infty} \bigl( 1+ z^{2i} \bigr)  = 1+ \sum_{k=1}^{\infty} L_{EDP,2k} \Bigl[1\Bigl]  z^{2k} \;\;,
\label{six-thirtythree4}\end{eqnarray}
where only even elements are to be considered in the even discrete partition operator, i.e $n_{2i+1} \!=\!0$ for
all $i$. Alternatively, we can replace the even discrete partition operator by the discrete partition operator since
\begin{eqnarray}
L_{EDP,2k} \bigl[ 1 \bigr] = L_{DP,k}\bigl[1\bigr] \;\;.
\label{six-thirtythree4a}\end{eqnarray}
Moreover, by equating like powers of $z$ in the power series on both rhs's of Equivalences\ (\ref{six-thirtythree1}) and 
(\ref{six-thirtythree3}), we arrive at 
\begin{eqnarray} 
L_{P,k} \Bigl[ (-1)^{N_k}\Bigr]= (-1)^k L_{ODP,k} \Bigl[ 1 \Bigl]\;\;.
\label{six-thirtythree5}\end{eqnarray}
From this result we see that when $k$ is even, the number of even partitions is greater than the number of odd partitions,
while for odd values of $k$, the opposite applies. Multiplying both sides by $(-1)^k$ results in taking the absolute value 
or modulus of the lhs. Thus, the above statement tells us that absolute value of the difference between the number of 
even and odd partitions is equal to the number of discrete partitions with only odd elements in them or the number of 
partitions with distinct odd parts, a result first proved by Euler according to p.\ 14 of Ref.\ \cite{and03}. 

\begin{figure}
\includegraphics{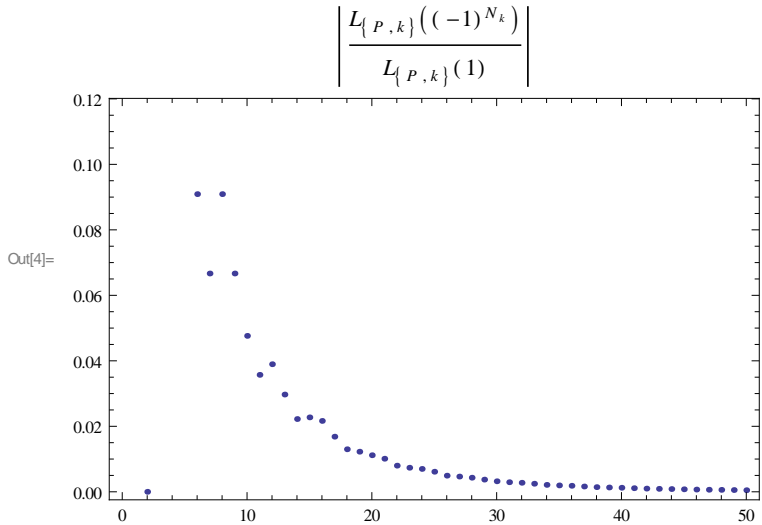}
\caption{The difference between even and odd partitions to the total number of partitions 
summing to $k$ versus $k$}
\label{figfour}
\end{figure}

As a consequence of the previous section, it is a relatively simple exercise to produce a code that evaluates the difference
between the number of partitions with an even number of elements and those with an odd number of elements. Two new global
variables are required, one for evaluating the difference as each partition is scanned and another that is either equal to 1
or -1 depending on whether there is an even number of elements or an odd number of elements. Once the second value is 
determined, it needs to be added to the first global variable in the main function prototype. The second global variable must 
be evaluated in the {\bf termgen} function prototype after the for loop has been altered to calculate the total number elements 
in the partition, which is determined by summing all components of {\it part}. As a consequence, one finds after running 
the code for several values of $k$ that
\begin{eqnarray}
\left| L_{P,k}\Bigl[(-1)^{N_k}\Bigr] \right| \geq \left|L_{P,j}\Bigl[(-1)^{N_j}\Bigr] \right| \;\;,
\label{six-thirtythreee5a}\end{eqnarray}
for $k \! \geq \! j$. Fig.\ \ref{figfour} presents the graph of the ratio of the absolute value of the difference between 
odd and even partitions summing to $k$ to the total number of partitions or $p(k)$ for $k \! \leq \! 50$. Whilst the absolute
value of the $L_{P,k}[(-1)^{N_k}]$ increases with $k$, we see that in relation to the total number of partitions the ratio 
decreases monotonically, once $k$ exceeds 15.  

From Eqs.\ (\ref{six-twentythree}) and (\ref{six-twentyfour}) we see that the infinite products, $\prod_{k=1}^{\infty}
(1+z^k )$ and $\prod_{k=1}^{\infty} \left( 1-z^k \right)$ can be expressed in terms of a series in successive powers 
of $z$, where the coefficients are equal to $L_{DP,k} \left[1 \right]$ and $L_{DP,k} \left[(-1)^{N_k} \right]$ respectively.
The first product can also be written as
\begin{eqnarray}
\prod_{k=1}^{\infty} \Bigl( 1+z^k \Bigr) = \prod_{k=1}^{\infty} \frac{1}{\left( 1-z^{2k-1} \right)} \;\;.
\label{six-thirtythree6}\end{eqnarray}  
Eq.\ (\ref{six-thirtythree6}) is easily obtained by manipulating the rhs after multiplying it by $\left(1-z^{2k}\right)/\left( 
1-z^{2k} \right)$. If we put $C_i \!=\! -1$ in Equivalence\ (\ref{six-thirtythreeb}), then according to Equivalence\ 
(\ref{six-thirtythreed}) and Eq.\ (\ref{six-thirtythreee}) we arrive at
\begin{eqnarray}
\prod_{k=1}^{\infty} \Bigl( 1+z^k \Bigr) \equiv 1+ \sum_{k=1}^{\infty} L_{OEP,k}\Bigl[1 \Bigr] \,z^k \;\;,
\label{six-thirtythree7}\end{eqnarray}
where $L_{OEP,k}[\cdot ]$ represents the odd element partition operator, which we have seen has two different forms
given by Eqs.\ (\ref{five-six}) and (\ref{five-sixa}) depending upon whether $k$ is an even or odd number. Because 
the generating function on the lhs also yields a power series whose coefficients represent the number of discrete 
partitions summing to $k$, we see immediately that
\begin{eqnarray}
L_{DP,k} \Bigl[ 1 \Bigr] = L_{OEP,k} \Bigl[ 1 \Bigr] \;\;.
\label{six-thirtythree8}\end{eqnarray}
Hence, the number of discrete partitions is equal to the number of partitions composed only of odd elements, another 
result attributed to Euler according to p.\ 5 of Ref.\ \cite{and03}. 

By putting $C_i \!=\!-1$ in Eq.\ (\ref{six-nineteen}), we obtained Eq.\ (\ref{six-twentyfour}). If we put $C_i \!=\! -1$
and $z \!=\! z^2$, then the product on the lhs of Eq.\ (\ref{six-twentyfour}) yields a series expansion in powers of 
$z^2$, but now the coefficients represent the number of distinct partitions summing to $2k$ with only even elements operating 
on $(-1)^{N_{2k}}$. Alternatively, this is equivalent to the number of distinct partitions summing to $k$ operating on 
$(-1)^{N_k}$. Therefore, we can write
\begin{align}
\prod_{k=1}^{\infty} \Bigl( 1- z^{2k} \Bigr) & = 1+ \sum_{k=1}^{\infty} L_{EDP,2k}\Bigl[ (-1)^{N_{2k}} \Bigr] \,z^{2k}
\nonumber\\ 
& = \;\; 1+ \sum_{k=1}^{\infty} L_{DP,k}\Bigl[ (-1)^{N_k} \Bigr]\, z^{2k} \;\;.
\label{six-thirtythree9}\end{align}
In Eq.\ (\ref{six-thirtythree9}) $L_{EDP,2k}[\cdot]$ denotes the even discrete partition operator, which acts on the number
of discrete partitions summing to $2k$ where the elements are only even integers. This is opposed to the odd discrete partition 
operator, $L_{ODP,k}[ \cdot ]$, where the elements are odd numbers and discrete, but can sum to both even and odd integers. 
Since the infinite product of $(1-z^{2k})$ is the product of two separate infinite products involving $(1-z^k)$ and $(1+z^k)$,
we have
\begin{eqnarray}
\prod_{k=1}^{\infty} \left( 1 -z^{2k} \right) =  \Bigl(1+ \sum_{k=1}^{\infty} L_{DP,k} \Bigl[ (-1)^{N_k} \Bigr] \, z^k \Bigr) 
\Bigl(1+ \sum_{k=1}^{\infty} L_{DP,k} \Bigl[ 1 \Bigr] \, z^k \Bigr) \;\;.
\label{six-thirtythree10}\end{eqnarray}
Equating like powers of $z$ in Eq.\ (\ref{six-thirtythree10}) with those in Eq.\ (\ref{six-thirtythree9}) yields
\begin{eqnarray}
L_{DP,k}\left[ (-1)^{N_k} \right] =\sum_{j=0}^{2k} L_{DP,j}\Bigl[ (-1)^{N_j}\Bigr] \; L_{DP,2k-j} \Bigl[ 1 \Bigr] \;\;,
\label{six-thirtythree11}\end{eqnarray} 
and
\begin{eqnarray}
\sum_{j=0}^{2k+1} L_{DP,j}\Bigl[ (-1)^{N_j}\Bigr] \; L_{DP,2k+1-j} \Bigl[ 1 \Bigr] =0 \;\;.
\label{six-thirtythree12}\end{eqnarray} 
Since $L_{DP,k}[(-1)^{N_k}]$ and $L_{DP,k}[1]$ are equal to $q(k)$ and $q(k,1)$ respectively, Eqs.\ (\ref{six-thirtythree11}) and
(\ref{six-thirtythree12}) can also be written as
\begin{eqnarray}
q(k) \Bigl( 1- q(k,1) \Bigr) = \sum_{j=0}^{k-1} \Bigl( q(j) q(2k-j,1) + q(2k-j) q(j,1) \Bigr) \;\;,
\label{six-thirtythree13}\end{eqnarray}
and 
\begin{eqnarray}
\sum_{j=0}^{k} \Bigl( q(j) q(2k+1-j,1) + q(2k+1-j) q(j,1) \Bigr) =0 \;\;.
\label{six-thirtythree14}\end{eqnarray}
In these results it should be borne in mind that $q(0) \!=\! q(0,1) \!=\! 1$. Isolating the $j \!=\!0$ terms in the
above equations yields
\begin{align}
q(2k,1) + q(2k) & =  q(k) \Bigl( 1- q(k,1) \Bigr) -  \sum_{j=1}^{k-1} \Bigl( q(j) \, q(2k-j,1) 
\nonumber\\
& + \;\; q(2k-j) \, q(j,1) \Bigr) \;\;,
\label{six-thirtythree15}\end{align}
and 
\begin{align}
q(2k+1,1) + q(2k+1) & =  \sum_{j=1}^{k} \Bigl( q(j) \, q(2k+1-j,1) 
\nonumber\\
& + \;\; q(2k+1-j) \, q(j,1) \Bigr)  \;\;.
\label{six-thirtythree16}\end{align}
Eqs.\ (\ref{six-thirtythree15}) and (\ref{six-thirtythree16}) represent the recurrence relations for determining
the number of discrete partitions or $q(k,1)$. Like the Euler/MacMahon recurrence relation given by Eq.\ (\ref{six-nine}),
they utilise the discrete partition numbers or $q(k)$ and consequently, most of the terms in the sums vanish
when the summation index $j$ is not equal to a pentagonal number.

It should also be noted that the analysis resulting in Eq.\ (\ref{six-twelvea}) can be adapted to provide another 
representation for the partition number polynomials or $p(k,\omega)$. First, we re-write the generalised product in 
Equivalence\ (\ref{six-thirtytwo}) as
\begin{eqnarray}
P(z,\omega) = \exp \Bigl( \sum_{m=1}^{\infty} \sum_{j=1}^{\infty} \omega ^j \, z^{mj}/j \Bigr) \;\;,
\label{six-thirtyfoura}\end{eqnarray}  
where now it is assumed that $|\omega z| \!<\! 1$. Consequently, the modified version of Eq.\ (\ref{six-twelvea})
becomes
\begin{align}
P(z,\omega) & = 1 +  \sum_{k=1}^{\infty} \frac{1}{k!} \Bigl( \omega z+(\omega+\omega^2/2)\,z^2+ 
(\omega+ \omega^3/3) \, z^3
\nonumber\\
& + (\omega+ \omega^2/2+\omega^4/4)\, z^4 + \cdots + \gamma_j(\omega) z^j+ \cdots \Bigr)^k \, ,
\label{six-thirtyfourb}\end{align} 
where $\gamma_j(\omega)= \sum_{d|j} (d/j)\,\omega^{j/d}$ and $d$ represents a divisor of $j$ as before. The first
few divisor polynomials are: $\gamma_0(\omega) =1$, $\gamma_1(\omega)= \omega$, $\gamma_2(\omega)=\omega +\omega^2/2$,
and $\gamma_3(\omega)= \omega+ \omega^3/3$. As expected, for $\omega=1$ the $\gamma_j(\omega)$ reduce to the 
$\gamma_j$ below Eq.\ (\ref{six-eleven}). This means that Eq.\ (\ref{six-thirteen}) can be generalised to
\begin{eqnarray}
q(k,-\omega) = L_{P,k} \Bigl[ (-1)^{N_k} \prod_{i=1}^k \frac{\gamma_i(\omega)^{n_i}}{n_i!} \Bigr] \;\;.
\label{six-thirtyfourb1}\end{eqnarray}
Furthermore, the $\gamma_j(\omega)$ are polynomials in $\omega$ whose highest and lowest orders are respectively
$j$ and unity. That is, deg $\gamma_j(\omega) \!=\! j$. We shall refer to these unusual polynomials as the divisor 
polynomials. By applying Theorem\ 1 to Eq.\ (\ref{six-thirtyfourb}) with the coefficients of the inner and outer 
series set equal to $\gamma_k(\omega)$ and $1/k!$ respectively, we arrive at
\begin{eqnarray}
p(k,\omega) = L_{P,k} \Bigl[ \prod_{i=1}^k \frac{\gamma_i(\omega)^{n_i}}{n_i!} \Bigr] \;\;.
\label{six-thirtyfourc}\end{eqnarray}  
Since $ {\rm deg}\;\gamma_i(\omega) \!= \! i$ and $\prod_{i=1}^k \omega^{i n_i} \!=\! \omega^k$, the highest 
order term in the $p(k,\omega)$ is $k$. On the other hand, since the lowest order term in the $\gamma_i(\omega)$ 
is $\omega$, $\prod_{i=1}^k \omega^{n_i}\!=\!\omega^{N_k}$ and $N_k$ ranges from unity to $k$, the lowest order 
term in the $p(k,\omega)$ is unity, again confirming that the partition number polynomials are polynomials in $\omega$
with deg $p(k,\omega) \!=\! k$.

The sixth and last program presented in the appendix is called ${\bf dispfnpoly}$. It prints out both $q(k,-\omega)$ 
and $p(k,\omega)$ in symbolic form so that it can be handled by Mathematica \cite{wol92}. To run this program, the 
user must specify the order $k$ of the polynomials. The program is different from the other programs in the
appendix since it is not required to determine the factorial of the total number of distinct parts, i.e.\ $N_k !$. 
When the global variable ${\it polytype}$ is equal to unity in the for loop in ${\bf main}$, the program 
determines the discrete partition polynomial for the order specified by the user. When it becomes equal to 2, the 
program determines the corresponding partition function polynomial. E.g., for $k \!=\! 6$, the following output is
generated:\newline
Q[6,-w$_{-}$]:= DP[6,w] (-1) + DP[1,w] DP[5,w] (-1)$^{\wedge}$(2) + \newline 
DP[1,w]$^{\wedge}$(2) DP[4,w] (-1)$^{\wedge}$(3)/2! + DP[1,w]$^{\wedge}$(3) DP[3,w] (-1)$^{\wedge}$(4)/3! + 
DP[1,w]{$^{\wedge}$(4) DP[2,w] (-1)$^{\wedge}$(5)/4! + \newline 
DP[1,w]$^{\wedge}$(6) (-1)$^{\wedge}$(6)/6! + DP[1,w]$^{\wedge}$(2) DP[2,w]$^{\wedge}$ (2) (-1)$^{\wedge}$(4)/(2! 2!) 
+ DP[1,w] DP[2,w] DP[3,w] (-1)$^{\wedge}$(3) + \newline
DP[2,w] DP[4,w] (-1)$^{\wedge}$(2) + DP[2,w]$^{\wedge}$(3) (-1)$^{\wedge}$(3)/3! + DP[3,w]$^{\wedge}$(2) 
(-1)$^{\wedge}$(2)/2! \newline
\newline
P[6,w$_{-}$]:= DP[6,w]  + DP[1,w] DP[5,w]  + \newline
DP[1,w]$^{\wedge}$(2) DP[4,w] /2! + DP[1,w]$^{\wedge}$(3) DP[3,w] /3! + DP[1,w]$^{\wedge}$(4) DP[2,w] /4! + \newline
DP[1,w]$^{\wedge}$(6) /6! + DP[1,w]$^{\wedge}$(2) DP[2,w]$^{\wedge}$(2) /(2! 2!) + DP[1,w] DP[2,w] DP[3,w]  + \newline
DP[2,w] DP[4,w]  + DP[2,w]$^{\wedge}$(3) /3! + DP[3,w]$^{\wedge}$(2) /2!\newline
Time taken to compute the coefficients is 0.000000 seconds \newline
The terms on the rhs denoted by DP[$k$,w] represent the divisor polynomials of order $k$. These can be obtained by
typing in the following line in Mathematica: \newline
\indent   DP[k$_{-}$,w$_{-}$]:= Sum[w$^{\wedge}$d/d,\{d,Divisors[k]\}]. 

From Eqs.\ (\ref{six-thirtyfourb1}) and (\ref{six-thirtyfourc}) we see that the discrete partition polynomials or 
rather the $q(k,-\omega)$ are almost identical to the partition function polynomials or $p(k,\omega)$ when they are both
expressed in terms of the divisor polynomials, the only difference being a phase factor that appears in the former.
This factor is positive when the number of elements in a partition is even, but is negative when there is an
odd number of elements in a partition. Hence, the only difference between the two sets of polynomials occurs
for odd partitions or those with an odd number of elements.

If we multiply Eq.\ (\ref{six-twentyfive}) with Eq.\ (\ref{six-thirtytwo}), then by equating like powers on both
sides of the resulting equation we obtain
\begin{eqnarray}
\sum_{j=0}^{k} q(j,-\omega) \, p(k-j,\omega) = 0 \;\;.
\label{six-thirtyfive} \end{eqnarray}
Alternatively, Eq.\ (\ref{six-thirtyfive}) can be expressed as
\begin{eqnarray}
\sum_{j=0}^{k} L_{P,j} \Bigl[ (-1)^{N_j} \prod_{i=1}^j \frac{\gamma_i(\omega)^{n_i}}{n_i!} \Bigr] 
L_{P,k-j} \Bigl[ \prod_{i=1}^{k-j} \frac{\gamma_i(\omega)^{n_i}}{n_i!} \Bigr] = 0 \;\;.
\label{six-thirtyfivea} \end{eqnarray}
These results represent the generalisation of the Euler/MacMahon recurrence relation given by Eq.\ (\ref{six-nine}).

We can also generalise the product in Equivalence\ (\ref{six-thirtytwo}) by introducing the arbitrary power of $\rho$ 
as we did in the discrete partition case of Equivalence (\ref{six-twentyeightb}). Thus, the generalised product
becomes
\begin{eqnarray}
P_{\rho}(z,\omega) = \prod_{k=1}^{\infty} \frac{1}{(1-\omega z^k)^{\rho}} \equiv 1+\sum_{k=1}^{\infty} p(k,\omega,\rho)
\, z^k \;\;.
\label{six-thirtyfiveb}\end{eqnarray}
In this case the coefficients $p(k,\omega,\rho)$ in the generating function can be determined by setting the $D_i$ 
equal to $p(i,\omega)$ in Eq.\ (\ref{thirtyfive}). Then one finds that
\begin{eqnarray}
p(k,\omega,\rho)= L_{P,k} \Bigl[ (-1)^{N_k} \, (-\rho)_{N_k} \prod_{i=1}^k \frac{p(i,\omega)^{n_i}}{n_i!} \Bigr] \;\;.
\label{six-thirtyfivec}\end{eqnarray}

As a consequence of the preceding analysis, we are now in a position to study more advanced products. In particular,
let us consider the following quotient: 
\begin{eqnarray}
P(z,\beta \omega, \alpha \omega) = Q(z,-\beta \omega) P(z,\alpha \omega)= \prod_{k=1}^{\infty} \frac{\left( 1- 
\beta \omega z^k \right)} {\left( 1 -\alpha \omega z^k\right)} \;\;. 
\label{six-thirtyfived}\end{eqnarray}
In deriving a power series expansion or generating function for the above product we expect the power of $\omega$ in
the coefficients to yield the total number of elements in the partitions summing to the power of $z$. Furthermore, the 
power of $\beta$ in the coefficients should represent the number of elements due to the discrete partitions, while the 
power of $\alpha$ should give the number of elements due to the standard partitions. By adopting the same approach as 
for the other infinite products that have already been presented in this section, we can express Eq.\ (\ref{six-thirtyfived}) as
\begin{align}
P(z,\beta \omega,\alpha \omega) & = \sum_{k=0}^{\infty} q(k,-\beta \omega) z^k \sum_{k=0}^{\infty} 
p(k,\alpha \omega) z^k 
\nonumber\\
& = \;\; \sum_{k=0}^{\infty} z^k \; QP_k(\omega, \beta,\alpha) 
\label{six-thirtysix}\end{align}
where 
\begin{eqnarray}
QP_k(\omega,\beta,\alpha) = \sum_{j=0}^k q(j,-\beta \omega) \, p(k-j,\alpha \omega) \;\;. 
\label{six-thirtyseven}\end{eqnarray}
From Eq.\ (\ref{six-thirtyseven}) we see that $QP_0(\omega,\beta,\alpha) \!=\! 1$. Furthermore, since $P(z,\alpha \omega,\alpha \omega)$ 
equals unity, it follows that $QP_k(\omega,\alpha,\alpha) \!=\! 0$ for $k \!>\!0$. 

\begin{table}
\begin{tabular}{|c|l|} \hline
$k$ &  $QP_k(\omega,\alpha,\beta)$  \\ \hline
$0$ & $ 1$ \\  
$1$ &  $ (\alpha-\beta) \omega $  \\ 
$2$ & $ (\alpha-\beta) (\omega + \alpha \omega^2)$ \\ 
$3$ & $(\alpha-\beta) (\omega + (\alpha-\beta) \omega^2+ \alpha^2 \omega^3)  $  \\ 
$4$ & $(\alpha-\beta) (\omega + (2\alpha-\beta) \omega^2+ \alpha (\alpha-1) \omega^3 +\alpha^3 \omega^4)   $ \\ 
$5$ & $(\alpha-\beta) (\omega +2 (\alpha-\beta) \omega^2+ 2\alpha (\alpha-\beta) \omega^3 +\alpha^2
(\alpha-\beta) \omega^4 +\alpha^4 \omega^5)  $ \\ 
$6$ & $(\alpha-\beta) (\omega + (3\alpha-2\beta) \omega^2+ (3\alpha^2-3\alpha \beta+\beta^2) \omega^3 +2\alpha^2
(\alpha-\beta) \omega^4 $\\ 
$ $ & $ + \;\; \alpha^3(\alpha-\beta) \omega^5+ \alpha^5 \omega^6) $ \\  
$7$ & $ (\alpha-\beta) (\omega +3 (\alpha-\beta) \omega^2+ (4\alpha- \beta)(\alpha-\beta) \omega^3 +(3\alpha-\beta)
(\alpha-\beta) \alpha\omega^4 $\\ 
$ $ & $ + \;\;2 \alpha^3(\alpha-\beta) \omega^5+ \alpha^4 (\alpha-\beta) \omega^6+ \alpha^6 \omega^7) $ 
\\ 
$8$ & $ (\alpha-\beta) (\omega +(4 \alpha-3\beta) \omega^2+ (5\alpha- 2 \beta)(\alpha-\beta) \omega^3 +(5\alpha^2-6 \alpha 
\beta +2 \beta^2) \alpha \omega^4 $\\
$ $ & $ + \;\;\alpha^2(3\alpha-\beta) (\alpha -\beta) \omega^5+ 2\alpha^4 (\alpha-\beta) \omega^6+ \alpha^5(\alpha-\beta) 
\omega^7 +\alpha^7 \omega^8) $ \\ \hline

\end{tabular}
\caption{Coefficients $QP_k(\omega,\beta,\alpha)$ in the power series expansion for the generating function given by Eq.\ 
(\ref{six-thirtyfive}).}
\label{table3}
\end{table}
 
Table\ \ref{table3} presents the coefficients $QP_k(\omega,\beta,\alpha)$ up till $k \!= \!8$. As can be seen
from the table they are polynomials of $O(k)$ in $\omega$. The power of $\omega$ in these polynomials gives the 
number of elements in the final partitions, which combine the elements from standard partitions with those from 
discrete partitions. As expected, the polynomials vanish when $\alpha \!=\! \beta$ since the product 
$P(z,\beta \omega,\alpha \omega)$ equals unity in this case. Furthermore, the highest power of $\alpha$ is 
$k$, which is also the highest power of $\omega$. This corresponds to the fact that the power of $\alpha$ 
represents the number of elements in the standard partitions. Therefore, the greatest number of elements in the final
partitions will be due to the partition \{1,1,...,$1_k$\} with no elements coming from a discrete partition. The 
highest power of $\beta$, however, is considerably lower since it is determined by the partition with the most number 
of discrete elements summing to $k$. In this instance the power of $\alpha$ will be zero. E.g., for $k \!=\! 8$, the 
highest power of $\beta$ is three, which is in accordance with Eq.\ (\ref{six-twentysevena}). When $\alpha \!=\!0$, 
the polynomials reduce to the polynomials arising from the generating function for discrete partitions, i.e.\ 
$q(k,-\beta \omega)$, while for $\beta \!=\!0$, they reduce to the partition function polynomials or $p(k,\alpha \omega)$. 
In addition, for $\alpha \!=\!0$ the coefficients in the resulting polynomials give the number of discrete partitions 
where the number of elements is equal to the power of $\beta$. For $\beta \!=\!0$ the coefficients of the resulting 
polynomials become the sub-partition numbers or $p_i(k)$, where $i$ represents the power of $\alpha$. 

The interesting terms in the polynomials displayed in Table\ \ref{table3} are the cross-terms involving $\alpha$ and 
$\beta$, which represent the mixture of the discrete partitions and standard partitions with the total number of the 
elements equal to the power of $\omega$. For example, in $QP_3(\omega,\alpha,\beta)$ the coefficient of $\omega^2$ has 
a term equal to $-2 \alpha \beta$, which tells us that one element in the partition has come from the generating function 
for standard partitions and the other has come from the generating function for discrete partitions. There are 
two instances where this can occur: either the one has come from the discrete part or numerator on the rhs of Eq.\ 
(\ref{six-thirtyfive}) and the two from the denominator or vice-versa. On the other hand, there is only one instance of 
the partition \{1,2\} emanating only from either the numerator or the denominator. Thus, the coefficients of $\omega^2$ 
in $QP_3(\omega,\alpha,\beta)$ for only standard and discrete partitions are respectively $\alpha^2$ and $\beta^2$. 
In this instance the partition maintains its discreteness when accepting an element from the discrete partitions 
and one from the standard partitions, but this will not always be the case. In addition, the power of $\alpha$ can be 
much higher than the power of $\beta$ reflecting the fact that the greatest number of elements in a discrete partition 
summing to a particular value is significantly less than the greatest number of elements in a standard partition summing
to the same value.

According to p.\ 23 of Ref.\ \cite{and03}, Gauss derived the following result:
\begin{eqnarray}
1+ 2 \sum_{k=1}^{\infty} (-1)^k\, z^{k^2} = \prod_{k=1}^{\infty} \frac{(1-z^k)}{(1+z^k)} \;\;.
\label{six-thirtysevena}\end{eqnarray}
The rhs of the above result is a special case of $P(z,\beta \omega,\alpha \omega)$, namely $P(z,1,-1)$. If the values
for $\alpha \omega$ and $\beta \omega$ are introduced into the rhs of Eq.\ (\ref{six-thirtysix}), then we can equate
like powers of $z$ with the lhs of the Eq.\ (\ref{six-thirtysevena}). Consequently, for $i$ equal to a positive integer
we arrive at
\begin{eqnarray}
QP_k(1,1,-1) = \sum_{j=0}^k q(j)\, p(k-j,-1)= \begin{cases} 2(-1)^i, & \;\; k=i^2, \\
0,  & \;\; {\rm otherwise} \;\;. \end{cases}
\label{six-thirtysevenb}\end{eqnarray}
Eq.\ (\ref{six-thirtysevenb}) can be checked with the results appearing in Table\ \ref{table3}.

We can also use the preceding analysis to derive a power series expansion or generating functon for the product of 
two specific forms of $P(z,x,y)$ involving the three parameters, $\omega$, $x$, and $y$, and the variable, $z$. This
product was first studied by Heine. According to p.\ 55 of Ref.\ \cite{knu005} he found that
\begin{eqnarray}
\prod_{k=1}^{\infty} \frac{(1-\omega x z^k)}{(1- \omega z^k)}  \frac{(1-\omega y z^k)}{(1- \omega x y z^k)} =
\sum_{k=0}^{\infty} \frac{(1/x;z)_k}{(z;z)_k}  \frac{(1/y;z)_k}{(\omega z ;z)_{k+1}} \,(\omega x y z)^k .
\label{six-thirtyeight}\end{eqnarray}
If we set $a \!=\!1/x$, $b \!=\! 1/y$, $c \!=\! \omega z$ and $q \!=\! z$ with $|c/ab| \!<\! 1$ and $|q| \!<\! 1$, which are
the conditions for guaranteeing absolute convergence, then the above result can be expressed as a q-hypergeometric series.  
This is perhaps the more familiar form for the product, where it is written as
\begin{eqnarray}
\sum_{k=0}^{\infty} \frac{(a;q)_k}{(q;q)_k}  \frac{(b;q)_k}{(c;q)_k} \Bigl(\frac{c}{ab} \Bigr)^k =
\prod_{k=0}^{\infty} \frac{\bigl(1-(c/a) q^k \bigr)}{\bigl(1- c q^k \bigr)}  \frac{\bigl(1-(c/b) q^k \bigr)}{\bigl(1- (c/ab) q^k \bigr)} \;\;.
\label{six-thirtyeighta}\end{eqnarray}
This result appears as Corollary 2.4 on p.\ 20 of Ref.\ \cite{and03}.

The lhs of Eq.\ (\ref{six-thirtyeight}) represents the product of $P(z,\omega x,\omega)$ and $P(z,\omega y,\omega x y)$. If we 
denote the product of $P(z,x,y)$ and $P(z,s,t)$ by $P_2(z,x,y,s,t)$ and introduce Eq.\ (\ref{six-thirtysix}), then we find that
\begin{align}
P_2(z,\omega x,\omega,\omega y,\omega x y) & = \prod_{k=1}^{\infty} \frac{(1-\omega x z^k)}{(1- \omega z^k)}  
\frac{(1-\omega y z^k)}{(1- \omega x y z^k)} 
\nonumber\\
&= \;\; \sum_{k=0}^{\infty} z^k \, HP_k(\omega,x,y) \;\;,
\label{six-thirtynine}\end{align}  
where  
\begin{eqnarray}
HP_k(\omega,x,y)= \sum_{j=0}^k QP_j(\omega,x,1) \, QP_{k-j}(\omega, y, xy) \;\;.
\label{six-forty}\end{eqnarray}  
The equals sign appears here because of the conditions given below Eq.\ (\ref{six-thirtyeight}).

\begin{table}
\begin{tabular}{|c|l|} \hline
$k$ &  $HP_k(\omega,x,y)$  \\ \hline
$0$ & $ 1$ \\  
$1$ &  $ (x-1)(y-1) \omega $  \\ 
$2$ & $ (x-1)(y-1)\omega \bigl(1+ (1 +x y )\omega \bigr)$ \\ 
$3$ & $ (x-1)(y-1)\omega \bigl(1+(x-1)(y-1) \omega +(1+ x y +x^2y^2) \omega^2 \bigr)$ \\ 
$4$ & $ (x-1)(y-1)\omega \bigl(1+(2-x-y+2xy) \omega +(x-1)(y-1)(1+xy) \omega^2$ \\
$$  & $ +(1+ x y +x^2y^2+x^3y^3) \omega^3 \bigr)$ \\ 
$5$ & $ (x-1)(y-1)\omega \bigl(1 +2(x-1)(y-1) \omega +2(x-1)(y-1)(1+xy) \omega^2 $\\
$$ & $ + (x-1)  (y-1)(1+ x y +x^2y^2+x^3y^3) \omega^3 +(1+ x y +x^2y^2+x^3y^3$ \\
$$ & $ +x^4 y^4) \omega^4 \bigr)$ \\ 
$6$ & $ (x-1)(y-1)\omega \bigl( 1+(3-2x-2y+3xy) \omega +(3-3x-3y+x^2+y^2$ \\
$$ & $ +21 x y-3x y^2-3 x^2 y+3 x^2 y^2) \omega^2 + (x-1) (y-1)(2+ 3 x y +2x^2y^2) \omega^3$\\
$$ & $ +(x-1)(y-1) (1+ x y +x^2y^2 +x^3y^3) \omega^4 +(1+ xy+ x^2 y^2+ x^3 y^3 $ \\
$$ & $ + x^4 y^4+x^5 y^5) \omega^5 \bigr)$ \\ \hline
\end{tabular}
\caption{Coefficients of the polynomials $HP_k(\omega,\alpha,\beta)$ arising from the three-parameter one variable generating 
function given in Eq.\ (\ref{six-thirtynine}).}
\label{table4}
\end{table}

Table\ \ref{table4} presents the coefficients $HP_k(\omega,x,y)$ up to $k \!=\! 6$. They have been obtained by implementing
Eq.\ (\ref{six-forty}) in Mathematica. From the table it can be seen that the $HP_k(\omega,x,y)$ are polynomials in $\omega$ 
of degree $k$. The coefficient of the leading order term is 
\begin{eqnarray}
C^{HP}_k= \frac{\left( 1- x^k y^k \right)}{(1-xy)} \; (x-1) (y-1)\;\;,
\label{six-fortyone}\end{eqnarray}
while that for the penultimate leading order term is found to be
\begin{eqnarray}
C^{HP}_{k-1}= \frac{\left( 1- x^{k-2}y^{k-2} \right)}{(1-xy)} \; (x-1)^2 (y-1)^2\;\;.
\label{six-fortytwo}\end{eqnarray}
In the above equation $k \!>\! 2$. The lowest order term in $\omega$ for these polynomials is linear and its coefficient is
equal to
\begin{eqnarray}
C^{HP}_1= (x-1)(y-1) \; \;.
\label{six-fortythree}\end{eqnarray}
As expected, the polynomials are zero when either $x \!=\! 1$ or $y \!=\!1$ since in these cases 
$P(z,\omega x,\omega,\omega y,\omega x y)$ is equal to unity. In addition, they are symmetrical in $x$ and $y$ 
in the sense that any power of $\omega$ with a term of $\alpha x^i y^j$, where $i \! \neq \!j$, in its 
coefficient will also possess the term of $\alpha x^j y^i$ with the same power of $\omega$.  

\section{Conclusion}
Originally, this work set out to devise a programming methodology on the partition method for a power series 
expansion, which has been used recently to solve important and intractable problems in applied mathematics
\cite{kow10}-\cite{kow11a}. In this method the coefficients of the resulting power series for a function are 
obtained by summing the contributions made by each partition summing to the order $k$. These contributions are 
evaluated by: (1) assigning values $p_i$ to each element $i$ in a partition, (2) multiplying by a multinomial factor 
composed of the factorial of the total number of elements  in the partition, $N_k!$ divided by the factorial of each 
element's frequency, $n_i!$ and (3) multiplying by the coefficient of an outer series for the total number of elements
in the partition, viz.\ $q_{N_k}$. In order to apply the method, it means that one needs to know the composition 
of all the partitions summing to $k$, which includes the frequencies or numbers of occurrences of all the 
elements in each partition. Therefore, an algorithm is required that is capable of scanning all these partitions, 
the number of which increase exponentially with $k$. Whilst Sec.\ 2 discusses various methods of generating 
partitions, it turns out that the novel bi-variate recursive central partition or BRCP algorithm is the most 
suitable method for implementation in the partition method for a power series expansion because it is based on 
the graphical representation of the partitions in the form of a non-binary tree diagram as depicted for $k \!=\! 6$ 
in Fig.\ \ref{figone}. As a consequence, the BRCP algorithm is able to print out partitions in the multiplicity 
representation more efficiently than the other algorithms discussed in Sec.\ 2, while the multiplicity representation 
turns out to be the minimum amount of information required for carrying out the method for a power series expansion.  

The theory behind the partition method for a power series expansion is presented in Sec.\ 3 as Theorem\ 1, which 
shows how power series expansions can be derived from a quotient of pseudo-composite functions. It also represents
the lynchpin of this work. Next the regularisation of the binomial series is presented in Lemma\ 1. With this result
a corollary to Theorem\ 1 is developed, whereby the partition method is adapted to the situation in which the 
quotient of the pseudo-composite functions can be taken to an arbitrary power. As a result of Theorem\ 1 and its
corollary, we observe that the process of evaluating the contributions made by each partition can be viewed as a 
discrete operation, giving rise to a partition operator denoted by $L_{P,k}[\cdot]$. While $L_{P,k}[1] \!=\! p(k)$ or 
the number of partitions summing to $k$, varying the argument inside the operator yields completely different 
identities. Moreover, the partition operator can be modified so that it only applies to specific types of partitions 
such as discrete or odd/even partitions, again resulting in further new and fascinating identities when the arguments 
are altered.

Because the number of partitions increases exponentially, it becomes rather onerous to apply the partition method
for a power series expansion when the order $k$ is greater than 10. This problem is overcome by modifying 
the BRCP algorithm to calculate the contribution due to each partition in symbolic form. Sec.\ 4 presents two such 
programs. The first {\bf partmeth} calculates all the coefficients $D_k$ and $E_k$ in Theorem\ 1 up to and including the 
value of $k$ specified by the user. Unfortunately, for $k \!>\!20$ the output files generated by this program become
too large. Then we require a code that only evaluates the coefficients $D_k$ and $E_k$ for a particular value of $k$, 
which is accomplished by the second code {\bf mathpm}. For much larger values of $k$, say for $k \!>\! 100$, storing the 
coefficients is no longer a viable option. In these cases the output files need to be divided into smaller files. Then
each file can be evaluated separately in Mathematica and the result stored. Once a result is stored, the file can be 
deleted and another file can be imported into the software package. Once it is evaluated and the result stored, it too 
can be deleted. Finally, the stored results can be summed to yield the value for the coefficient. 

As a result of the success in developing a programming methodology for the partition method for a power series expansion, 
we have been able to create programs that can determine various types of integer partitions such as those with either a 
fixed number of elements or specific elements, doubly-restricted partitions and discrete/distinct partitions. Normally, 
different programming approaches are required to solve each of these problems, but as explained in Sec.\ 5, they can all 
be solved by introducing minor modifications to the program {\bf partgen} in Sec.\ 2. In the process new operators such 
as the discrete partition operator $L_{DP,k}[ \cdot]$ and the odd- and even-element partition operators, $L_{OEP,k}[ \cdot]$ 
and $L_{EEP,k}[\cdot]$ are defined. In particular, the number of discrete partitions summing to $k$ or $q(k,1)$ is equal 
to $L_{DP,k}[1]$. Another interesting application in this section is the development of the program {\bf transp}, which 
determines conjugate partitions by means of Ferrers diagrams. These are created by the dynamic memory allocation of 
two-dimensional arrays in the C/C++ code appearing in the appendix. 

In Secs.\ 6-8 the operator approach of Secs.\ 3 and 5 is employed in the derivation of new power series expansions or
generating functions for numerous infinite products of increasing complexity that arise in the theory of partitions. 
Sec.\ 6 begins by studying the product $P(z)$ defined by Equivalence\ (\ref{fortyfoura}), which produces a generating 
function or power series expansion whose coefficients are the partition function or $p(k)$. Theorem\ 2 shows that the
generating function of this important product is absolutely convergent for $|z| \!<\! 1$ and divergent elsewhere. That is,
unlike the geometric series, there is no region in the complex plane where the generating function is conditionally
convergent. Instead, there is a ring of singularity separating the absolutely convergent unit disk from the rest of the
divergent complex plane. As a result of Theorem\ 2, Equivalence\ (\ref{fortyfoura}) is only an equation when $|z| \!<\!1$. 
For these values of $z$, the generating function can be inverted and Theorem\ 1 can then be applied to the ensuing result. 
The coefficients of the resulting power series expansion or generating function are referred to as the discrete partition 
numbers $q(k)$, which equal $(-1)^j$ when $k$ is a pentagonal number or equal to $(3j^2\pm j)/2$ and zero, otherwise. We 
also find that they can be expressed in terms of the partition operator with each element $i$ assigned to the partition 
function value of $p(i)$ as given by Eq.\ (\ref{six-five}). On the other hand, when Theorem\ 1 is applied to the inverse 
of $P(z)$, we find that the discrete partition numbers can be expressed in terms of the partition operator acting with 
each element $i$ assigned to a value of $p(i)$ as given by Eq.\ (\ref{six-fifteen}). 

Sec.\ 6 continues with the derivation of alternative representations for the generating functions of $P(z)$ and its inverse,
which are given by Eqs.\ (\ref{six-twelvea}) and (\ref{six-eleven}), respectively. In these results the coefficients for 
each power $j$ of the inner series are expressed in terms of a sum over the divisors $d$ of $j$ divided by $j$ and are 
denoted by $\gamma_j$. When Theorem\ 1 is applied to these new forms for the products, the coefficients of the generating 
functions now have the partition operator acting with the elements $i$ are assigned to values of $\gamma_i$ as in Eqs.\ 
(\ref{six-thirteen}) and (\ref{six-sixteena}). The difference between these results is the appearance of the phase 
factor $(-1)^{N_k}$ in the case of the discrete partition numbers.

Because a result like Eq.\ (\ref{six-fifteen}) gives the partition numbers in terms of the partition operator acting 
with the elements $i$ assigned to the discrete partition numbers or $q(i)$, it becomes necessary to develop a program that 
excludes those partitions in which the $q(i)$ vanish. This represents a completely different type or class of partition 
studied in Sec.\ 5. Thus, Sec.\ 6 presents the program called ${\bf partfn}$, which describes the modifications that need 
to be made to ${\bf partgen}$ in Sec.\ 2 so that only those partitions in which all the elements are pentagonal numbers 
are printed out. 

Although the programs in Sec.\ 4 were developed with the partition operator acting on all partitions, it was stated 
at the end of Sec.\ 5 that they could be adapted to handle situations where only a subset of the total number
of partitions is required. That is, the programming methodology in Sec.\ 4 is not restricted to the partition operator,
but as a result of the material in Sec.\ 5, it can be adapted to handle situations in which the coefficients are
expressed in terms of different operators or specific types of partitions. Such a situation occurs with Eq.\ 
(\ref{six-fifteen}). Consequently, Sec.\ 6 concludes by presenting the program ${\bf pfn}$, which determines the 
contributions made by the partitions whose elements yield non-zero discrete partition numbers. This program expresses 
the partition function or $p(k)$ in two symbolic forms, both of which must be imported into Mathematica to obtain 
the actual values for $p(k)$. In the first form the partition function is expressed directly in terms of the discrete 
partition numbers or $q(i)$. This means that in order to obtain the values for the partition function, another module 
for calculating the discrete partition numbers must be created in Mathematica. In the second form the non-zero 
values of the $q(i)$ are replaced by their symbolic form of (-1)$^{\wedge}$j for $i \!=\! (3j^2\pm j)/2$. When the 
second form is entered as input into Mathematica, it immediately computes the partition function in integer form.

Sec.\ 7 begins with the presentation of Theorem\ 3, which generalises the infinite product of $1/P(z)$. Instead of the 
coefficients of the powers of $z$ being equal to -1, they are now assumed to be equal to general values $C_k$ as in Eq.\
(\ref{six-nineteen}). The theorem is proved by adapting the proof of Theorem\ 1 and means that the coefficients of the 
generating series for this product are given in terms of the discrete partition operator acting with each element $i$ 
assigned the value of $C_i$ as in Eq.\ (\ref{six-twentyone}). Conversely, the theorem implies that any power series 
expansion for a function can be expressed as an infinite product. After some elementary examples are studied, viz.\
the geometric series and the exponential power series, the $C_i$ are set equal to $\omega$, whereupon we find that the 
coefficients of the resulting generating function become the polynomials $q(k,\omega)$, which are referred to as the 
discrete partition polynomials. As expected, for $\omega \!=\!-1$ they reduce to the discrete partition numbers, 
while for $\omega \!=\! 1$, they yield the number of discrete partitions summing to $k$. In fact, the coefficients 
of the discrete partition polynomials represent the number of discrete partitions where the number of elements is equal to  
the power of $\omega$. Then a brief description appears on how the program ${\bf dispart}$ presented in Sec.\ 5 can be 
adapted to evaluate these numbers, which can also be written as $L_{DP,k,i}[1]$.
 
The infinite product of Theorem\ 3 is further generalised by the introduction of an arbitrary power $\rho_k$ as given by
Eq.\ (\ref{six-twentyeightc1}) in the corollary to Theorem\ 3. In this case the coefficients of the resulting generating
function are expressed in terms of the partition operator acting with the elements $i$ assigned a value of $-C_i$ multiplied
by $(-\rho_i)_{n_i}$, where $n_i$ is the number of occurrence of $n_i$ and $(\rho)_k$ denotes the Pochhammer symbol.
With the aid of other results appearing in the corollary, we are able to study the generating functions when $1/P(z)$ is 
squared or cubed. For these cases $\rho_i$ has been set equal to the uniform value of $\rho$ and the discrete partition 
polynomials $q(k,\omega)$ are extended to become $q(k,\omega,\rho)$ as defined by Eq.\ (\ref{six-twentyeightc10}). Several 
identities involving the $q(k,\omega,\rho)$ are also derived.   

Sec.\ 8 is also devoted to deriving generating functions for other forms of infinite products using the material from the
previous two sections. First of all, the coefficients of the powers of $z$ in $P(z)$ are set equal to $-\omega$ instead of -1.
Consequently, the coefficients of the generating function obtained after applying Theorem\ 1 to the modified version of 
$P(z)$ or $P(z,\omega)$ become polynomials in $\omega$, which are given by the partition operator acting with the elements $i$ 
assigned to the value of $\omega$ as in Eq.\ (\ref{six-thirtytwob}). As expected, these polynomials denoted by $p(k,\omega)$ 
reduce to the partition function $p(k)$ for $\omega \!=\! 1$. Their coefficients are referred to as the sub-partition numbers 
$p_i(k)$, but are also equal to $\begin{vmatrix} k \\ m\end{vmatrix}$, where the latter notation was introduced in Sec.\ 
2 to denote the number of partitions summing to $k$ with $m$ elements.

Sec.\ 8 continues with further generalisations of $P(z)$, where the coefficients of the powers of $z$ are set equal to $C_k$. 
This represents not only the inversion of Theorem\ 2, but is also a special case of the corollary to Theorem\ 3. The specific 
case of the $C_k$ equalling unity has the interesting property that the coefficients of the generating function represent the 
difference between the number of even- and odd-element partitions summing to $k$. The absolute values of these coefficients 
are shown to equal the number of discrete partitions with only odd elements, a result first obtained by Euler. Then recurrence 
relations are developed for the number of discrete partitions summing to $2k$ and $2k+1$ given by Eqs.\ (\ref{six-thirtythree15}) 
and (\ref{six-thirtythree16}), respectively. Like the Euler/MacMahon recurrence relation or Eq.\ (\ref{six-nine}), which evaluates
the values of $p(k)$, the new recurrence relations also require the discrete partition numbers. This means that not many of the
previous values of $q(k,1)$ are required to determine the highest successive value.

The approach that resulted in the alternative representations for $P(z)$ and its inverse given by Eqs.\ (\ref{six-eleven}) and 
(\ref{six-twelvea}) is then applied to $P(z,\omega)$. Replacing the $\gamma_i$ are the divisor polynomials $\gamma_i(\omega)$, 
where each coefficient represents a divisor $d$ of $i$ divided by $i$ while the power of $\omega$ equals the inverse or 
reciprocal of this value. As a result, the discrete partition and partition function polynomials can be expressed in terms of 
the partition operator acting with the elements $i$ assigned to $-\gamma_i(\omega)$ and $\gamma_i(\omega)$, respectively. Next,
the program ${\bf dispfnpoly}$ is presented. The purpose of this code is to express both $q(k,-\omega)$ and $p(k,\omega)$ in
symbolic form so that they can be imported into Mathematica, where the Divisors routine can be exploited to yield the final 
forms as polynomials in $\omega$ of degree $k$.

Sec.\ 8 concludes by applying the the preceding material to derive the generating functions for more advanced infinite products.
First, the product $P_{\rho}(z,\omega)$ is studied, in which an arbitrary power of $\rho$ is applied to $P(z,\omega)$. The 
generating function for the new product has coefficients $p(k,\omega,\rho)$, which are expressed in terms of the partition 
operator acting on elements $i$ assigned values of $p(i,\omega)$, but now the multinomial factor is altered as described in
the proof of Corollary\ 1 to Theorem\ 1. The next example is the product of $Q(z,-\beta \omega)$ as defined by Eq.\ 
(\ref{six-twentysixa}) and $P(z,\alpha \omega)$. The coefficients of the generating function become special polynomials
$QP_k(\omega,\beta,\alpha)$ as defined by Eq.\ (\ref{six-thirtyseven}) and are tabulated in Table\ \ref{table3}. In these coeffcients
the power of $\omega$ represents the number of elements in the final partitions, which are composed of the elements from
both discrete and standard partitions, while the powers of $\alpha$ and $\beta$ represent the numbers of elements emanating
from the standard and discrete partitions, respectively. The final example in the section is the derivation of the generating 
function for the famous three-parameter plus one variable product first studied by Heine and given by Eq.\ (\ref{six-thirtynine}). 
This product, which represents the product of $P(z,\omega x,\omega)$ and $P(z,\omega y,\omega x y)$, is found to possess a 
generating function whose coefficients are polynomials denoted by $HP_k(\omega,x,y)$. Like the $QP_k(\omega,\beta,\alpha)$, they
are of degree $k$ in $\omega$. General expressions for some of the coefficients are derived, while the first seven polynomials are 
presented in Table\ \ref{table4}.

\section{Acknowledgements}
The author is grateful to Dr. D.X. Balaic of Australian X-ray Capillary Optics Pty. Ltd. (http://dxb@axco.com.au),
for creating computer programs based on the BRCP algorithm that served as the basis for some of the programming 
material appearing in this work. He also thanks Dr. J. Kelleher, University of Edinburgh, for very informative
discussions, Prof. C. Krattenthaler, University of Vienna, for incisive remarks and Prof. G.E. Andrews, Pennsylvania
State University, for his interest and advice.

\newpage
\section{Appendix}
The first code presented in this appendix is the program ${\bf partmeth}$ discussed in detail at the beginning of 
Sec.\ 4. This code employs the BRCP alogrithm of Sec.\ 2 to generate the coefficients arising from the partition 
method for a power series expansion according to Theorem\ 1 in symbolic form so that the values of the coefficients 
can be evaluated either in integer form or as algebraic expressions in Mathematica. 

\begin{lstlisting}{}
/* This code deals with the application of the partition 
   method for a power series expansion to the pseudo-
   composite function g(af(z)). Here it is assumed that 
   g(z)= h(z)(1+ q_1 z + q_2 z^2+ ...+q_k z^k +...), where 
   h(z) can be any function, but is usually equal to unity 
   or some factor multiplied by a non-integer power of z. 
   In addition, f(z) is assumed to be a power series 
   expansion of the form  (p_0 + p_1 y+ p_2 y^2 + ...+ 
   p_k y^k+...) with y=z^alpha. This code is only valid
   for the case of p_0=0. The coefficients DS[k,n] and 
   ES[k,n] are those in Theorem 1 and are computed in a 
   format suitable for processing in Mathematica. */

#include <stdio.h>
#include <memory.h>
#include <stdlib.h>
#include <math.h>
#include <time.h>

int dim,sum,*part,inv_case;
time_t init_time, end_time;

/* If inv_case equals unity, then the code will evaluate 
   the coefficients of the inverse power series, i.e. 
   the power series for h(af(z))/g(af(z)). */ 

void termgen(int p)
{ 
int freq,i,num_parts=0,l,spacing=0,num_dis_parts=0,dis_parts=0;
double sign,dnum_parts;

/* num_parts is the total # of parts in the partition, 
   while num_dis_parts is the number of distinct parts 
   with greater than unity frequency snd is required 
   for the multinomial factor.    */

if(p==sum) printf((inv_case >0)?"ES[%i,n_]:= ":"DS[%i,n_]:= ",p);
if(inv_case==0  && p!=sum) printf("+"); 
for(i=1; i<=dim; i++){
          freq=part[i];
          if(freq){
               dis_parts++;
               if(inv_case==0){
                       printf("p[%i,n]",i);
                       if(freq>1) printf("^(%i) ",freq);
                       else printf(" ");
                              }
               num_dis_parts += (freq>1);  
               num_parts=num_parts+freq;
                  }
                     }
if(inv_case>0){
          dnum_parts=(double) num_parts;
          sign=pow(-1.0,dnum_parts);
          printf((sign>0.0)?"+":"-");
              }


if(inv_case>0){
           printf("DS[0,0]^(-%i) ",num_parts+1);
              }
else{
           printf("q[%i] a",num_parts);
           printf((num_parts>1)? "^(%i) " : " ",num_parts);
    }

if(inv_case>0){
          for(i=1;i<=dim;i++){
             freq=part[i];
             if(freq==1) printf("DS[%i,n] ",i);
             else if(freq>1) printf("DS[%i,n]^(%i) ",i,freq);
                             }
              }
if(num_parts>1 && dis_parts>1){
             printf("%i!",num_parts);
             if(num_dis_parts){
                    printf((num_dis_parts>1)? "/(" : "/");
                    for(i=1; i<=dim; i++){
                             freq=part[i];
                             if(freq>1){
                                  if(spacing++) printf(" ");
                                  printf("%i!",freq);
                                       }
                                         }
                    if (num_dis_parts>1) printf(")");
                              }
                    printf(" ");   
                            }
printf("\n");
}

void idx(int p,int q)
{
part[p]++;
termgen(p);
part[p]--;
p -= q;
while(p >= q){
        part[q]++;
        idx(p--,q);
        part[q++]--;
        }
}

int main(int argc, char *argv[])
{
int i;
double delta_t;

time(&init_time);
if(argc != 2) printf("usage: ./partmeth <#partitions>\n" );
else{
        dim=atoi(argv[1]);
        for(sum=1; sum<=dim; sum++){
               part=(int *) malloc((dim+1)*sizeof(int));
               inv_case=0;
               for(i=1; i<=dim; i++) part[i]=0;
               idx(sum,1);
               free(part);
               printf("\n");
               inv_case=1;
               part=(int *) malloc((dim+1)*sizeof(int));
               for(i=1;i<=dim; i++) part[i]=0;
               idx(sum,1);
               free(part);
               printf("\n");
                                   }
     }       
printf("\n");
time(&end_time);
delta_t= difftime(end_time,init_time);
printf("Computation time is %f seconds\n",delta_t);
return(0);
}
\end{lstlisting}

The above code represents the implementation of Theorem\ 1 for the case where $p_0$ vanishes. Once the order of
coefficients becomes sufficiently large, viz.\ for $k \geq 20$, the code needs to be adapted so that only specific 
values for one of the two types of coefficient are evaluated in symbolic form. This requires: (1) separating the 
inverse case or the $E_k$ from the $D_k$, and (2) removing the first for loop in ${\bf main}$ so that only $i \!=\! dim$ 
is computed. The modified code called ${\bf mathpm}$, which determines only the $D_k$ in symbolic form, is presented 
below.

\begin{lstlisting}
/* The code mathpm determines the coefficients of the
   power series for the pseudo-composite function 
   g(af(z)), where g(z)= h(z) (1+ q_1 z + q_2 z^2+ ...
   +q_k z^k +...) and  h(z) can be any function. The 
   function f(z) must be expressed as (p_0 + p_1 y+ 
   p_2 y^2 + ... + p_k y^k+...) where y=z^alpha and 
   p_0 =0.   */

#include <stdio.h>
#include <memory.h>
#include <stdlib.h>
#include <time.h>

int dim,*part;
long unsigned int term=1;
time_t init_time, end_time;

void termgen(int p)
{
int f,i,num_parts=0,dis_part_cnt=0,l,num_dis_parts=0;
/* num_parts is the total # of parts in the partition, 
   while num_dis_parts is the number of distinct parts. 
   The latter is required for the multinomial factor. */

if(p==dim) printf("DS[%i,n_]:= ",p);
else{
     printf("+ ");
     if (term %3 == 0) printf("\n");
     term++;
     }
for(i=1; i<=dim; i++){
        f=part[i];
        if(f){
                printf("p[%i]",i);
                if(f>1) printf("^(%i) ", f);
                else printf(" ");
                num_parts += l = f;
                num_dis_parts += (f>1);
             }
                     }
   printf("q[%i] a",num_parts);
   printf((num_parts>1)? "^(%i) " : " ",num_parts);
if( num_parts > l ){
        printf("%i!",num_parts);
        if(num_dis_parts){
                printf((num_dis_parts>1)? "/(" : "/");
                for(i=1; i<=dim; i++){
                      f=part[i];
                      if(f>1){
                           if(dis_part_cnt++) printf(" ");
/* dis_part_cnt counts the # of distinct parts with 
   greater than unity frequency It is required to 
   insert blank spaces in the denominator of the 
   multinomial factor  */
                           printf("%i!", f);
                             }
                                     }
                if(num_dis_parts>1) printf(")");
                          }
        printf(" ");
                   }
}

void idx(int p,int q)
{
part[p]++;
termgen(p);
part[p]--;
p -= q;
while(p >= q){
        part[q]++;
        idx(p--,q);
        part[q++]--;
             }
}

int main( int argc, char *argv[])
{
int i;
double delta_t;
FILE *ptr;
char filename[10]="times";

time(&init_time);
if(argc != 2) printf("usage: ./mathpm <#partitions>\n" );
else{
    dim=atoi(argv[1]);
    part=(int *) malloc((dim+1)*sizeof(int));
    if(part==NULL) printf("unable to allocate array\n\n");
    else{
        free(part);
        }
    }
printf("\n");
time(&end_time);
delta_t= difftime(end_time,init_time);
ptr=fopen(filename,"a");
fprintf(ptr,"Time to compute p[%i,n] is %f seconds\n",
	dim,delta_t);
fclose(ptr);
\end{lstlisting} 

In Sec.\ 5 we discussed the problem of generating discrete partitions or partitions in which the elements
only appear once. Since discrete partitions are important in the theory of partitions as can be seen from
Secs. 6-8, the entire code called ${\bf dispart.cpp}$ is presented below.  

\begin{lstlisting}{}
#include <stdio.h>
#include <memory.h>
#include <stdlib.h>
#include <time.h>

int tot, *part;
long unsigned int term=1;
time_t init_time, end_time;

void termgen()
{
     int freq,i;

     for(i=0;i<tot;i++){
                freq=part[i];
                if(freq>1) goto end;
                       }
     printf("%ld: ",term++);
     for(i=0;i<tot;i++){
                freq=part[i];
                if(freq) printf("%i(%i) ",freq,i+1);
                       }
     printf("\n");
end:;
}

void idx(int p,int q )
{
part[p-1]++;
termgen();
part[p-1]--;
p -= q;
while(p >= q){
        part[q-1]++;
        idx(p--,q);
        part[q++ - 1]--;
             }
}

int main( int argc, char *argv[] )
{
int i;
double delta_t;
FILE *ptr;
char filename[10]="times";

time(&init_time);
if(argc != 2) printf("usage: ./dispart <#partitions>\n");
else{
  tot=atoi( argv[1]);
  part= (int * ) malloc(tot*sizeof(int));
  if(part == NULL) printf("unable to allocate array\n\n");
  else{
      for(i=0;i<tot;i++) part[i]=0;
      idx(tot,1);
      free(part);
      }
     }
printf("\n");
time(&end_time);
delta_t=difftime(end_time,init_time);
ptr=fopen(filename,"a");
fprintf(ptr,"Time taken to compute discrete partitions 
	summing to %i is %f secs.\n",tot,delta_t);
fclose(ptr);
return(0);


\end{lstlisting} 

The next code presented here determines the conjugate partition as described in Sec.\ 5. In evaluating 
conjugate partitions, the programme called ${\bf transp}$ casts the original partition in the form of a 
Ferrers diagram, but, instead of being composed of dots, the Ferrers diagram is composed of ones. Thus, 
the conjugate partition is determined by summing the ones in each column of the Ferrers diagram.

\begin{lstlisting}{}
/* This program evaluates the partitions and their 
conjugates for any integer greater than or equal to 
2. Conjugates are determined by summing the columns 
in the Ferrers diagram for each partition. */

#include <stdio.h>
#include <memory.h>
#include <stdlib.h>

int tot, *part;
long unsigned int term=1;

void termgen()
{
int freq,i,j,k,next_el,index,prev_el,*part2,
        *sum_col,row_cnt=0;
int *ferrers,**rptr; 
/* The Ferrers array and its array of pointers */

printf("Partition %ld is: ",term++);
ferrers=(int *) malloc(tot*tot*sizeof(int));
rptr=(int **) malloc(tot*sizeof(int *));

/* Get the pointers to the rows of ferrers */
for (i=0;i<tot;i++) rptr[i]= ferrers+(i*tot);
/* Here tot refers to number of columns */

for(j=0;j<tot;j++){
      for(i=0;i<tot;i++){
             rptr[j][i]=0;
                        }
                  }
/* Creation of the Ferrers diagram. Rather than 
   being composed of dots the Ferrers diagram is 
   composed of unit values since these will be 
   used to determine the conjugate partition. */

for(j=0;j<tot;j++){
     freq=part[j];
     if(freq){
           for(i=row_cnt;i<row_cnt+freq;i++){
                  for(k=0;k<=j;k++){
                       rptr[i][k]=1;
                                   }
                                            }
             }
     row_cnt=row_cnt+freq;
                  }
/* Summation of the columns in the Ferrers 
   diagram yielding a new array called sum_col. */

sum_col=(int *) malloc(tot*sizeof(int));
/* Initialising the array elements to zero */

for(i=0;i<tot;i++) sum_col[i]=0;
/* Summation of the columns now occurs */

for (j=0;j<tot;j++){
     for(i=0;i<tot;i++) sum_col[j]=sum_col[j]+rptr[i][j];
                   }

/* The array sum_col is reduced to the conjugate 
   of the original partition through another array 
   called part2 */

part2=(int *) malloc(tot*sizeof(int));
/* Initialisation of the array elements to zero */

for(i=0;i<tot;i++) part2[i]=0;

/* Now the conjugate partition is arranged as 
in the order of part */

prev_el=sum_col[0]; 
/* the highest part is the value of sum_col[0] */
index=prev_el-1; 
/* the index in the partition array is one less */
part2[index]=1;
for(i=1;i<=tot;i++){
     next_el=sum_col[i];
     if(next_el==0) goto out;
     if(next_el==prev_el) part2[index]=part2[index]+1;
     else{
          index=next_el-1;
          part2[index]=1;
          prev_el=next_el;
         }
                   }
out: for(i=0;i<tot;i++){
     freq=part[i];
     if(freq) printf("%i(%i) ",freq,i+1);
                       }
/* The conjugate is actually the reverse order 
   of part2 */
printf(" and its conjugate is: ");
for(i=0;i<tot;i++){
     freq=part2[tot-i-1];
     if(freq) printf("%i(%i) ",freq,tot-i);
                  }
free(sum_col);
free(part2);
free(ferrers);
free(rptr);
printf("\n");
}

void idx(int p,int q)
{
part[p-1]++;
termgen();
part[p-1]--;
p -= q;
while(p >= q){
        part[q-1]++;
        idx(p--,q);
        part[q++ -1]--;
             }
}

int main(int argc,char *argv[])
{
int i;
if(argc != 2) printf("usage:./transp <#partitions>\n");
else{
        tot=atoi(argv[1]);
        part=(int *) malloc(tot*sizeof(int));
        if(part == NULL) printf("unable to allocate 
				array\n\n");
        else{
                for(i=0;i<tot;i++) part[i]=0;
                idx(tot,1);
                free(part);
            }
    }
printf("\n");
return(0);
}
\end{lstlisting}  

In Sec.\ 6 the partition method for a power series expansion was applied to an exponential form of the generating function of 
the partition function $p(k)$, which yielded Eq.\ (\ref{six-fifteen}). Although this result represents a sum over partitions 
involving the special numbers $q(i)$ called the discrete partition numbers, it incorporates much redundancy because these
numbers are often zero except when $i$ can be written as a pentagonal number or as $(3j^2\pm j)/2$, where $j$ is a non-negative 
integer. Then they are equal to $(-1)^j$. Appearing below is the program ${\bf pfn}$ which gives the partition function based 
on the properties of the discrete partition numbers (represented by {\it Q}[i] in the code) in symbolic form. The actual values 
of the partition function can be evaluated by importing the output into Mathematica. 

\begin{lstlisting}
#include <stdio.h>
#include <memory.h>
#include <stdlib.h>
#include <time.h>
#include <math.h>

int dim,*part,limit,freq,first_term=0;
long unsigned int term=1;
time_t init_time, end_time;

void termgen(int p)
{
int f,i,num_parts=0,j,jval,dis_part_cnt=0,l,
	num_dis_parts=0;
/* num_parts is the total # of parts in the partition 
   while num_dis_parts is the number of distinct parts. 
   The latter is required for the multinomial factor. */
/* jval=0; */

if(p==dim){
      printf("p[%i]:= ",p);
      for (i=1;i<=dim;i++){
                   jval=0;
                   freq=part[i];
                   if (freq >0){
                         for (j=1;j<=limit;j++){
                             if(i == (3*j-1)*j/2) jval=j;
                             if(i == (3*j+1)*j/2) jval=j;
                                               }
                         if((jval==0) && (freq>0)) goto end;
                             first_term++;
                               }
                          }
          }
else {
     for (i=1;i<=dim;i++){
                   jval=0;
                   freq=part[i];
                   if (freq >0){
                       for (j=1;j<=limit;j++){
                           if(i == (3*j-1)*j/2) jval=j;
                           if(i == (3*j+1)*j/2) jval=j;
                                             }
                       if((jval==0) && (freq>0)) goto end;
                               }
                         }
     if(first_term !=0) printf("+ ");
     if(first_term ==0) first_term++;

     if (term %3 == 0) printf("\n");
     term++;
     }


for(i=1; i<=dim; i++){
        f=part[i];
        if(f ) {
/*                printf("Q[%i]",i);   */
             for (j=1;j<=limit;j++){
                        if (i == (3*j-1)*j/2) jval=j;
                        if (i == (3*j+1)*j/2) jval=j;
                                   }
             printf("((-1)^(%i))",jval);
             if(f>1) printf("^(%i) ", f);
             else printf(" ");
             num_parts += l = f;
             num_dis_parts += (f>1);
               }
                    }
printf("(-1)");
printf((num_parts>1)? "^(%i) " : " ",num_parts);
if( num_parts > l ){
        printf("%i!",num_parts);
        if(num_dis_parts){
              printf((num_dis_parts>1)? "/(" : "/");
              for(i=1; i<=dim; i++){
                    f=part[i];
                    if(f>1){
                        if(dis_part_cnt++) printf(" ");
/* dis_part_cnt counts the # of distinct parts with greater
   than unity frequency. It is required for inserting blank 
   spaces in the denominator of the multinomial factor */
                        printf("%i!", f);
                           }
                                   }
              if (num_dis_parts>1) printf(")");
                         }
        printf(" ");
                  }
end: ;
}
void idx(int p,int q)
{
part[p]++;
termgen(p);
part[p]--;
p -= q;
while(p >= q){
        part[q]++;
        idx(p--,q);
        part[q++]--;
        }
}

int main( int argc, char *argv[] )
{
int i;
double delta_t;

time(&init_time);
if(argc != 2) printf("usage: ./pfn <#partitions>\n" );
else{
        dim=atoi(argv[1]);
        limit= floor(1+sqrt((1+ 24 * dim))/6);

        part=(int *) malloc((dim+1)*sizeof(int));
        if(part==NULL)printf("unable to allocate array\n\n");
        else{
/*                for(sum=1; sum <=dim; sum++){    */
                        for(i=1; i<=dim; i++) part[i] = 0;
                        idx(dim,1);
/*                        printf("\n");
                        }    */
                free(part);
                }
        }
printf("\n");
time(&end_time);
delta_t= difftime(end_time,init_time);
printf("Time taken to compute the coefficient is 
	%f seconds\n",delta_t);
return(0);
}

\end{lstlisting}

Also in Sec.\ 6, the discrete partition polynomials $q(k,\omega)$ were found to be the coefficients of
the generating function for by the product $Q(z,\omega)=\prod_{k=1}^{\infty} 1/(1+\omega z^k)$ (see 
Eq.\ (\ref{six-twentyfive})). Later, it was found that these polynomials of degree $i$ in $\omega$ could be 
expressed in terms of the partition operator acting with the elements $i$ equal to special polynomials $\gamma_i(\omega)$. 
The latter were referred to as divisor polynomials since their coefficients are divisors or factors of $i$. The 
relationship between both types of polynomials is given by Eq.\ (\ref{six-thirtyfourb1}), while another result 
below it relates the partition function polynomials $p(k,\omega)$ to another sum over partitions acting on
the divisor polynomials. Appearing below is the program called ${\bf dispfnpoly}$, which prints out both
the discrete partition and partition function polynomials for a specified order in symbolic form so that they can
be handled by Mathematica \cite{wol92}. 
\begin{lstlisting}
#include <stdio.h>
#include <memory.h>
#include <stdlib.h>
#include <time.h>

int dim,*part,polytype;
long unsigned int term=1;
time_t init_time, end_time;

void termgen(int p)
{
int f,i,num_parts=0,dis_part_cnt=0,l,num_dis_parts=0;
/* num_parts is the total # of parts in the partition 
   while num_dis_parts is the number of distinct parts. 
   The latter is required for the multinomial factor */

if(p==dim) printf((polytype==1)?"Q[%i,-w_]:= ":
	"P[%i,w_]:= ",p);

for(i=1; i<=dim; i++){
        f=part[i];
        if(f){
                num_parts += l = f;
                num_dis_parts += (f>1);
                }
        }
    if(term !=1) printf(" + ");
    if(term % 3 ==0) printf("\n");
    term++;
    for(i=1; i<=dim; i++){
        f=part[i];
        if(f){
                printf("DP[%i,w]",i);
                if(f>1) printf("^(%i) ",f);
                else printf(" ");
                }
        }
     printf((polytype==1)?"(-1)":"");
     if (num_parts>1) printf((polytype==1)?"^(%i)":"",
             num_parts);
     if(num_parts > 1 ){
          if(num_dis_parts){
                printf((num_dis_parts>1)? "/(" : "/");
                for(i=1; i<=dim; i++){
                        f=part[i];
                        if(f>1){
                             if(dis_part_cnt++) printf(" ");
/* dis_part_cnt counts the # of distinct parts with greater
   than unity frequency. It is required for inserting blank 
   spaces in the denominator of the multinomial factor */
                             printf("%i!",f);
                               }
                        }
                if(num_dis_parts>1) printf(")");
                }
        }
}

void idx(int p,int q)
{
part[p]++;
termgen(p);
part[p]--;
p -= q;
while(p >= q){
        part[q]++;
        idx(p--,q);
        part[q++]--;
        }
}

int main(int argc, char *argv[])
{
int i;
double delta_t;

time(&init_time);
if(argc != 2) printf("usage: ./dispfnpoly 
		<#partitions>\n");
else{
      dim=atoi(argv[1]);
      part=(int *) malloc((dim+1)*sizeof(int));
      if(part==NULL) printf("unable to allocate 
			array\n\n");
      else{
                for(polytype=1; polytype<=2;polytype++){
                      for(i=1; i<=dim; i++) part[i] = 0;
                      idx(dim,1);
                      if(polytype==1) printf("\n\n");
                      term=1;
                                                       }
                free(part);
          }
     }
printf("\n");
time(&end_time);
delta_t= difftime(end_time,init_time);
printf("Time taken to compute the coefficients is 
	%f seconds\n",delta_t);
return(0);
}





\end{lstlisting}

\end{document}